\documentclass{article}
\usepackage{amsmath,amssymb,amscd}
\usepackage{tikz}
\usepackage{tikz-cd}
\newtheorem{Thm}{Theorem}
\newtheorem{Def}[Thm]{Definition}
\newtheorem{Lem}[Thm]{Lemma}
\newtheorem{Rem}[Thm]{Remark}
\newtheorem{Prop}[Thm]{Proposition}
\newtheorem{Cor}[Thm]{Corollary}

\newtheorem{Conjecture}[Thm]{Conjecture}
\newtheorem{Claim}[Thm]{Claim}

\begin{document}
\title{On the Hofer-Zehnder conjecture for non-contractible periodic orbits in Hamiltonian dynamics}
\author{Yoshihiro Sugimoto}
\date{}
\maketitle

\begin{abstract}
In this paper, we treat non-contractible periodic orbits in Hamiltonian dynamics on symplectic manifolds. We prove that any Hamiltonian diffeomorphism has infinitely many simple non-contractible periodic orbits provided that the Hamiltonian diffeomorphism has at least one periodic orbit of infinite order in the first homology group and the orbit has non-trivial local Floer cohomology. Our proof is an application of the equivariant Hamiltonian Floer cohomology.
\end{abstract}

\section{Introduction and main results}
In this section, we briefly explain the main theme of this paper. The precise definitions and notations are given in the next section. In this paper, we treat periodic orbits in Hamiltonian dynamics. The Hofer-Zehnder conjecture is a conjecture about the number of contractible periodic orbits in Hamiltonian dynamics and it states that ``every Hamiltonian map on a compact symplectic manifold $(M,\omega)$ possessing more fixed points than necessarily required by the V. Arnold conjecture possesses always infinitely many periodic points"(\cite{HZ} p. 263). This is a generalization of Franks' theorem \cite{Fr,Fr2}, which states that any area-preserving homeomorphism on the sphere ${S^2}$ with more than two fixed points has infinitely many periodic points. We have to clarify the meaning of ``more fixed points than necessarily required by the V. Arnold conjecture". The Floer cohomology group of a Hamiltonian function is generated by its contractible periodic orbits and it is isomorphic to the singular homology of the underlying manifold. So, the number of contractible periodic orbits required by the Arnold conjecture is the sum of the Betti numbers. In summary, the Hofer-Zehnder conjecture can be stated in the following form:

\begin{Conjecture}[Hofer-Zehnder conjecture] \label{ConjectureB}
Let ${\phi\in \textrm{Ham}(M,\omega)}$ be a Hamiltonian diffeomorphism with more simple contractible periodic orbits than the total Betti number of $M$. Then, ${\phi}$ has infinitely many simple contractible periodic orbits.
\end{Conjecture}

Shelukhin proved this conjecture for spherically monotone symplectic manifolds with semi-simple quantum cohomology ring under the assumption that periodic orbits are counted homologically (\cite{Sh}). Atallah and Lou generalized Shelukhin's theorem to weakly monotone symplectic manifolds (\cite{AL}). Bai and Xu proved the Hofer-Zehnder conjecture for all symplectic toric manifolds (without assuming that the quantum cohomology ring is semi-simple) under the assumption that periodic orbits are counted homologically (\cite{BX}). Readers may consider the relationship between the Conley conjecture and the Hofer-Zehnder conjecture. We say that the Conley conjecture holds on a closed symplectic manifold ${(M,\omega)}$ if every Hamiltonian diffeomorphism on ${(M,\omega)}$ has infinitely many simple contractible periodic orbits. Roughly speaking, the Hofer-Zehnder conjecture states that any Hamiltonian diffeomorphism with finitely many contractible periodic orbits is a so-called pseudo-rotation. A pseudo-rotation is a Hamiltonian diffeomorphism with the minimum number of contractible periodic orbits. The Conley conjecture and the Hofer-Zehnder conjecture imply that it is important to know the necessary and sufficient conditions for the existence of a pseudo-rotation on a symplectic manifold. The Conley conjecture implies that symplectic manifolds with a pseudo-rotation are very rare. In \cite{CGG,CGG2,GG6,GG7,Sh2,Sh3}, the importance of the existence of non-trivial pseudo-holomorphic curves was pointed out. They proved that the quantum Steenrod square is deformed if the symplectic manifold is monotone and has a pseudo-rotation. This is strong evidence for the Chance-McDuff conjecture which states that the Conley conjecture holds if some Gromov-Witten invariants vanish or if the quantum cohomology ring is undeformed (\cite{GG3,GG4}). A slightly different probrem from the Conley conjecture is the generic Conley conjecture, which states that $C^{\infty}$-generic Hamiltonian diffeomorphisms have infinitely many simple contractible periodic orbits. The generic Conley conjecture is proved for a broad class of symplectic manifolds (\cite{GG9,Sug}). It is also proved that $C^{\infty}$-generic Hamiltonian diffeomorphisms have infinitely many hyperbolic periodic points for a broad class of symplectic manifolds (\cite{CGG4}).
 
 We can state an analogue of the Hofer-Zehnder conjecture for non-contractible periodic orbits. The set of non-contractible periodic orbits may be empty and the Floer cohomology of non-contractible periodic orbits is always trivial. This implies that the necessary number of non-contractible periodic orbits is zero. So, one is tempted to conjecture that the existence of one non-contractible periodic orbit implies the existence of infinitely many simple non-contractible periodic orbits. G\"{u}rel, Ginzburg-G\"{u}rel and Orita  proved that any Hamiltonian diffeomorphism has infinitely many simple non-contractible periodic orbits provided that the Hamiltonian diffeomorphism has at least one periodic orbit of infinite order in the first homology group for some closed simplectic manifolds \cite{G,GG8,O,O2}. In this paper, we apply equivariant Floer theory \cite{Se,Sh,SZ} and a key idea of G\"{u}rel \cite{G}. We prove this conjecture for very wide classes of symplectic manifolds.

 Our main result is stated as follows. A ${2n}$-dimensional symplectic manifold ${(M,\omega)}$ is called a weakly monotone symplectic manifold if it satisfies one of the following conditions. We explain the precise meaning of terminologies in the next section.

\begin{enumerate}
\item ${(M,\omega)}$ is monotone 
\item ${c_1(A)=0}$ for every ${A\in \pi_2(M)}$
\item The minimal Chern number ${N>0}$ is greater than or equal to ${n-2}$
\end{enumerate}
Note that weakly monotone symplectic manifolds cover wide classes of symplectic manifolds. For example, every symplectic manifold whose dimension is less than or equal to $6$ is a weakly monotone symplectic manifold. We prove the following theorem.

\begin{Thm}
Let ${(M,\omega)}$ be a closed weakly monotone symplectic manifold, and let ${\phi \in \textrm{Ham}(M,\omega)}$ be a Hamiltonian diffeomorphism generated by a periodic Hamiltonian function ${H\in C^{\infty}(S^1\times M)}$. Suppose that the number of ${1}$-periodic orbits of $\phi$ in a non-trivial class ${\gamma \neq 0 \in H_1(M:\mathbb{Z})/\textrm{Tor}}$ is finite and non-zero. Further suppose that the local Floer cohomology ${HF^{loc}(\phi, x)}$ of at least one of these orbits is non-trivial.
Then, for every sufficiently large prime $p$, ${\phi}$ must have at least one $p$-periodic or $p'$-periodic simple orbit in the class ${p\cdot \gamma}$. Here, $p'$ is the smallest prime greater than $p$. In particular, there are infinitely many simple periodic orbits of $\phi$ with classes contained in the set of classes ${\mathbb{N}\cdot \gamma}$.
\end{Thm}

The assumption ``weakly monotone" is purely technical. Our proof is based on the ${\mathbb{Z}_p}$-equivariant Floer cohomology theory and we need ${\mathbb{Z}_p}$-coefficient Floer theory. For general closed symplectic manifolds, we need the so-called virtual technique to define Floer theory (\cite{FO2,LT}), but the virtual technique works over $\mathbb{Q}$-coefficient in general (see also \cite{FO}). We are not sure we can construct an equivariant theory on general closed symplectic manifolds.

\section{Preliminaries}
In this section, we explain notations and terminologies used in this paper.

\subsection{Elementary notations}
Let ${(M,\omega)}$ be a symplectic manifold, so $M$ is a finite-dimensional $C^{\infty}$-manifold and ${\omega \in \Omega^2(M)}$ is a symplectic form on $M$. In this paper, we always assume that $M$ is a closed manifold.

For any $C^{\infty}$-function $H\in C^{\infty}(M)$, we define the Hamiltonian vector field ${X_H}$ by the following relation:

\begin{equation*}
\omega(X_H, \cdot)=-dH.
\end{equation*}

We can also consider an ${S^1}$-dependent (in other words, $1$-periodic) Hamiltonian function $H$ and a Hamiltonian vector field $X_H$ by the same formula. The time-one map of the flow of $X_H$ is called a Hamiltonian diffeomorphism generated by $H$ and is denoted by $\phi_H$. The set of all Hamiltonian diffeomorphisms is called the Hamiltonian diffeomorphism group and we denote the Hamiltonian diffeomorphism group of ${(M,\omega)}$ by ${\textrm{Ham}(M,\omega)}$. That is,

\begin{equation*}
\textrm{Ham}(M,\omega)=\{\phi_H \ | \ H\in C^{\infty}(S^1\times M)\}.
\end{equation*}
We also consider ``iterations" of $H$ and ${\phi_H}$. For any integer ${k\in \mathbb{N}}$, we define ${H^{(k)}}$ as follows:

\begin{equation*}
H^{(k)}=kH(kt, x).
\end{equation*}
It is straightforward to see that ${\phi_{H^{(k)}}=(\phi_H)^k}$ holds. Let ${P^l(H)}$ be the space of $l$-periodic orbits of $X_H$. Set ${S^1_l=\mathbb{R}/l\cdot \mathbb{Z}}$. Then,

\begin{gather*}
P^l(H)=\{x:S_l^1\rightarrow M \ | \ \dot{x}(t)=X_{H_t}(x(t))  \}
\end{gather*}
holds. It is also straightforward to see that there is a one-to-one correspondence between ${P^k(H)}$ and ${P^1(H^{(k)})}$. We abbreviate ${P^1(H)}$ to ${P(H)}$. An $l$-periodic orbit ${x\in P^l(H)}$ is called simple if there is no $l'$-periodic orbit ${y\in P^{l'}(H)}$ which satisfies the following conditions:

\begin{gather*}
l=l'\cdot m  \ \ \ \ (l', m\in \mathbb{N}) \\ 
x(t)=y(\pi_{l,l'}(t)).
\end{gather*} 
Here, ${\pi_{l,l'}:S_l\rightarrow S_{l'}}$ is the natural projection. So, a periodic orbit is simple if and only if it is not an iteration of a periodic orbit of a lower period.

Next, we explain the definition of the minimal Chern number $N$. ${(M,\omega)}$ becomes an almost complex manifold after a choice of a compatible almost complex structure ${J}$ on the tangent bundle ${TM}$. This gives the tangent bundle the structure of a complex vector bundle, hence allowing us to define its first Chern class ${c_1(TM)\in H^2(M:\mathbb{Z})}$. Moreover, the first Chern class does not depend on the choice of $J$ because the space of compatible almost complex structures on $M$ is contractible. The minimal Chern number $N\in \mathbb{N}\cup \{+\infty\}$ is the positive generator of the image of ${c_1(TM)|_{\pi_2(M)}}$. Note that if the image is zero, ${N}$ is defined by ${N=+\infty}$. A symplectic manifold ${(M,\omega)}$ is called monotone if the cohomology class of the symplectic form ${[\omega]}$ over ${\pi_2(M)}$ is a non-negative multiple of the first Chern class. In other words, there is a constant ${\lambda \ge 0}$ such that 
\begin{equation*}
\int_{S^2}v^*\omega=\lambda\int_{S^2}v^*c_1
\end{equation*}  
holds for any smooth ${v:S^2\rightarrow M}$. As we mentioned in the previous section, a weakly monotone symplectic manifold is defined as follows:

\begin{Def}
A ${2n}$-dimensional symplectic manifold ${(M,\omega)}$ is called weakly monotone if and only if it satisfies one of the following conditions:
\begin{enumerate}
\item $(M,\omega)$ is a monotone symplectic manifold.
\item $c_1(A)=0$ holds for all ${A\in \pi_2(M)}$.
\item The minimal Chern number $N$ is greater than or equal to ${n-2}.$
\end{enumerate}
\end{Def}

\subsection{Floer cohomology theory}
 
 In this subsection, we explain Floer cohomology for non-contractible periodic orbits. References for this section are \cite{F1,F2,F3,FH,FHS,MS,FO2,SZ,AD}. Essentially, there is nothing new in the non-contractible case, but in this paper, we need a non-contractible Floer cochain complex over the universal Novikov ring. Let ${\mathbb{K}}$ be the ground field (In this paper, we consider the case of ${\mathbb{K}=\mathbb{F}_p}$ where ${\mathbb{F}_p}$ is a prime field of characteristic $p$). We also assume that ${(M,\omega)}$ is a weakly monotone symplectic manifold. The universal Novikov ring ${\Lambda}$ is defined as follows:

\begin{equation*}
\Lambda=\bigg\{\sum_{i=1}^{\infty}a_i\cdot T^{\lambda_i} \ \bigg| \ a_i\in \mathbb{K},\lambda_i\in \mathbb{R}, \lambda_i\to +\infty  \bigg\}.
\end{equation*}
We need the subring ${\Lambda_0\subset \Lambda}$ defined as follows:

\begin{equation*}
\Lambda_0=\bigg\{\sum a_i\cdot T^{\lambda_i}\in \Lambda \ \bigg| \ \lambda_i\ge 0 \bigg\}.
\end{equation*}
We need non-contractible Floer theory over ${\Lambda}$ and ${\Lambda_0}$. We fix a non-zero element ${\gamma\neq0\in H_1(M:\mathbb{Z})/\textrm{Tor}}$ and we denote the set of $1$-periodic orbits of ${H\in C^{\infty}(S^1\times M)}$ in ${\gamma}$ by ${P(H,\gamma)}$.

\begin{equation*}
P(H,\gamma)=\{x\in P(H) \ | \ [x]=\gamma  \}
\end{equation*}
A $1$-periodic orbit $x$ is called non-degenerate if the differential map
\begin{gather*}
d\phi_H|_{x(0)}:T_{x(0)}M \longrightarrow T_{x(0)}M
\end{gather*}
does not have eigenvalue $1$. The Floer cochain complex over ${\Lambda}$ and ${\Lambda_0}$ is defined as follows. We assume that every periodic orbit in ${P(H,\gamma)}$ is non-degenerate and set:

\begin{gather*}
CF(H,\gamma:\Lambda)=\bigoplus_{x\in P(H,\gamma)}\Lambda \cdot x  \\
CF(H,\gamma:\Lambda_0)=\bigoplus_{x\in P(H,\gamma)}\Lambda_0 \cdot x.
\end{gather*}

Note that above ${CF(H,\gamma:\Lambda)}$ and ${CF(H,\gamma:\Lambda_0)}$ are not graded over ${\mathbb{Z}}$ (It is possible to give a grading over ${\mathbb{Z}_2}$.). We fix a compatible almost complex structure ${J_t}$ parametrized by ${t\in S^1}$. The coboundary operator ${d_F}$ is defined as follows:

\begin{equation*}
d_F(x)=\sum_{y\in P(H,\gamma),\lambda\ge 0}n_{\lambda}(x,y)T^{\lambda}\cdot y,
\end{equation*}
where ${n_{\lambda}(x,y)\in \mathbb{K}}$ is the number of the solutions of the following Floer equation modulo the natural ${\mathbb{R}}$-action:

\begin{gather*}
u:\mathbb{R}\times S^1\longrightarrow M \\
\partial_su(s,t)+J_t(u(s,t))(\partial_tu(s,t)-X_{H_t}(u(s,t)))=0  \\
\lim_{s\to -\infty}u(s,t)=x(t), \lim_{s\to +\infty}u(s,t)=y(t)   \\
\lambda=\int_{\mathbb{R}\times S^1}u^*\omega+\int_0^1H(t,x(t))-H(t,y(t))dt.
\end{gather*}

\begin{Rem}
We have to achieve transversality of the linearized operator of this equation and determine an orientation of the moduli space to count the number of the solutions. They are very complex issues, but they are standard today. See the above references.
\end{Rem}

The Floer cohomology ${HF(H,\gamma:\Lambda)}$ is defined to be the cohomology of ${(CF(H,\gamma:\Lambda),d_F)}$. Note that ${\lambda>0}$ holds for all non trivial solution of the above Floer equation. So, the coboundary operator $d_F$ is defined over ${\Lambda_0}$. ${HF(H,\gamma:\Lambda_0)}$ is given by the cohomology of ${(CF(H,\gamma:\Lambda_0),d_F)}$. The Floer cohomology over ${\Lambda}$ does not depend on the choice of Hamiltonian function ${H}$ because we have the continuation map
\begin{gather*}
HF(H_1,\gamma:\Lambda)\longrightarrow HF(H_2,\gamma:\Lambda)
\end{gather*}
for any Hamiltonian functions ${H_1}$ and ${H_2}$ and this is an isomorphism. This is a very standard fact in the Floer cohomology theory (see for instance \cite{HS}). This implies that ${HF(H,\gamma:\Lambda)}$ equals to zero because ${P(f,\gamma)}$ is empty for a $C^1$-small  function $f$ and hence ${HF(f,\gamma:\Lambda)}$ is trivial. However, ${HF(H,\gamma:\Lambda_0)}$ does not equal to zero in general. In fact, there is a sequence ${0<\beta_1\le \cdots \le \beta_k}$ such that 

\begin{equation*}
HF(H,\gamma:\Lambda_0)\cong \bigoplus_{i=1}^k \Lambda_0/T^{\beta_i}\cdot\Lambda_0
\end{equation*}
holds \cite{UZ} (see also \cite{FOOO}). We have the following filtration on ${(CF(H,\gamma:\Lambda),d_F)}$. For any ${c\in \mathbb{R}}$, we have the following subcomplex of ${(CF(H,\gamma:\Lambda),d_F)}$:

\begin{equation*}
CF^c(H,\gamma:\Lambda)=\bigg\{ \sum_{x\in P(H:\gamma)}\Big(\sum_{\lambda_i\ge c}a_i\cdot T^{\lambda_i}\Big)\cdot x \in CF(H,\gamma:\Lambda) \bigg\}.
\end{equation*}
Then for any ${c<d}$, we define the following Floer complex with action window ${[c,d)}$:
\begin{equation*}
CF^{[c,d)}(H,\gamma:\Lambda)=CF^d(H,\gamma:\Lambda)/CF^c(H,\gamma:\Lambda).
\end{equation*}
We define the Floer cohomology with action window ${[c,d)}$ as follows:

\begin{equation*}
HF^{[c,d)}(H,\gamma:\Lambda)=H(CF^{[c,d)}(H,\gamma:\Lambda),d_F)).
\end{equation*}

\subsection{Local Floer cohomology}

First, we prove the following two lemmas to treat possibly degenerate Hamiltonian functions. We apply Lemma 6 to define local Floer cohomology. Lemma 5 enable us to apply the local Floer cohomology to global problems in later chapters.

\begin{Lem}
We fix ${\gamma \in H_1(M:\mathbb{Z})/\mathrm{Tor}}$. Assume that $H$ is a Hamiltonian function such that the number of $1$-periodic orbits in ${\gamma}$ is finite. Let ${P(H,\gamma)=\{y_1,\cdots,y_m\}}$ be the set of $1$-periodic orbits in ${\gamma}$ and ${U_i\subset S^1\times M}$ be a sufficiently small open neighborhood of ${\{(t,y_i(t))\}}$.
\begin{enumerate}
\item We denote the energy of solutions of the Floer equation 
\begin{gather*}
u:\mathbb{R}\times S^1\rightarrow M \\
\partial_su(s,t)+J_t(\partial_tu(s,t)-X_{H_t}(u(s,t)))=0
\end{gather*}
by $E(u)$. $E(u)$ is defines as follows:
\begin{gather*}
E(u)=\int_{\mathbb{R}\times S^1}|\partial_su(s,t)|^2dsdt=
\int_{\mathbb{R}\times S^1}\omega(\partial_su(s,t),J_t(\partial_su(s,t)))dsdt.
\end{gather*}
Let $\delta_1$ be the infimum of the energy of solutions of the Floer equation:
\begin{gather*}
\delta_1=\inf\Big\{E(u)>0 \ \Big| \ \begin{matrix}u:\mathbb{R}\times S^1\rightarrow M, \ [u(0,t)]=\gamma \\ \partial_su(s,t)+J_t(\partial_tu(s,t)-X_{H_t}(u(s,t)))=0\end{matrix}\Big\}.
\end{gather*}
Then, ${\delta_1>0}$ holds.
\item Let ${\delta_2}$ be the infimum of the energy of pseudoholomorphic spheres as follows:
\begin{gather*}
\delta_2=\inf \Big\{\int_{\mathbb{CP}^1}v^*\omega>0 \ \Big| \ \begin{matrix}v:\mathbb{CP}^1\rightarrow M \\ J_t\circ dv=dv\circ j_{\mathbb{CP}^1} \end{matrix}\Big\}.
\end{gather*}
Note that ${\delta_2>0}$ follows from the Gromov compactness theorem of pseudoholomorphic maps (see \cite{H}). Let ${\{H_i\}}$ be a sequence of Hamiltonian functions and let ${\{u_i\}\subset C^{\infty}(\mathbb{R}\times S^1,M)}$ be a sequence of solutions of the Floer equation such that
\begin{gather*}
H_i(t,x)\longrightarrow H(t,x) \ \ (\mathrm{in} \ C^{\infty}\mathrm{ \ topology})\\
\partial_su_i(s,t)+J_t(\partial_tu_i(s,t)-X_{{(H_i)}_t}(u_i(s,t)))=0 \\
[u_i(0,t)]=\gamma  \\
E(u_i)\le \alpha<\delta_2
\end{gather*}
holds for some ${\alpha<\delta_2}$. Assume that there is a sequence ${\{(s_i,t_i)\}}$ such that
\begin{gather*}
u_i(s_i,t_i)\notin \bigcup_{1\le j\le m} U_j
\end{gather*}
holds. Then ${\delta_1\le \alpha}$ holds. 
\end{enumerate} 
\end{Lem}
\vspace{5mm}
\textbf{proof}(Lemma 5):
\begin{enumerate}
\item We can prove this lemma by applying the arguments used in the Gromov compactness theorem (see Chapter 6 in \cite{AD}). Assume that ${\{u_i\}}$ is a sequence of solutions such that ${E(u_i)>0}$ and ${E(u_i)\rightarrow 0}$ hold. Then, the positive and the negative end of solutions ${\lim_{s\rightarrow \pm \infty}u_i(s,t)}$ converge to $1$-periodic orbits (Theorem 6.5.6 in \cite{AD}. $H$ is assumed to be non-degenerate in Theorem 6.5.6, but only the finiteness of the number of $1$-periodic orbits is used in the proof of Theorem 6.5.6). We denote these positive and negative periodic orbits by ${y_{\pm}^{i}(t)}$:
\begin{gather*}
\lim_{s\rightarrow \pm \infty}u_i(s,t)=y_{\pm}^i(t).
\end{gather*}
Assume that ${y_+^i=y_-^i}$ holds and a family of loops ${z_s(t)=u_i(s,t)}$ is contained in some ${U_j}$ for all ${s\in \mathbb{R}}$. Then,
\begin{gather*}
E(u_i)=\int_{\mathbb{R}\times S^1}u_i^*\omega+\int_0^1H_i(t,y_-^i(t))dt-\int_0^1H_i(t,y_+^i(t))dt \\ =\int_{\mathbb{R}\times S^1}u_i^*\omega=0
\end{gather*}
holds. So, ${E(u_i)>0}$ implies that there is a sequence ${(a_i,b_i)}$ such that
\begin{gather*}
u_i(a_i,b_i)\in \bigcup_{1\le j\le m}\partial U_j
\end{gather*}
holds. We define shifted solutions ${\{v_i\}}$ as follows:
\begin{gather*}
v_i(s,t)=u_i(a_i+s,t).
\end{gather*}
By taking a subsequence, we assume that ${b_i}$ converges to ${b\in S^1}$. Note that ${E(u_i)\rightarrow 0}$ implies that sphere bubbles do not occur in the limit. So we can apply Theorem 6.5.4 in \cite{AD} (Note that the symplectically aspherical condition (${\omega|_{\pi_2(M)}=0}$) is assumed in \cite{AD}, but we can remove this condition here because we exclude sphere bubbles). Theorem 6.5.4 in \cite{AD} implies that there is a solution $v$ as follows:
\begin{gather*}
v_i(s,t)\longrightarrow v(s,t) \ \ \mathrm{in} \ C_{\mathrm{loc}}^{\infty}(\mathbb{R}\times S^1,M)  \\
\partial_sv(s,t)+J_t(\partial_tv(s,t)-X_{H_t}(v(s,t)))=0  \\
v(0,b)\in \bigcup_{1\le j\le m}\partial U_j  \\
[v(0,t)]=\gamma  .
\end{gather*}
Let ${U\subset \mathbb{R}\times S^1}$ be any relatively compact subset. Then, 
\begin{gather*}
    0\le \int_{U}|\partial_sv(s,t)|^2dsdt=\lim_{i\rightarrow\infty}\int_{U}|\partial_sv_i(s,t)|^2dsdt\le\lim_{i\rightarrow \infty}E(v_i)=0
\end{gather*}
holds. In particular, 
\begin{gather*}
    E(v)=\int_{\mathbb{R}\times S^1}|\partial_sv(s,t)|^2dsdt=0
\end{gather*}
holds. ${E(v)=0}$ implies that $v$ is a trivial solution. In other words, ${v(s,t)}$ does not depend on $s$ and ${l(t)=v(s,t)}$ is a periodic orbit. This is a contradiction because the third and the fourth conditions imply that ${v(0,t)}$ is not a periodic orbit in $\gamma$.
\item
${E(u_i)\le \alpha\le \delta_2}$ implies that there is no sphere bubbles in the limit of $u_i$. So we can apply Theorem 6.5.4 in \cite{AD} again. Note that Theorem 6.5.4 in \cite{AD} treats fixed Hamiltonian function $H$ (we treat a ${C^{\infty}}$-convergent sequence of Hamiltonian functions ${\{H_i\}}$). However, the proof of the theorem can be applied to a ${C^{\infty}}$-convergent sequence of Hamiltonian functions (and a $C^{\infty}$-convergent sequence of almost complex structures) because $C^{\infty}$-convergence implies that the Sobolev norms ${W^{k,p}(H_i)}$ are uniformely bounded for fixed $p$ and $k$. See Chapter ${12}$ of \cite{AD} where analytic foundations of Theorem 6.5.4 are proved (in particular, see the statement of Theorem 12.1.5 and its proof).

Let ${w_i}$ be shifted solutions as follows:
\begin{gather*}
w_i(s,t)=u_i(s_i+s,t).
\end{gather*}
By taking a subsequence, we assume that $t_i$ converges to $T$. Then, there is a solution ${w:\mathbb{R}\times S^1\rightarrow M}$ such that
\begin{gather*}
w_i(s,t)\longrightarrow w(s,t) \ \ \mathrm{in} \ C_{\mathrm{loc}}^{\infty}(\mathbb{R}\times S^1,M)  \\
w(0,T)\in \overline{M\backslash \cup U_j} \\
[w(0,t)]=\gamma \\
E(w)\le \alpha
\end{gather*}
holds. So, $w$ is a non-trivial solution (${E(w)>0}$). This implies that
\begin{gather*}
\delta_1\le E(w)\le \alpha <\delta_2
\end{gather*}
holds.
\end{enumerate}
\begin{flushright} $\Box$ \end{flushright}

\begin{Lem}
Let $x$ be an isolated $1$-periodic orbit of $H$ and let ${U\subset S^1\times M}$ be a sufficiently small open neighborhood of $x$. We fix a sequence of Hamiltonian functions ${H_i}$ which converges to $H$ smoothly so that $x$ splits into non-degenerate $1$-periodic orbits ${\{x_1^{(i)},\cdots,x_{m_i}^{(i)}\}}$. Assume that ${u_i}$ and ${(s_i,t_i)}$ satisfy the following conditions:
\begin{gather*}
u_i:\mathbb{R}\times S^1\longrightarrow M \\
\partial_su_i(s,t)+J_t(\partial_tu_i(s,t)-X_{(H_i)_t}(u_i(s,t)))=0 \\
u_i(s_i,t_i)\in \partial U  \\
\lim_{s\rightarrow - \infty}u_i(s,t)=x_{k(i)}^{(i)}(t) \ \ \ (k(i)\in \{1,\cdots,m_i\}).
\end{gather*}
Then, there is ${\epsilon>0}$ such that ${E(u_i)\ge \epsilon}$ holds for all $u_i$.
\end{Lem}
\vspace{5mm}
\textbf{proof}(Lemma 6):
We assume that ${E(u_i)\rightarrow 0}$ holds. Without loss of generality, we assume that a loop ${l^{(i)}_s(t)=u_i(s,t)}$ is contained in ${U}$ for all ${s<s_i}$. We consider the following half cylinders $v_i$:
\begin{gather*}
v_i:(-\infty,0]\times S^1\rightarrow M \\
v_i(s,t)=u_i(s_i+s,t).
\end{gather*}
By taking a subsequence, we assume that ${t_i}$ converges to $T$. As in the proof of Lemma 5, we can see that there is a half cylinder ${v}$ such that
\begin{gather*}
v_i(s,t)\longrightarrow v(s,t) \ \ \mathrm{in} \ C_{\mathrm{loc}}^{\infty}((-\infty,0]\times S^1,M)  \\
\partial_sv(s,t)+J_t(\partial_tv(s,t)-X_{H_t}(v(s,t)))=0  \\
v(0,T)\in \partial U \\
E(v)=0
\end{gather*}
holds. This is a contradiction.
\begin{flushright} $\Box$ \end{flushright}

Next, we apply Lemma 6 to construct the local version of the Floer cohomology theory, the so-called local Floer cohomology (see \cite{GG5}). Let ${x\in P(H,\gamma)}$ be an isolated periodic orbit of a Hamiltonian function ${H\in C^{\infty}(S^1\times M)}$. In this subsection, we explain the local Floer cohomology ${HF^{loc}(H,x)}$. Note that $x$ is not necessarily a non-degenerate periodic orbit. Let ${U\subset S^1\times M}$ be a sufficiently small open neighborhood of the embedded circle ${\{(t,x(t))\}\subset S^1\times M}$ and let ${\widetilde{H}}$ be a non-degenerate ${C^{\infty}}$-small perturbation of $H$. Then ${x}$ splits into non-degenerate periodic orbits ${\{x_1,\cdots,x_k\}}$ of ${\widetilde{H}}$ in ${U}$. The local Floer cochain complex ${CF^{loc}(H,x)}$ is generated by these perturbed periodic orbits:

\begin{gather*}
CF^{loc}(H,x)=\bigoplus_{i=1}^k \mathbb{K}\cdot x_i.
\end{gather*}
The coboundary operator ${d_F^{loc}}$ on ${CF^{loc}(H,x)}$ is defined by counting solutions of the Floer equation in ${U}$. So, ${d_F^{loc}}$ can be written by ${d_F^{loc}(x_i)=\sum n(x_i,x_j)x_j}$ where the coefficient ${n(x_i,x_j)\in \mathbb{K}}$ is determined by the number of solutions of the following equation modulo the natural ${\mathbb{R}}$-action:

\begin{gather*}
u:\mathbb{R}\times S^1\longrightarrow M \\
\partial_su(s,t)+J_t(u(s,t))(\partial_tu(s,t)-X_{\widetilde{H}_t}(u(s,t)))=0  \\
\lim_{s\to -\infty}u(s,t)=x_i(t), \lim_{s\to +\infty}u(s,t)=x_j(t)   \\
(t,u(s,t))\in U.
\end{gather*}

\begin{Rem}
Lemma 6 implies that any solution of the above Floer equation with small energy does not escape the isolating neighborhood of periodic orbits ${U}$ if the perturbation ${\widetilde{H}}$ is sufficiently close to the original $H$. 
\end{Rem}

The local Floer cohomology ${HF^{loc}(H,x)}$ is the cohomology of the cochain complex ${(CF^{loc}(H),d_F^{loc})}$. The following properties hold (see \cite{GG5}).

\begin{enumerate}
\item If $x$ is a non-degenerate periodic orbit of ${H}$, then ${HF^{loc}(H,x)\cong \mathbb{K}}$ holds.
\item ${HF^{loc}(H,x)}$ does not depend on the choice of the perturbation ${\widetilde{H}}$ if it is sufficiently small.
\item Let ${\{\theta_1,\cdots,\theta_l\}\subset S^1\setminus \{1\}}$ be the set of eigenvalues of the differential map 
\begin{equation*}
d\phi_H|_{x(0)}:T_{x(0)}M\longrightarrow T_{x(0)}M
\end{equation*}
on ${S^1\setminus \{1\}\subset \mathbb{C}}$. An integer ${k\in \mathbb{N}}$ is called admissible if ${\theta_i^k\neq 1}$ holds for all ${1\le i \le l}$. If $k$ is admissible, the local Floer cohomology of ${(H,x)}$ and ${(H^{(k)}, x^{(k)})}$ are isomorphic. In other words,
\begin{equation*}
HF^{loc}(H^{(k)},x^{(k)})\cong HF^{loc}(H,x) 
\end{equation*}
holds.
\end{enumerate}

Note that the local Floer cochain complex ${(CF^{loc}(H,x),d_F^{loc})}$ is a cochain complex over a finite-dimensional vector space over ${\mathbb{K}}$ (not over $\Lambda$ nor $\Lambda_0$).

\section{$\mathbb{Z}_p$-equivariant Floer theory}

Our proof of Theorem 2 is an application of the ${\mathbb{Z}_p}$-equivariant Floer cohomology. In this section, we briefly review constructions and basic properties of the equivariant theory which we will use in our proof of Theorem 2. The readers can find the detailed constructions of the equivariant Floer theory in \cite{SZ} (see also \cite{Se}). The first attempt in this direction is due to Seidel in \cite{Se}, where he constructed ${\mathbb{Z}_2}$-equivariant Floer theory. After that, Seidel's construction was generalized to ${\mathbb{Z}_p}$-equivariant Floer theory by Shelukhin-Zhao and Shelukhin in \cite{SZ,Sh}. In these papers, they gave various applications of the equivariant Floer cohomology. In particular, the equivariant Floer theory is very useful when we study the relationship between ${HF(\phi)}$ and the Floer cohomology ${HF(\phi^p)}$ associated to a prime iterate of ${\phi}$. As in \cite{Sh}, we will focus on the behavior of the torsion of ${HF(\phi^p)}$ as $p$ becomes bigger and bigger.

We fix a prime number $p$ and assume that ${H\in C^{\infty}(S^1\times M)}$ is a Hamiltonian function such that ${H^{(p)}}$ is non-degenerate. In this section, we assume that ${(M,\omega)}$ is a toroidally monotone symplectic manifold. 
\begin{Def}
A symplectic manifold ${(M,\omega)}$ is called toroidally monotone if there is a constant ${\lambda \ge 0}$ such that 
\begin{gather*}
\int_{\mathbb{T}^2}v^*\omega=\lambda \int_{\mathbb{T}^2}v^*c_1
\end{gather*}
holds for every smooth map ${v:\mathbb{T}^2\rightarrow M}$.
\end{Def}
We will give a construction of the equivariant theory for weakly monotone symplectic manifolds in the last section, where we explain what the technical difficulty is in the weakly monotone case and how we can overcome this problem. The equivariant Floer cohomology is a mixture of the Floer theory on ${(M,\omega)}$ and the Morse theory on the classifying space. 

Let ${\mathbb{C}^{\infty}}$ be the infinite-dimensional complex vector space
\begin{equation*}
\mathbb{C}^{\infty}=\{z=(z_k)_{k\in \mathbb{Z}_{\ge 0}} \ | \ z_k\in \mathbb{C}, z_k=0 \textrm{ \ for\ sufficiently\ large\ }k\}
\end{equation*}
and let ${S^{\infty}\subset \mathbb{C}^{\infty}}$ be the infinite-dimensional sphere defined by
\begin{equation*}
S^{\infty}=\{z\in \mathbb{C}^{\infty} \ | \ \sum|z_k|^2=1\}.
\end{equation*}
There is a natural ${\mathbb{Z}_p}$-action and a shift operator ${\tau}$ on ${S^{\infty}}$ as follows:

\begin{gather*}
(m\cdot z)_k=\exp \bigg(\frac{2\pi im}{p}\bigg)\cdot z_k  \ \ (m\in \mathbb{Z}_p)  
\end{gather*}
\begin{gather*}
(\tau(z))_k=\begin{cases}  0 \  \ \ \ (k=0) \\ z_{k-1} \ \ \ \  (k\ge 1). \end{cases}
\end{gather*}
Let ${f:S^{\infty}\rightarrow \mathbb{R}}$ be the following Morse-Bott function:
\begin{equation*}
f(z)=\sum_{k=0}^{\infty}k|z_k|^2.
\end{equation*}
This $f$ satisfies ${\tau^*f=f+1}$ and ${f(z)=f(m\cdot z)}$ for ${m\in \mathbb{Z}_p}$. The critical submanifolds are ${S_l^1=\{z\in S^{\infty} \ | \ |z_l|=1\}}$ and their indexes are ${2l}$. Next, we perturb $f$ to a Morse function with the above properties. For example, we can use the following explicit perturbation (see \cite{SZ}):

\begin{equation*}
\widetilde{F}=f(z)+\epsilon \cdot \sum_k\textrm{Re}(z_k^p).
\end{equation*}
Then ${\widetilde{F}}$ is a Morse function and ${S_l^1}$ splits into critical points ${\{Z_{2l}^0,\cdots,Z_{2l}^{p-1}\}}$ and ${\{Z_{2l+1}^0,\cdots,Z_{2l+1}^{p-1}\}}$. The Morse index of ${Z_j^m}$ is $j$. It is straightforward to see that ${\tau^*\widetilde{F}=\widetilde{F}+1}$ and ${\widetilde{F}}(z)=\widetilde{F}(m\cdot z)$ hold. We also fix a Riemannian metric ${\widetilde{g}}$ such that ${\widetilde{g}}$ is invariant under ${\mathbb{Z}_p}$-action and $\tau$ and ${(\widetilde{F},\widetilde{g})}$ is a Morse-Smale pair.

The Morse coboundary operator ${d_M}$ of the pair ${(\widetilde{F},\widetilde{g})}$ can be written in the following form:

\begin{gather*}
d_M(Z_{2l}^i)=Z_{2l+1}^i-Z_{2l+1}^{i+1} \\
d_M(Z_{2l+1}^i)=Z_{2l+2}^0+Z_{2l+2}^1+\cdots +Z_{2l+2}^{p-1} \ \ \ (i\in \mathbb{Z}_p).
\end{gather*}
This ${d_M}$ is equivariant under the ${\mathbb{Z}_p}$ action ${m\cdot Z_i^j=Z_i^{j+m(\textrm{mod} p)}}$. Let ${(H^{(p)}, J_t)}$ be a pair consisting of a Hamiltonian function ${H^{(p)}}$ and an ${S^1}$-dependent almost complex structure on ${(M,\omega)}$. We want to extend ${J_t}$ to a family of ${S^1}$-dependent almost complex structures parametrized by ${w\in S^{\infty}}$. So, we consider a family of almost complex structures ${\{J_{w,t}\}_{(w,t)\in S^{\infty}\times S^1}}$ which satisfies the following conditions:

\begin{itemize}
\item (locally constant at critical points) For all $w$ in a small neighborhood of ${Z_i^m}$, 
\begin{equation*}
J_{w,t}=J_{t-\frac{m}{p}} .
\end{equation*}
\item (${\mathbb{Z}_p}$-equivariance) ${J_{m\cdot w,t}=J_{w,t-\frac{m}{p}}}$ holds for any ${m\in \mathbb{Z}_p}$ and ${w\in S^{\infty}}$.
\item (invariance under the shift ${\tau}$) ${J_{\tau(w),t}=J_{w,t}}$ holds.
\end{itemize}

We consider the following equation, which is a mixture of the Floer equation and the Morse equation for ${x,y\in P(H^{(p)})}$, ${m\in \mathbb{Z}_p}$, ${\lambda \ge 0}$, ${\alpha \in \{0,1\}}$, ${i\in \mathbb{Z}_{\ge 0}}$:

\begin{gather*}
(u,v)\in C^{\infty}(\mathbb{R}\times S^1, M)\times  C^{\infty}(\mathbb{R}, S^{\infty})  \\
\partial _su(s,t)+J_{v(s),t}(u(s,t))(\partial_tu(s,t)-X_{H^{(p)}}(u(s,t)))=0  \\
\frac{d}{ds}v(s)-\textrm{grad}\widetilde{F}=0  \\
\lim_{s\to -\infty}v(s)=Z_{\alpha}^0, \lim_{s\to +\infty}v(s)=Z_i^m, \lim_{s\to -\infty}u(s,t)=x(t), \lim_{s\to +\infty}u(s,t)=y(t-\frac{m}{p}) \\
\int_{\mathbb{R}\times S^1}u^*\omega+\int_0^1H^{(p)}(t,x(t))-H^{(p)}(t,y(t))dt=\lambda.
\end{gather*}
We denote the space of solutions of this equation by ${\mathcal{M}_{\alpha,i,m}^{\lambda}(x,y)}$.
\begin{equation*}
\mathcal{M}_{\alpha,i,m}^{\lambda}(x,y)=\{(u,v) \ | \ (u,v) \ \textrm{satisfies\ the\ above\ equation}\}/\sim
\end{equation*}
The equivalence relation ${\sim}$ is given by the natural ${\mathbb{R}}$-action on the solution space.

\begin{Rem}
We have to perturb the above equation to achieve transversality of the linearized operator. This can be achieved by a perturbation in a compact region ${S^N\subset S^{\infty}}$ because we assumed that ${(M,\omega)}$ is toroidally monotone. See \cite{Se,Sh,SZ}.
\end{Rem}

We define ${d_{\alpha}^{i,m}}$ (${\alpha \in \{0,1\}}$, ${i\in \mathbb{Z}_{\ge 0}}$, ${m\in \mathbb{Z}_p}$) as follows:

\begin{gather*}
d_{\alpha}^{i,m}:CF(H^{(p)},\gamma:\Lambda_0)\longrightarrow CF(H^{(p)},\gamma:\Lambda_0) \\
x\mapsto \sum_{y\in P(H^{(p)}), \lambda \ge 0}\sharp \mathcal{M}_{\alpha,i,m}^{\lambda}(x,y)\cdot T^{\lambda}y.
\end{gather*}
Let ${d_{\alpha}^i}$ be the sum ${d_{\alpha}^i=\sum_{m\in \mathbb{Z}_p}d_{\alpha}^{i,m}}$. Note that ${d_0^0=d_1^1=d_F}$ holds. We define the ${\mathbb{Z}_p}$-equivariant Floer cochain complex as follows:
\begin{equation*}
CF_{\mathbb{Z}_p}(H^{(p)},\gamma)=CF(H^{(p)},\gamma:\Lambda_0)\otimes \Lambda_0[[u]]\langle \theta\rangle.
\end{equation*}
Here, ${\langle \theta \rangle}$ is the exterior algebra on the formal variable ${\theta}$ and ${\Lambda_0[[u]]}$ is the ring of formal power series of the formal variable ${u}$. So, any element of ${CF_{\mathbb{Z}_p}(H^{(p)},\gamma)}$ is written in the following form for some ${k\in \mathbb{Z}}$:

\begin{equation*}
(x_k+y_k\theta)u^k+\sum_{i=k+1}^{\infty}(x_i+y_i\theta)u^i \ \ \ (x_j,y_j\in CF(H^{(p)},\gamma:\Lambda_0)).
\end{equation*}
The equivariant Floer coboundary operator ${d_{eq}}$ is a ${\Lambda_0[[u]]}$-linear map defined as follows:

\begin{gather*}
d_{eq}(x\otimes 1)=\sum_{i=0}^{\infty}d_0^{2i}(x)\otimes u^i+\sum_{i=0}^{\infty}d_0^{2i+1}(x)\otimes u^i\theta \\
d_{eq}(x\otimes \theta)=\sum_{i=0}^{\infty}d_1^{2i+1}(x)\otimes u^i\theta+\sum_{i=1}^{\infty}d_1^{2i}(x)\otimes u^i.
\end{gather*}
${d_0^0=d_1^1=d_F}$ implies that ${d_{eq}}$ is a sum of ${d_F}$ and higher terms. The ${\mathbb{Z}_p}$-equivariant Floer cohomology ${HF_{\mathbb{Z}_p}(H^{(p)},\gamma)}$ is defined by
\begin{equation*}
HF_{\mathbb{Z}_p}(H^{(p)},\gamma)=H((CF_{\mathbb{Z}_p}(H^{(p)},\gamma),d_{eq})).
\end{equation*}

We also consider the coefficient extension of the ${\mathbb{Z}_p}$-equivariant Floer cochain complex. Let ${\Lambda_0[u^{-1},u]]}$ be the ring of the completion of Laurent polynomials:
\begin{gather*}
\Lambda_0[u^{-1},u]]=\bigg\{\sum_{i=k}^{\infty}a_iu^i \ \bigg| \ k\in \mathbb{Z},a_i\in \Lambda_0 \bigg\}.
\end{gather*}
The ${\mathbb{Z}_p}$-equivariant Tate Floer cochain complex of ${H^{(p)}}$ is defined as follows:
\begin{gather*}
\widehat{CF}_{\mathbb{Z}_p}(H^{(p)},\gamma)=CF(H^{(p)},\gamma:\Lambda_0)\otimes \Lambda_0[u^{-1},u]]\langle \theta \rangle  \\
\widehat{d}_{eq}:\widehat{CF}_{\mathbb{Z}_p}(H^{(p)},\gamma) \longrightarrow \widehat{CF}_{\mathbb{Z}_p}(H^{(p)},\gamma).
\end{gather*}
Here, ${\widehat{d}_{eq}}$ is the natural extension of ${d_{eq}}$. The ${\mathbb{Z}_p}$-equivariant Tate Floer cohomology of ${H^{(p)}}$ is the cohomology of this cochain complex ${(\widehat{CF}_{\mathbb{Z}_p}(H^{(p)},\gamma),\widehat{d}_{eq})}$ as follows:
\begin{gather*}
\widehat{HF}_{\mathbb{Z}_p}(H^{(p)},\gamma)=H(\widehat{CF}_{\mathbb{Z}_p}(H^{(p)},\gamma),\widehat{d}_{eq}).
\end{gather*}

\begin{Rem}
We consider coefficient extensions (in other words, Tate complexes) because the ${\mathbb{Z}_p}$-equivariant pair of pants product gives a local isomorphism between the local Floer cohomology and the local ${\mathbb{Z}_p}$-equivariant Tate Floer cohomology (see \cite{Se,SZ}). We will explain this below and we will apply this local isomorphism.
\end{Rem}

\section{Local (${\mathbb{Z}_p}$-equivariant Tate) Floer cohomology and homological perturbation theory}

In the proof of Theorem 2, we have to treat Hamiltonian diffeomorphisms with finitely many periodic orbits, which are not necessarily non-degenerate. One natural way to overcome this difficulty is to just perturb the original Hamiltonian function ${H}$ to a non-degenerate Hamiltonian function ${\widetilde{H}}$ and consider ${HF(\widetilde{H}:\Lambda_0)}$ instead of ${HF(H:\Lambda_0)}$. However, it is not sufficient because the structure of ${HF(\widetilde{H}:\Lambda_0)}$ does depend on the perturbation ${\widetilde{H}}$. We have to construct a cochain complex and a cohomology theory in a homologically canonical way. In this section, we explain the construction of Floer theory for possibly degenerate Hamiltonian diffeomorphisms (see also \cite{Sh}). This is an application of the homological perturbation theory in \cite{M}.

Assume that ${H}$ is a possibly degenerate Hamiltonian function and ${P(H,\gamma)}$ is finite for ${\gamma \neq 0\in H_1(M:\mathbb{Z})/\textrm{Tor}}$. Let ${\widetilde{H}}$ be a small perturbation of ${H}$ such that ${\phi_{\widetilde{H}}}$ is non-degenerate. Then every element ${x_i\in P(H,\gamma)}$ splits into finitely many non-degenerate periodic orbits ${Q(\widetilde{H},x_i)=\{x_i^{1},\cdots,x_i^{l_i}\}\subset P(\widetilde{H},\gamma)}$. We also fix a small isolating neighborhood of each ${x_i\in P(H,\gamma)}$. Let ${U_i\subset S^1\times M}$ be a sufficiently small open neighborhood of ${(t,x_i(t))}$. First, we deform the complex ${(CF(\widetilde{H},\gamma:\Lambda_0),d_F)}$ in a canonical way. For any ${x_i\in P(H,\gamma)}$ and ${x_i^j\in Q(\widetilde{H},x_i)}$, we fix a sufficiently small connecting cylinder between ${x_i}$ and ${x_i^j}$ which is contained in the isolating neighborhood ${U_i}$ as follows:
\begin{gather*}
v_i^j:[0,1]\times S^1\longrightarrow U_i \\
v_i^j(0,t)=x_i(t), v_i^j(1,t)=x_i^j(t).
\end{gather*}
The ``gap" of the action functional
\begin{equation*}
c(x_i,x_i^j)=\int_{[0,1]\times S^1}(v_i^j)^*\omega+\int_0^1H(t,x_i(t))-\widetilde{H}(t,x_i^j(t))dt
\end{equation*}
does not depend on the choice of ${v_i^j}$. We define a correction map ${\kappa}$ as follows:

\begin{gather*}
\kappa: CF(\widetilde{H},\gamma:\Lambda)\rightarrow CF(\widetilde{H},\gamma:\Lambda) \\
x_i^j\mapsto T^{c(x_i,x_i^j)}\cdot x_i^j.
\end{gather*}
Note that ${\kappa}$ is defined over ${\Lambda}$ (not over ${\Lambda_0}$). Then we can define the modified coboundary operator ${\widetilde{d_F}}$ by ${\widetilde{d_F}=\kappa^{-1}\circ d_F\circ \kappa}$. It is easy to see that ${\widetilde{d_F}}$ can be defined over ${\Lambda_0}$ (see Lemma 11). Note that ${\kappa}$ gives an isomorphism between Floer cochain complexes over ${\Lambda}$ (but not over ${\Lambda_0}$) as follows:
\begin{gather*}
\kappa_{\widetilde{H}}: (CF(\widetilde{H},\gamma:\Lambda), \widetilde{d_F}) \longrightarrow (CF(\widetilde{H},\gamma:\Lambda), d_F)  \\
x\mapsto \kappa(x).
\end{gather*}
Let ${\widetilde{H}'}$ be another non-degenerate perturbation of ${H}$ and we also denote the modified coboundary operator on ${CF(\widetilde{H}',\gamma:\Lambda)}$ by ${\widetilde{d_F}}$. It is straightforward to see that the natural continuation map between the Floer cochain complexes
\begin{equation*}
\Phi: (CF(\widetilde{H},\gamma:\Lambda),d_F)\longrightarrow (CF(\widetilde{H}',\gamma:\Lambda),d_F)
\end{equation*}
descends to a continuation map over ${\Lambda_0}$ (not only over ${\Lambda}$) as follows:
\begin{gather*}
\widetilde{\Phi}:  (CF(\widetilde{H},\gamma:\Lambda_0),\widetilde{d_F})\longrightarrow (CF(\widetilde{H}',\gamma:\Lambda_0),\widetilde{d_F}) \\
x\mapsto (\kappa_{\widetilde{H}'})^{-1}\circ \Phi \circ \kappa_{\widetilde{H}}.
\end{gather*}

\begin{Lem}
The modified coboundary operator ${\widetilde{d_F}}$ is defined over $\Lambda_0$ if ${\widetilde{H}}$ is a sufficiently small perturbation of ${H}$. Moreover, the modified continuation map ${\widetilde{\Phi}}$ is defined over ${\Lambda_0}$ if ${\widetilde{H}}$ and ${\widetilde{H}'}$ are sufficiently close to $H$.
\end{Lem}
\vspace{5mm}
\textbf{proof}(Lemma 11):

Let ${\widetilde{H}_i}$ be a series of perturbations of ${H}$ satisfying the following conditions.
\begin{itemize}
\item ${\widetilde{H}_i}$ converges to ${H}$ in ${C^{\infty}}$-topology
\item ${\widetilde{d_F}}$ is not defined over ${\Lambda_0}$ for each ${\widetilde{H}_i}$
\end{itemize}
The second assumption implies that we can choose a sequence of cylinders $u_i$ as follows:
\begin{gather*}
u_i:\mathbb{R}\times S^1\longrightarrow M \\
\partial_su_i(s,t)+J_t(\partial_tu_i(s,t)-X_{\widetilde{H}_i}(u_i(s,t)))=0 \\
\exists (s_i,t_i)\in \mathbb{R}\times S^1 \ \ s.t. \ \ u_i(s_i,t_i)\in \bigcup_j\partial U_j  \\
E(u_i)=\int_{\mathbb{R}\times S^1}|\partial_su_i(s,t)|^2dsdt\le 2\max\big\{|c(x_i,x_i^j)| \ | \ x_i^j\in Q(\widetilde{H}_i,x_i)\big\}.
\end{gather*}
Let $v_i$ be a shifted solution:
\begin{gather*}
v_i(s,t)=v(s_i+s,t).
\end{gather*}
By taking a subsequence, we assume that ${t_i}$ converges to ${T\in S^1}$. As in the proof of Lemma 5 and Lemma 6, we can see that there is a cylinder $v$ as follows:
\begin{gather*}
v_i(s,t)\longrightarrow v(s,t) \ \ \mathrm{in} \ C_{\mathrm{loc}}^{\infty}(\mathbb{R}\times S^1,M)  \\
\partial_sv(s,t)+J_t(\partial_tv(s,t)-X_H(s,t))=0  \\
v(0,T)\in \bigcup_j \partial U_j \\
E(v)=\int_{\mathbb{R}\times S^1}|\partial_su(s,t)|^2dsdt=0.
\end{gather*}
This is a contradiction because ${E(v)=0}$ implies that ${v}$ is a ``constant" solution (${\partial_sv(s,t)\equiv 0}$), but the image of $v$ is not contained in the union of isolating neighborhoods ${U_i}$. A similar argument can be used to prove that ${\widetilde{\Phi}}$ is defined over ${\Lambda_0}$ if perturbations are sufficiently close to $H$.
\begin{flushright}
    $\Box$
\end{flushright}

The continuation map in the inverse direction 
\begin{gather*}
\Psi: (CF(\widetilde{H}',\gamma:\Lambda),d_F)\longrightarrow (CF(\widetilde{H},\gamma:\Lambda),d_F)
\end{gather*}
also descends to a continuation map 
\begin{equation*}
\widetilde{\Psi}:  (CF(\widetilde{H}',\gamma:\Lambda_0),\widetilde{d_F})\longrightarrow (CF(\widetilde{H},\gamma:\Lambda_0),\widetilde{d_F})
\end{equation*}
in the same way. The cochain homotopy maps between the identity and ${\Psi \circ \Phi}$, ${\Phi \circ \Psi}$ also descend to cochain homotopy maps between the identity and ${\widetilde{\Psi}\circ \widetilde{\Phi}}$, ${\widetilde{\Phi}\circ \widetilde{\Psi}}$. This implies that the perturbed and modified Floer cochain complex ${(CF(\widetilde{H},\gamma:\Lambda_0),\widetilde{d_F})}$ is unique up to cochain homotopy equivalence. So we can define the Floer cohomology of ${H}$ by the cohomology of this cochain complex and this is well defined:

\begin{equation*}
HF(H,\gamma:\Lambda_0)=H(CF(\widetilde{H},\gamma:\Lambda_0),\widetilde{d_F}).
\end{equation*}

However, this construction is not sufficient for our purpose. The cochain complex ${(CF(\widetilde{H},\gamma:\Lambda_0),\widetilde{d_F})}$ is not strict in the following sense. ``Strict" means that there is ${\epsilon>0}$ such that 
\begin{equation*}
\widetilde{d_F}(z)\in CF^{\epsilon}(\widetilde{H},\gamma:\Lambda_0)
\end{equation*}
holds for all ${z\in CF(\widetilde{H},\gamma:\Lambda_0)}$. Note that we can identify ${CF(\widetilde{H},\gamma:\Lambda_0)}$ and ${\bigoplus_{x\in P(H,\gamma)}CF^{loc}(H,x)\otimes \Lambda_0}$. Then ${\widetilde{d_F}}$ can be decomposed into the sum of ${\widetilde{d_F}=d_F^{loc}+D}$ where ${D}$ is a higher energy term:
\begin{equation*}
D(CF(\widetilde{H},\gamma:\Lambda_0))\subset CF^{\epsilon}(\widetilde{H},\gamma:\Lambda_0)  \ \ \ \ (\epsilon>0).
\end{equation*}

The existence of such a decomposition for sufficiently small perturbation ${\widetilde{H}}$ follows from the following arguments. Let ${\delta_1>0}$ be the smallest energy of non-constant solutions of the Floer equations and let ${\delta_2}$ be the smallest energy of non-constant pseudo-holomorphic spheres as follows:
\begin{gather*}
\delta_1=\min\Big\{E(u)>0 \ \Big| \ \begin{matrix}u:\mathbb{R}\times S^1\longrightarrow M, \ [u(0,t)]=\gamma   \\ \partial_su+J_t(\partial_tu(s,t)-X_H(s,t))=0\end{matrix}\Big\}  \\
\delta_2=\min\Big\{E(v)=\int_{\mathbb{CP}^1}v^*\omega>0 \ \Big|  \ \begin{matrix} v:\mathbb{CP}^1\longrightarrow M \\ J_t\circ dv=dv\circ j_{\mathbb{CP}^1} \end{matrix} \Big\}.
\end{gather*}
Note that $\delta_1$ and $\delta_2$ are strictly positive (Lemma 5). We fix a constant ${0<\epsilon<\min\{\delta_1,\delta_2\}}$. If the existence of such decomposition is false, we can find a sequence of perturbations ${\widetilde{H}_i}$ and a sequence of solutions of the Floer equations ${\{u_i\}}$ as follows:
\begin{itemize}
\item $\widetilde{H}_i$ converges to $H$ in ${C^{\infty}}$-topology
\item ${\partial_su_i(s,t)+J_t(\partial_tu_i(s,t)-X_{\widetilde{H}_i})(u_i(s,t))=0}$
\item $E(u_i)\le \epsilon$
\item $\exists (s_i,t_i)\in \mathbb{R}\times S^1$ s.t. ${u_i(s_i,t_i)\in \bigcup_j \partial U_j}$
\end{itemize}
As in the proof of Lemma 5 and Lemma 6, we can see that there is a solution of the Floer equation as follows:
\begin{gather*}
u:\mathbb{R}\times S^1\longrightarrow M \\
\partial_su(s,t)+J_t(\partial_tu(s,t)-X_H(u(s,t)))=0  \\
\exists T\in S^1 \ \ \mathrm{s.t} \ \ u(0,T)\in \bigcup_j\partial U_j  \\
[u(0,t)]=\gamma \\
E(u)=\int_{\mathbb{R}\times S^1}|\partial_su(s,t)|^2dsdt\le \epsilon.
\end{gather*}
This is a contradiction because the above conditions imply that $u$ is a non-constant solution whose energy is smaller than $\delta_1$.

Note that ${(CF(\widetilde{H},\gamma:\Lambda_0),\widetilde{d_F})}$ is strict if and only if ${d_F^{loc}=0}$ holds. Our aim in the rest of this section is to construct a homologically canonical coboundary operator on  ${\bigoplus_{x\in P(H,\gamma)}HF^{loc}(H,x)\otimes \Lambda_0}$ instead of on  ${\bigoplus_{x\in P(H,\gamma)}CF^{loc}(H,x)\otimes \Lambda_0}$ and prove that the cochain complex is cochain homotopy equivalent to the original cochain complex ${(CH(\widetilde{H},\gamma:\Lambda_0),\widetilde{d_F})}$. We define the Floer cochain complex of $H$ by 
\begin{equation*}
CF(H,\gamma:\Lambda_0)=\bigoplus_{x\in P(H,\gamma)}HF^{loc}(H,x)\otimes \Lambda_0 .
\end{equation*}
We choose a basis $\{X_i^a,Y_i^b,Z_i^c\}$ of ${CF^{loc}(H,x_i)=\bigoplus \mathbb{F}_p\cdot x_i^j}$ which satisfies the following relations:
\begin{gather*}
d_F^{loc}(X_i^a)=0 \ \ \ (1\le a\le \textrm{dim}_{\mathbb{F}_p}HF^{loc}(H,x_i)) \\
d_F^{loc}(Y_i^b)=0 \ \ \ (1\le b\le \frac{1}{2}(l_i-\textrm{dim}_{\mathbb{F}_p}HF^{loc}(H,x_i))) \\
d_F^{loc}(Z_i^c)=Y_i^c \ \ \ (1\le c\le \frac{1}{2}(l_i-\textrm{dim}_{\mathbb{F}_p}HF^{loc}(H,x_i))).
\end{gather*} 
Here, $l_i$ is the number of the perturbed periodic orbits (${Q(\widetilde{H},x_i)=\{x_i^1,\cdots,x_i^{l_i}\}}$). Then ${HF^{loc}(H,x_i)}$ can be identified with ${\textrm{span}_{\mathbb{F}_p}\langle X_i^1,\cdots,X_i^{\textrm{dim}_{\mathbb{F}_p}HF^{loc}(H,x_i)}\rangle}$. We define two operators ${\Pi_i}$ and ${\Theta_i}$ on ${CF^{loc}(H,x_i)}$ as follows:

\begin{gather*}
\Pi_i(\sum \alpha_aX_i^a+\sum \beta_bY_i^b+\sum \gamma_cZ_i^c)=\sum \alpha_aX_i^a\\
\Theta_i(\sum \alpha_aX_i^a+\sum \beta_bY_i^b+\sum \gamma_cZ_i^c)=\sum \beta_cZ_i^c.
\end{gather*}
Let ${\Pi=\sum \Pi_i}$ and ${\Theta=\sum \Theta_i}$ be the sums of the above operators. Next, we apply the perturbation theory in \cite{M}. First, we explain what perturbation theory is and how we can apply it to our case.

Let ${M=(M,d_M)}$ be a chain complex with a decreasing filtration ${\{F^pM\}_{p\in \mathbb{Z}_{\ge 0}}}$:
\begin{equation*}
M=F^0M\supset F^1M\supset F^2M\supset \cdots.
\end{equation*}
We assume that this filtration is complete. In other words, $M$ is complete with respect to the ${F^p}$-adic topology and an infinite sum ${\sum_{p} m_p}$ such that ${m_p\in F^pM}$ holds converges uniquely to an element of $M$. Let ${(N,d_N)}$ be another chain complex with a complete and decreasing filtration ${\{F^pN\}}$. Morphisms between filtered chain complexes are maps that preserve filtrations. Let ${f:M\rightarrow N}$ be a morphism. Another morphism ${g:M\rightarrow N}$ is called a perturbation of $f$ if 
\begin{equation*}
(f-g)F^pM\subset F^{p+1}N
\end{equation*}
holds for all ${p\in \mathbb{Z}_{\ge 0}}$. Perturbations of boundary operators $d_M$ and ${d_N}$ are defined in the same way. In \cite{M}, Markl studied the following problem. Let 
\begin{gather*}
F:(M,d_M)\longrightarrow (N,d_N) \\
G:(N,d_N)\longrightarrow (M,d_M)
\end{gather*}
be chain maps that preserve filtrations and assume that we also have chain homotopies between the identity and ${GF}$, ${FG}$ as follows:
\begin{gather}
H:M\rightarrow M, \ L:N\rightarrow N \\
GF-\textrm{Id}_M=d_MH+Hd_M   \\
FG-\textrm{Id}_N=d_NL+Ld_N.
\end{gather}
Next, we perturb the original boundary operator ${d_M}$ to a new boundary operator ${\widetilde{d_M}=d_M+D_M}$ (${D_M(F^pM)\subset F^{p+1}M}$ holds for any $p$). Then, can we perturb ${d_N}$,${F}$,${G}$,${H}$ and ${L}$ so that the above equations ${(4.1)\thicksim(4.3)}$ hold? Markl gave a complete answer to this problem. It is always possible if the ``obstruction class" vanishes (Ideal perturbation lemma). However, we do not need the full generality of Markl's theorem. We only treat the simplest case of these perturbation problems. We treat the case that ${(N,d_N)}$ is a strong deformation retract of ${(M,d_M)}$. Let ${(M,d_M)}$ and ${(N,d_N)}$ be two chain complexes with complete and decreasing filtrations. Assume that chain maps 
\begin{gather*}
F:(M,d_M)\longrightarrow (N,d_N)  \\
G:(N,d_N)\longrightarrow (M,d_M)
\end{gather*}
preserve filtrations and there is a morphism ${H:M\longrightarrow N}$ which satisfies the following conditions:
\begin{gather*}
GF-\textrm{Id}_M=d_MH+Hd_M  \\
FG=\textrm{Id}_N  \\
HH=0, HG=0, FH=0 \ \ (\textrm{annihilation \ properties}).
\end{gather*}
Let ${\widetilde{d_M}=d_M+D_M}$ be a perturbation of ${d_M}$. Assume that $H$ also preserves the filtration. Then the Basic perturbation lemma (\cite{M}) states that there are perturbations ${\widetilde{d_N}}$, $\widetilde{F}$, $\widetilde{G}$ and ${\widetilde{H}}$ which satisfy the following conditions:
\begin{gather*}
\widetilde{G}\widetilde{F}-\textrm{Id}_M=\widetilde{d_M}\widetilde{H}+\widetilde{H}\widetilde{d_M} \\
\widetilde{F}\widetilde{G}=\textrm{Id}_N.
\end{gather*}
So we can perturb ${(N,d_N)}$ so that ${(M,\widetilde{d_M})}$ and ${(N,\widetilde{d_N})}$ are also chain homotopy equivalent. We also have explicit formulas for these perturbations as follows:
\begin{gather*}
\widetilde{d_N}=d_N+F\sum_{l\ge 0}(D_MH)^lD_MG \\
\widetilde{F}=F\sum_{l\ge 0}(D_MH)^l \\
\widetilde{G}=\sum_{l\ge 0}(HD_M)^lG  \\
\widetilde{H}=H\sum_{l\ge 0}(D_MH)^l.
\end{gather*}
\begin{Rem}
We assumed that $M$ and $N$ are complete and $H$ preserves the filtration. Note that
\begin{gather*}
(D_MH)^l(F^p)\subset F^{p+l}, \ (HD_M)^l(F^p)\subset F^{p+l}
\end{gather*}
holds for all $p$ and $l$. This implies that the above power series converge. 
\end{Rem}

Next, we apply the Basic perturbation lemma to construct a strict Floer cochain complex of $H$. In our case, ${(M,d_M)}$, ${(N,d_N)}$ and ${\widetilde{d_M}}$ correspond to the following:
\begin{gather*}
(M,d_M)=(\bigoplus_{x\in P(H,\gamma)}CF^{loc}(H,x)\otimes \Lambda_0, d_F^{loc})  \\
(N,d_N)=(\bigoplus_{x\in P(H,\gamma)}HF^{loc}(H,x)\otimes \Lambda_0, 0)   \\
\widetilde{d_M}=\widetilde{d_F}=d_F^{loc}+D  \\
F=\Pi, H=\Theta .
\end{gather*}
Filtrations on $M$ and $N$ are defined as follows:
\begin{gather*}
F^pM=T^{p\epsilon}M \\
F^pN=T^{p\epsilon}N.
\end{gather*}
The ``inclusion"
\begin{equation*}
G:N\longrightarrow M
\end{equation*}
is defined  as follows:
\begin{gather*}
HF^{loc}(H,x_i)\longrightarrow CF^{loc}(H,x_i)  \\
\bigg[\sum{\alpha_aX_i^a}\bigg]\mapsto \sum{\alpha_aX_i^a}.
\end{gather*}

Note that ${\Theta}$ preserves the energy filtration. Then according to the Basic perturbation lemma and the explicit construction of ${\widetilde{d_N}}$, we can define the coboundary operator ${d_F:CF(H,\gamma:\Lambda_0)\rightarrow CF(H,\gamma:\Lambda_0)}$ by the following formula:
\begin{equation}
d_F=\Pi\circ \sum_{l=0}^\infty (D\Theta)^l D.
\end{equation}
Note that we identify ${HF^{loc}(H,x_i)}$ with ${\textrm{span}_{\mathbb{F}_p}\langle X_i^1,\cdots,X_i^{\textrm{dim}_{\mathbb{F}_p}HF^{loc}(H,x_i)}\rangle}$ in this formula. Then, ${(CF(H,\gamma:\Lambda_0),d_F)}$ is a strict cochain complex which is also cochain homotopy equivalent to the original cochain complex ${(CF(\widetilde{H},\gamma:\Lambda_0),\widetilde{d_F})}$. So we defined the Floer cochain complex ${(CF(H,\gamma:\Lambda_0),d_F)}$ in a homologically canonical way.

Next, we introduce the notion of the local equivariant (Tate) Floer cohomology. Assume that ${x\in P(H^{(p)},\gamma)}$ is an isolated $p$-periodic orbit of ${H^{(p)}}$. Recall that there is a ${\mathbb{Z}_p}$-action on ${P(H^{(p)},\gamma)}$ as follows:

\begin{gather*}
(m\cdot x)(t)=x(t+\frac{m}{p}) \ \ \ \ (m\in \mathbb{Z}_p).
\end{gather*}
The local ${\mathbb{Z}_p}$-equivariant (Tate) Floer cohomology is defined for 
the orbit ${\mathbb{Z}_px}$ of ${x\in P(H^{(p)},\gamma)}$ under the ${\mathbb{Z}_p}$-action. Let ${\widetilde{H}}$ be a small perturbation of ${H}$ such that 
${\widetilde{H}^{(p)}}$ is non-degenerate. Each ${m\cdot x}$ splits into non-degenerate $p$-periodic orbits of ${\widetilde{H}^{(p)}}$. So we have the local Floer cochain complex ${CF^{loc}(H^{(p)},m\cdot x)}$ for each ${m\cdot x}$. The local ${\mathbb{Z}_p}$-equivariant Floer cochain complex of ${\mathbb{Z}_px}$ is defined as follows:

\begin{gather*}
CF_{\mathbb{Z}_p}^{loc}(H^{(p)},\mathbb{Z}_px)=\bigoplus_{y\in \mathbb{Z}_px}CF^{loc}(H^{(p)},y)\otimes \mathbb{F}_p[[u]]\langle \theta \rangle.
\end{gather*}
Let ${U_y\subset S^1\times M}$ be a small isolating neighborhoods of ${y\in \mathbb{Z}_px}$. By counting solutions of the equivariant Floer equation contained in $\cup_{y\in \mathbb{Z}_px}U_y$, we can define the local ${\mathbb{Z}_p}$-equivariant Floer coboundary operator ${d_{eq}^{loc}}$ as follows. Let $\mathcal{P}$ be the set of $1$-periodic orbits of ${\widetilde{H}^{(p)}}$ whose images are contained in ${\cup_{y\in\mathbb{Z}_px}U_y}$. So, each $m\cdot x$ splits into periodic orbits contained in $\mathcal{P}$. Note that $\mathcal{P}$ is the generator of ${\bigoplus_{y\in \mathbb{Z}_px}CF^{loc}(H^{(p)},y)}$. For ${y,z\in \mathcal{P}}$, ${m\in \mathbb{Z}_p}$, ${\alpha \in \{0,1\}}$, ${i\in \mathbb{Z}_{\ge 0}}$, we consider the following equation:
\begin{gather*}
(u,v)\in C^{\infty}(\mathbb{R}\times S^1,M)\times C^{\infty}(\mathbb{R},S^{\infty}) \\
\partial_su(s,t)+J_{v(s,t),t}(\partial_tu(s,t)-X_{\widetilde{H}^{(p)}}(u(s,t)))=0  \\
\frac{d}{ds}v(s)-\mathrm{grad}(\widetilde{F})=0 \\
\lim_{s\rightarrow -\infty}v(s)=Z_{\alpha}^0,\lim_{s\rightarrow +\infty}v(s)=Z_i^m,\lim_{s\rightarrow -\infty}u(s,t)=y(t),\lim_{s\rightarrow +\infty}u(s,t)=z(t-\frac{m}{p})\\
u(s,t)\in \bigcup_{m\in \mathbb{Z}_p}U_{m\cdot x}.
\end{gather*}
We denote the space of solutions of this equation modulo ${\mathbb{R}}$-action by ${\mathcal{M}_{\alpha,i,m}^{loc}(y,z)}$. We define ${d_{\alpha,loc}^{i,m}}$ (${\alpha\in \{0,1\}}$, ${i\in \mathbb{Z}_{\ge0}}$, ${m\in \mathbb{Z}_p}$) as follows:
\begin{gather*}
d_{\alpha,loc}^{i,m}:CF_{\mathbb{Z}_p}^{\mathrm{loc}}(H^{(p)},\mathbb{Z}_px)\longrightarrow CF_{\mathbb{Z}_p}^{loc}(H^{(p)},\mathbb{Z}_px)\\
y\mapsto \sum_{z\in \mathcal{P}}\sharp\mathcal{M}_{\alpha,i,m}^{loc}(y,z)\cdot z.
\end{gather*}
Let ${d_{\alpha,loc}^{i,m}}$ be the sum ${d_{\alpha,loc}^{i}=\sum_{m\in \mathbb{Z}_p}d_{\alpha,loc}^{i,m}}$. Then, ${d_{eq}^{loc}}$ is defined as follows:
\begin{gather*}
d_{eq}^{loc}(y\otimes 1)=\sum_{i=0}^{\infty}d_{0,loc}^{2i}(y)\otimes u^i+\sum_{i=0}^{\infty}d_{0,loc}^{2i+1}(y)\otimes u^i\theta \\
d_{eq}^{loc}(y\otimes \theta)=\sum_{i=0}^{\infty}d_{1,loc}^{2i+1}(y)\otimes u^i\theta+
\sum_{i=1}^{\infty}d_{1,loc}^{2i}(y)\otimes u^i.
\end{gather*}
The local ${\mathbb{Z}_p}$-equivariant Floer cohomology of ${\mathbb{Z}_px}$ is defined as the cohomology of this cochain complex:

\begin{gather*}
HF_{\mathbb{Z}_p}^{loc}(H^{(p)},\mathbb{Z}_px)=H(CF_{\mathbb{Z}_p}^{loc}(H^{(p)},\mathbb{Z}_px),d_{eq}^{loc}).
\end{gather*}
The local ${\mathbb{Z}_p}$-equivariant Tate Floer cochain complex is defined as follows:

\begin{gather*}
\widehat{CF}_{\mathbb{Z}_p}^{loc}(H^{(p)},\mathbb{Z}_px)=CF_{\mathbb{Z}_p}^{loc}(H^{(p)},\mathbb{Z}_px)\otimes \mathbb{F}_p[u^{-1},u]]\langle \theta \rangle.
\end{gather*}
The local ${\mathbb{Z}_p}$-equivariant Tate Floer cohomology is the cohomology of this cochain complex.

\begin{gather*}
\widehat{HF}_{\mathbb{Z}_p}^{loc}(H^{(p)},\mathbb{Z}_px)=H(\widehat{CF}_{\mathbb{Z}_p}^{loc}(H^{(p)},\mathbb{Z}_px),\widehat{d}_{eq}^{loc})
\end{gather*}

Next, assume that ${(M,\omega)}$ is a toroidally monotone symplectic manifold and ${P(H^{(p)},\gamma)}$ is finite. We construct a strict coboundary operator ${\widehat{d}_{eq}}$ on the ${\mathbb{Z}_p}$-equivariant Tate Floer cochain complex
\begin{equation*}
\widehat{CF}_{\mathbb{Z}_p}(H^{(p)},\gamma)=\bigoplus_{\mathbb{Z}_px\in P(H^{(p)},\gamma)/\mathbb{Z}_p}\widehat{HF}_{\mathbb{Z}_p}^{loc}(H^{(p)},\mathbb{Z}_px)\otimes \Lambda_0
\end{equation*}
by applying the Basic perturbation lemma. For every ${\mathbb{Z}_px_i\in P(H^{(p)},\gamma)/\mathbb{Z}_p}$, we choose a basis ${\{U_i^a,V_i^b,W_i^c\}}$ of ${\widehat{CF}_{\mathbb{Z}_p}^{loc}(H^{(p)},\mathbb{Z}_px_i)}$ over ${\mathbb{F}_p[u^{-1},u]]}$ as follows:
\begin{gather*}
\widehat{d}_{eq}^{loc}(U_i^a)=0  \ \ \ (1\le a\le \textrm{dim}_{\mathbb{F}_p[u^{-1},u]]}\widehat{HF}_{\mathbb{Z}_p}^{loc}(H^{(p)},\mathbb{Z}_px_i))\\
\widehat{d}_{eq}^{loc}(V_i^b)=0  \ \ \ (1\le b\le \frac{1}{2}(\widetilde{l_i}-\textrm{dim}_{\mathbb{F}_p[u^{-1},u]]}\widehat{HF}_{\mathbb{Z}_p}^{loc}(H^{(p)},\mathbb{Z}_px_i))\\
\widehat{d}_{eq}^{loc}(W_i^c)=V_i^c \ \ \ (1\le c\le \frac{1}{2}(\widetilde{l_i}-\textrm{dim}_{\mathbb{F}_p[u^{-1},u]]}\widehat{HF}_{\mathbb{Z}_p}^{loc}(H^{(p)},\mathbb{Z}_px_i)).
\end{gather*}
Here, $\widetilde{l_i}$ is the number of perturbed periodic orbits of ${\mathbb{Z}_px_i}$. We define ${\Pi_{i,eq}}$, ${\Theta_{i,eq}}$, ${\Pi_{eq}}$ and ${\Theta_{eq}}$ as follows:
\begin{gather*}
\Pi_{i,eq}\big(\sum_a\alpha_aU_i^a+\sum_b\beta_bV_i^b+\sum_c\gamma_c
W_i^c\big) 
=\sum_a\alpha_aU^a  \\
\Theta_{i,eq}\big(\sum_a\alpha_aU_i^a+\sum_b\beta_bV_i^b+\sum_c\gamma_c
W_i^c\big) 
=\sum_c\gamma_cW_i^c  \\
\Pi_{eq}=\sum_i\Pi_{i,eq}  \\
\Theta_{eq}=\sum_{i}\Theta_{i,eq}.
\end{gather*}
Note that we can identify ${\widehat{CF}_{\mathbb{Z}_p}(\widetilde{H}^{(p)},\gamma)}$ with ${\bigoplus_{\mathbb{Z}_px\in P(H^{(p)},\gamma)/\mathbb{Z}_p}\widehat{CF}_{\mathbb{Z}_p}^{loc}(H^{(p)},\mathbb{Z}_px)\otimes \Lambda_0}$. So, the equivariant Floer coboundary operator $d_{eq}$ on ${\widehat{CF}_{\mathbb{Z}_p}(\widetilde{H}^{(p)},\gamma)}$ induces a coboundary operator on ${\bigoplus_{\mathbb{Z}_px\in P(H^{(p)},\gamma)/\mathbb{Z}_p}\widehat{CF}_{\mathbb{Z}_p}^{loc}(H^{(p)},\mathbb{Z}_px)\otimes \Lambda_0}$. We denote this coboundary operator by ${\overline{d}_{eq}}$. Next, we deform ${\overline{d}_{eq}}$. Every element ${y_i\in P(H^{(p)},\gamma)}$ splits into finitely many non-degenerate periodic orbits ${\{y_i^1,\cdots,y_i^{l_i}\}\subset P(\widetilde{H}^{(p)},\gamma)}$. Let ${v_i^j:[0,1]\times S^1\rightarrow M}$ be a small cylinder which connects $y_i$ and ${y_i^j}$:
\begin{gather*}
v_i^j(0,t)=y_i(t), \ v_i^j(1,t)=y_i^j(t) \\
v_i^j(s,t)\in U_{y_i}.
\end{gather*}
Here ${U_{y_i}}$ is a small open neighborhood of ${y_i}$. We define a constant ${c(y_i,y_i^j)}$ as follows:
\begin{gather*}
c(y_i,y_i^j)=\int_{[0,1]\times S^1}(v_i^j)^*\omega+\int_0^1H^{(p)}(t,y_i(t))-\int_0^1\widetilde{H}^{(p)}(t,y_i^j(t))dt.
\end{gather*}
We define the following correction map ${\kappa_{eq}}$:
\begin{gather*}
\kappa_{eq}:\bigoplus_{\mathbb{Z}_px\in P(H^{(p)},\gamma)/\mathbb{Z}_p}\widehat{CF}_{\mathbb{Z}_p}^{loc}(H^{(p)},\gamma)\otimes \Lambda
\longrightarrow \bigoplus_{\mathbb{Z}_px\in P(H^{(p)},\gamma)/\mathbb{Z}_p}\widehat{CF}_{\mathbb{Z}_p}^{loc}(H^{(p)},\gamma)\otimes \Lambda \\
y_i^j\longmapsto T^{c(y_i,y_i^j)}y_i^j.
\end{gather*}
We define ${\widetilde{d}_{eq}}$ by ${\widetilde{d}_{eq}=\kappa_{eq}^{-1}\circ \overline{d}_{eq}\circ \kappa_{eq}}$. Note that ${\kappa_{eq}}$ is defined over ${\Lambda}$ (not over ${\Lambda_0}$). However, ${\widetilde{d}_{eq}}$ is defined over ${\Lambda_0}$. The proof of this fact is the same as that of Lemma 11. Since this is repetitive, we omit it here.

This $\widetilde{d}_{eq}$ is decomposed as follows:
\begin{gather*}
\widetilde{d}_{eq}=\widehat{d}_{eq}^{loc}+D_{eq}  \\
D_{eq}\Big(\oplus \widehat{CF}_{\mathbb{Z}_p}^{loc}(H^{(p)},\mathbb{Z}_px)\otimes \Lambda_0\Big)\subset T^{\epsilon}\Big(\oplus\widehat{CF}_{\mathbb{Z}_p}^{loc}(H^{(p)},\mathbb{Z}_px)\otimes \Lambda_0 \Big).
\end{gather*}
The existence of such ${\epsilon>0}$ follows from the following argument. Let $\mathcal{S}$ be the following set:
\begin{gather*}
\mathcal{S}=\Bigg\{\int_{[0,1]\times S^1}u^*\omega+\int_0^1H^{(p)}(t,x(t))-H^{(p)}(t,y(t))dt\ \Bigg|\ \begin{matrix}x,y\in P(H^{(p)},\gamma)\\ u:[0,1]\times S^1\rightarrow M\\ u(0,t)=x(t), \ u(1,t)=y(t)\end{matrix}\Bigg\}.
\end{gather*}
This $\mathcal{S}$ is a discrete subset of ${\mathbb{R}}$ because ${(M,\omega)}$ is toroidally monotone. So, it suffices to choose $\epsilon>0$ so that 
\begin{gather*}
\epsilon<\min \{|a|>0 \ | \ a\in \mathcal{S}\}
\end{gather*}
holds.

Then, the coboundary operator ${\widehat{d}_{eq}}$ is defined as follows:
\begin{gather}
\widehat{d}_{eq}:\widehat{CF}_{\mathbb{Z}_p}(H^{(p)},\gamma)\longrightarrow \widehat{CF}_{\mathbb{Z}_p}(H^{(p)},\gamma) \\
\widehat{d}_{eq}= \Pi_{eq}\circ \sum_{l=0}^{\infty}(D_{eq}\Theta_{eq})^lD_{eq}.
\end{gather}
We can see that this sequence converges because ${\Theta_{eq}}$ preserves the filtration and 
\begin{gather*}
\textrm{Im}(\Pi_{eq}(D_{eq}\Theta_{eq})^lD_{eq})\subset T^{(l+1) \epsilon}\cdot \widehat{CF}_{\mathbb{Z}_p}^{loc}(H^{(p)},\gamma)
\end{gather*}
holds (Of course, it also follows from the Basic perturbation lemma. The conditions of the Basic perturbation lemma are ``${\widehat{CF}_{\mathbb{Z}_p}^{loc}(H^{(p)},\gamma)}$ is complete with respect to the filtration" and ``$\Theta_{eq}$ preserves the filtration".).

\section{Proof of the main theorem: toroidally monotone case}

In section 5 and section 6, we prove Theorem 2. We divide the proof into two parts, toroidally monotone cases, and weakly monotone cases. The reason for this is that we have not constructed ${\mathbb{Z}_p}$-equivariant Floer theory for weakly monotone symplectic manifolds yet. As we mentioned before, we have to overcome some technical difficulties in the weakly monotone case. Once we establish ${\mathbb{Z}_p}$-equivariant Floer theory on weakly monotone symplectic manifolds, the rest of the proof is almost the same as in the toroidally monotone case. So the essential part of our proof is given in the toroidally monotone case. 
 
In this section, we prove Theorem 2 for toroidally monotone symplectic manifolds. Let ${(M,\omega)}$ be a closed toroidally monotone symplectic manifold. We fix ${H\in C^{\infty}(S^1\times M)}$ and ${\gamma \neq 0\in H_1(M:\mathbb{Z})/\textrm{Tor}}$. We also assume that ${P(H,\gamma)=\{x_1,\cdots,x_k\}}$ and ${P(H^{(p)},p\gamma)=\{y_1,\cdots,y_k\}}$ and ${y_i=x_i^{(p)}}$ holds. In other words, every $p$-periodic orbit of $H$ in ${p\gamma}$ is not simple. Our purpose is to prove that there is a simple ${p'}$-periodic orbit in ${p\gamma}$. Here ${p'}$ is the smallest prime number greater than $p$.

In the previous section, we defined the ${\mathbb{Z}_p}$-equivariant Tate Floer cochain complex ${(\widehat{CF}(H^{(p)},p\gamma), \widehat{d}_{eq})}$ and the ${\mathbb{Z}_p}$-equivariant Tate Floer cohomology ${\widehat{HF}_{\mathbb{Z}_p}(H^{(p)},p\gamma)}$. Recall that we applied the perturbation theory to make ${(\widehat{CF}(H^{(p)},p\gamma), \widehat{d}_{eq})}$ a strict cochain complex. 

In the definition of ${\widehat{d}_{eq}}$, we have seen that ${\widehat{d}_{eq}}$ is written as ${d_F+\textrm{higer\ terms}}$ if ${H^{(p)}}$ is non-degenerate as follows:
\begin{gather*}
\widehat{d}_{eq}(x\otimes 1)=d_F(x)\otimes 1+x_1\otimes \theta+\sum_{k=1}\sum_{k'=0,1}x_{2k+k'}\otimes u^k\theta^{k'}  \\
\widehat{d}_{eq}(x\otimes \theta)=d_F(x)\otimes \theta+\sum_{k=1}\sum_{k'=0,1}y_{2k+k'}\otimes u^k\theta^{k'}.
\end{gather*}

Next we prove that this is also true when ${H^{(p)}}$ is possibly degenerate. Let ${\widetilde{H}}$ be a perturbation of $H$ such that ${\widetilde{H}^{(p)}}$ is non-degenerate. In this case, each ${y_i\in P(H^{(p)},p\gamma)}$ splits into non-degenerate periodic orbits ${\{y_i^1,\cdots,y_i^{l_i}\}\subset P(\widetilde{H}^{(p)},p\gamma)}$. We prove the following lemma.

\begin{Lem}
Assume that ${y_i=x_i^{(p)}}$ holds and $p$ is an admissible prime number. We can choose a basis ${\{\widetilde{X}_i^a, \widetilde{Y}_i^b, \widetilde{Z}_i^c, \widetilde{X}_{i,\theta}^a, \widetilde{Y}_{i,\theta}^b, \widetilde{Z}_{i,\theta}^c\}}$ of 
\begin{gather*}
\widehat{CF}_{\mathbb{Z}_p}^{loc}(H^{(p)}, y_i)=CF^{loc}(H^{(p)},y_i)\otimes \mathbb{F}_p[u^{-1},u]]\langle \theta \rangle
\end{gather*}
over ${\mathbb{F}_p[u^{-1},u]]}$ and a basis ${\{X_i^a,Y_i^b,Z_i^c\}}$ of ${CF^{loc}(H^{(p)},y_i)}$ which satisfies the following conditions:
\begin{gather*}
\widetilde{X}_i^a=X_i^a\otimes 1+(a_1^{(a)}\otimes \theta+\sum_{k=1}^{\infty}\sum_{k'=0,1}a_{2k+k'}^{(a)}\otimes u^k\theta^{k'})  \\
 \widetilde{Y}_i^b=Y_i^b\otimes 1+(b_1^{(b)}\otimes \theta+\sum_{k=1}^{\infty}\sum_{k'=0,1}b_{2k+k'}^{(b)}\otimes u^k\theta^{k'})  \\
 \widetilde{Z}_i^c=Z_i^c\otimes 1+(c_1^{(c)}\otimes \theta+\sum_{k=1}^{\infty}\sum_{k'=0,1}c_{2k+k'}^{(c)}\otimes u^k\theta^{k'})   \\
 \widetilde{X}_{i,\theta}^a=X_i^a\otimes\theta+\sum_{k=1}^{\infty}\sum_{k'=0,1}d_{2k+k'}^{(a)}\otimes u^k\theta^{k'} \\
 \widetilde{Y}_{i,\theta}^b=Y_i^b\otimes\theta+\sum_{k=1}^{\infty}\sum_{k'=0,1}e_{2k+k'}^{(b)}\otimes u^k\theta^{k'} \\
 \widetilde{Z}_{i,\theta}^c=Z_i^c\otimes\theta+\sum_{k=1}^{\infty}\sum_{k'=0,1}f_{2k+k'}^{(c)}\otimes u^k\theta^{k'} \\
d_F^{loc}(X_i^a)=d_F^{loc}(Y_i^b)=0, \ d_F^{loc}(Z_i^c)=Y_i^c, \ d_{eq}^{loc}(\widetilde{X}_i^a)=d_{eq}^{loc}(\widetilde{Y}_i^b)=0, \ d_{eq}^{loc}( \widetilde{Z}_i^c)=\widetilde{Y}_i^c \\
d_{eq}^{loc}(\widetilde{X}_{i,\theta}^a)=d_{eq}^{loc}(\widetilde{Y}_{i,\theta}^b)=0, \ d_{eq}^{loc}( \widetilde{Z}_{i,\theta}^c)=\widetilde{Y}_{i,\theta}^c \\
1\le a\le \textrm{dim}_{\mathbb{F}_p}HF^{loc}(H^{(p)},y_i), \ 1\le b,c \le \frac{1}{2}(l_i-\textrm{dim}_{\mathbb{F}_p}HF^{loc}(H^{(p)},y_i))   \\
\{a_*^{(a)},b_*^{(b)},c_*^{(c)}, d_*^{(a)},e_*^{(b)},f_*^{(c)}\}\subset CF^{loc}(H^{(p)},y_i).
\end{gather*}
\end{Lem}
\vspace{5mm}
\textbf{proof}(Lemma 13):
Note that ${CF^{loc}(H^{(p)},y_i)=\oplus_j \mathbb{F}_p\cdot y_i^j}$ holds. The existence of a basis ${\{X_i^a,Y_i^b,Z_i^c\}}$ is trivial and we can define ${ \widetilde{Y}_i^b, \widetilde{Z}_i^c, \widetilde{Y}_{i,\theta}^b,\widetilde{Z}_{i,\theta}^c}$ by ${\widetilde{Z}_i^c=Z_i^c\otimes 1}$, ${\widetilde{Z}_{i,\theta}^c=Z_i^c\otimes \theta}$, ${\widetilde{Y}_i^b=d_{eq}^{loc}(\widetilde{Z}_i^b)}$ and ${\widetilde{Y}_{i,\theta}^b=d_{eq}^{loc}(\widetilde{Z}_{i,\theta}^b)}$. So what we have to prove is the existence of ${\widetilde{X}_i^a}$ and ${\widetilde{X}_{i,\theta}^a}$.

Let ${\mathcal{C}}$ be the set of all cocycles in ${(\widehat{CF}_{\mathbb{Z}_p}^{loc}(H^{(p)},y_i),d_{eq}^{loc})}$. We consider two projections 
\begin{equation*}
\pi_1,\pi_2:\mathcal{C} \longrightarrow \textrm{span}_{\mathbb{F}_p} \langle X_i^1,\cdots,X_i^{d}\rangle \ \ \ \  (d=\textrm{dim}_{\mathbb{F}_p}HF^{loc}(H^{(p)},y_i)) .
\end{equation*}
For ${z=(a_k\otimes 1+b_k\otimes \theta)u^k+\sum_{l=k+1}(a_l\otimes 1+b_l\otimes \theta)u^l}$, we define ${\pi_1(z)}$ and ${\pi_2(z)}$ by the following formula. Let ${\Pi_i:CF(H^{(p)},y_i)\rightarrow  \textrm{span}_{\mathbb{F}_p} \langle X_i^1,\cdots,X_i^{d}\rangle}$ be the projection as before:
\begin{gather*}
\Pi_i(\sum \alpha_aX_i^a+\sum \beta_bY_i^b+\sum \gamma_cZ_i^c)=\sum \alpha_aX_i^a.
\end{gather*}
Then we define
\begin{gather*}
\pi_1(z)=\Pi_i(a_k) \\
\pi_2(z)=\begin{cases}  \Pi_i(b_k)  & a_k=0   \\   0& a_k\neq 0. \end{cases}
\end{gather*}
Note that ${\pi_1}$ and ${\pi_2}$ are not linear maps. However, their images are vector subspaces of ${\textrm{span}_{\mathbb{F}_p} \langle X_i^1,\cdots,X_i^{d}\rangle}$. This follows from the following arguments. Let ${c_1}$ and ${c_2}$ be elements of ${\textrm{Im}(\pi_1)}$. Then we can choose cocycles ${\{\widetilde{z}^{(1)},\widetilde{z}^{(2)}\}\subset \mathcal{C}}$ as follows:

\begin{gather*}
\widetilde{z}^{(1)}=a_0^{(1)}\otimes u^{k_1}+b_0^{(1)}\otimes u^{k_1}\theta+\sum_{l\ge k_1+1}(a_l^{(1)}\otimes 1+b_l^{(1)}\otimes \theta)u^l  \\
\widetilde{z}^{(2)}=a_0^{(2)}\otimes u^{k_2}+b_0^{(2)}\otimes u^{k_2}\theta+\sum_{l\ge k_2+1}(a_l^{(2)}\otimes 1+b_l^{(2)}\otimes \theta)u^l  \\
\pi_1(\widetilde{z}^{(1)})=\Pi_i(a_0^{(1)})=c_1 \\
\pi_1(\widetilde{z}^{(2)})=\Pi_i(a_0^{(2)})=c_2.
\end{gather*}
Note that ${z^{(1)}=u^{-k_1}\widetilde{z}^{(1)}}$ and ${z^{(2)}=u^{-k_2}\widetilde{z}^{(2)}}$ are also cocycles such that
\begin{gather*}
z^{(1)}=a_0^{(1)}\otimes 1+b_0^{(1)}\otimes \theta+\sum_{l\ge 1}(a_{l+k_1}^{(1)}\otimes 1+b_{l+k_1}^{(1)}\otimes \theta)u^l  \\
z^{(2)}=a_0^{(2)}\otimes 1+b_0^{(2)}\otimes \theta+\sum_{l\ge 1}(a_{l+k_2}^{(2)}\otimes 1+b_{l+k_2}^{(2)}\otimes \theta)u^l  \\
\pi_1(z^{(1)})=\Pi_i(a_0^{(1)})=c_1 \\
\pi_1(z^{(2)})=\Pi_i(a_0^{(2)})=c_2
\end{gather*}
holds. Then, for any ${\alpha_1}$, ${\alpha_2\in \mathbb{F}_p}$, 
\begin{gather*}
\pi_1(\alpha_1z^{(1)}+\alpha_2z^{(2)})=\Pi_i(\alpha_1a_0^{(1)}+\alpha_2a_0^{(2)})=\alpha_1c_1+\alpha_2c_2
\end{gather*}
holds if ${\alpha_1a_0^{(1)}+\alpha_2a_0^{(2)}\neq 0}$ holds. This implies that ${\textrm{Im}(\pi_1)}$ is a vector subspace. A similar argument can be used to prove that ${\textrm{Im}(\pi_2)}$ is a vector subspace.

Assume that ${\textrm{Im}(\pi_1)}$ is generated by ${\{\widetilde{V_j}\}}$ and ${\textrm{Im}(\pi_2)}$ is generated by ${\{\widetilde{W_j}\}}$:

\begin{gather*}
\textrm{Im}(\pi_1)=\textrm{span}_{\mathbb{F}_p}\langle \widetilde{V_1},\cdots,\widetilde{V_{\alpha}}\rangle  \\
\textrm{Im}(\pi_2)=\textrm{span}_{\mathbb{F}_p}\langle \widetilde{W_1},\cdots,\widetilde{W_{\beta}}\rangle.
\end{gather*}
We choose ${\{V_1,\cdots,V_{\alpha},W_1,\cdots,W_{\beta}\}\subset \mathcal{C}}$ such that ${\pi_1(V_j)=\widetilde{V_j}}$ and ${\pi_2(W_j)=\widetilde{W_j}}$ hold. Next we prove that ${\{[V_1],\cdots,[V_{\alpha}],[W_1],\cdots,[W_{\beta}]\}}$ generates ${\widehat{HF}_{\mathbb{Z}_p}^{loc}(H^{(p)},y_i)}$ over ${\mathbb{F}_p[u^{-1},u]]}$. We fix ${z\in \mathcal{C}}$ as follows:
\begin{gather*}
z=\sum_{l\ge m}(a_l\otimes 1+b_l\otimes \theta)u^l.
\end{gather*}
It suffices to find ${z_m=(c_m\otimes 1+d_m\otimes \theta)u^m}$ and ${\zeta_m^{(i)}}$, ${\eta_m^{(j)}\in \mathbb{F}_p}$ such that
\begin{gather*}
z-\widehat{d}_{eq}^{loc}(z_m)-\sum_{1\le i\le \alpha}(\zeta_m^{(i)}u^m)V_{i}-\sum_{1\le j\le \beta}(\eta_m^{(j)}u^m)W_{j}=\sum_{l\ge m+1}(a'_l\otimes 1+b'_l\otimes \theta)u^l
\end{gather*}
holds. If this is true, we can construct $z'$ and ${\zeta^{(i)},\eta^{(j)}\in \mathbb{F}_p[u^{-1},u]]}$ so that 
\begin{gather*}
z-\widehat{d}_{eq}^{loc}(z')=\sum_{1\le i\le \alpha}\zeta^{(i)}V_i+\sum_{1\le j\le \beta}\eta^{(j)}W_j
\end{gather*}
holds. We choose ${c_m\in CF^{loc}(H^{(p)},y_i)}$ so that 
\begin{gather*}
a_m-d_F^{loc}(c_m)\in \textrm{span}_{\mathbb{F}_p}\langle X_i^1,\cdots,X_i^d \rangle
\end{gather*}
holds. Note that ${a_m-d_F^{loc}(c_m)\in \textrm{span}_{\mathbb{F}_p}\langle \widetilde{V_1},\cdots,\widetilde{V_{\alpha}} \rangle}$ holds because 
\begin{gather*}
z-\widehat{d}_{eq}^{loc}(c_m\otimes u^m)=\big((a_m-d_F^{loc}(c_m))\otimes 1+\widetilde{b}_m\otimes \theta\big)u^m+\sum_{l\ge m+1}(\widetilde{a}_l\otimes 1+\widetilde{b}_l\otimes \theta)u^l
\end{gather*}
implies that ${\pi_1(z-\widehat{d}_{eq}^{loc}(c_m\otimes u^m))=a_m-d_F^{loc}(c_m)}$ holds. So we can choose ${\{\zeta_m^{(i)}\}_{i\le \alpha}}$ so that
\begin{gather*}
a_m-d_F^{loc}(c_m)=\sum_{1\le i\le \alpha}\zeta_m^{(i)}\widetilde{V_i}
\end{gather*}
holds. This implies that
\begin{gather*}
z-\widehat{d}_{eq}^{loc}(c_m\otimes u^m)-\sum_{1\le i\le \alpha}(\zeta_m^{(i)}u^m)V_i=b''_m\otimes u^m\theta+\sum_{l\ge m+1}(a''_l\otimes 1+b''_l\otimes \theta)u^l
\end{gather*}
holds. A similar argument can be used to find ${d_m\in CF^{loc}(H^{(p)},y_i)}$ and ${\{\eta_m^{(j)}\}_{j\le \beta}}$ such that 
\begin{gather*}
z_m=(c_m\otimes 1+d_m\otimes \theta)u^m \\
z-\widehat{d}_{eq}^{loc}(z_m)-\sum_{1\le i\le \alpha}(\zeta_m^{(i)}u^m)V_{i}-\sum_{1\le j\le \beta}(\eta_m^{(j)}u^m)W_{j}=\sum_{l\ge m+1}(a'_l\otimes 1+b'_l\otimes \theta)u^l
\end{gather*}
holds. So we proved that ${\widehat{HF}_{\mathbb{Z}_p}^{loc}(H^{(p)},y_i)}$ is generated by ${\{V_1,\cdots,V_{\alpha},W_1,\cdots,W_{\beta}\}}$.

Note that there is an isomoprhism between the local Floer cohomology and the local $\mathbb{Z}_p $-equivariant Tate Floer cohomology (\cite{SZ}).

\begin{equation*}
\widehat{HF}_{\mathbb{Z}_p}^{loc}(H^{(p)},y_i)\cong HF^{loc}(H^{(p)},y_i)\otimes \mathbb{F}_p[u^{-1},u]]\langle \theta \rangle
\end{equation*}  
Here we use the fact that $p$ is admissible and ${HF^{loc}(H^{(p)},y_i)\cong HF^{loc}(H,x_i)}$ holds. This implies that ${\alpha=\beta=\textrm{dim}_{\mathbb{F}_p}HF^{loc}(H^{(p)},y_i)}$ and ${\textrm{Im}(\pi_i)=\textrm{span}_{\mathbb{F}_p}\langle X_1,\cdots,X_d\rangle}$. So we can choose ${\widetilde{X}_i^a}$ and ${\widetilde{X}_{i,\theta}^a}$ and we proved Lemma 13.
\begin{flushright}
    $\Box$
\end{flushright}

\begin{Cor}
Assume that ${P(H^{(p)},p\gamma)}$ consists of finitely many non-simple periodic orbits. We fix a basis ${\{X_i^a,Y_i^b,Z_i^c\}}$ of ${CF^{loc}(H^{(p)},y_i)}$ and a basis ${\{\widetilde{X}_i^a,\widetilde{Y}_i^b,\widetilde{Z}_i^c,\widetilde{X}_{i,\theta}^a,\widetilde{Y}_{i,\theta}^b,\widetilde{Z}_{i,\theta}^c\}}$ of ${\widehat{CF}_{\mathbb{Z}_p}^{loc}(H^{(p)},y_i)}$ so that they satisfy the properties ensured by Lemma 13. Let ${(CF(H^{(p)},p\gamma),d_F)}$ be the Floer cochain complex as in ${(4.4)}$ and let ${\widehat{CF}_{\mathbb{Z}_p}^{loc}(H^{(p)},p\gamma)}$ be the ${\mathbb{Z}_p}$-equivariant Tate Floer cochain complex of ${H^{(p)}}$ as in ${(4.5)}$. Let ${\iota}$ and ${\iota_{\theta}}$ be the following injections:
\begin{gather*}
\iota:\bigoplus_iCF^{loc}(H^{(p)},y_i)\otimes \Lambda_0 \longrightarrow \bigoplus_i\widehat{CF}_{\mathbb{Z}_p}^{loc}(H^{(p)},y_i)\otimes \Lambda_0 \\
\sum_{i,a}\alpha_{i,a}X_i^a+\sum_{i,b}\beta_{i,b}Y_i^b+\sum_{i,c}\gamma_{i,c}Z_i^c\mapsto \sum_{i,a}\alpha_{i,a}\widetilde{X}_{i}^a+\sum_{i,b}\beta_{i,c}\widetilde{Y}_i^b+\sum_{i,c}\gamma_{i,c}\widetilde{Z}_i^c
\end{gather*}
\begin{gather*}
\iota_{\theta}:\bigoplus_iCF^{loc}(H^{(p)},y_i)\otimes \Lambda_0 \longrightarrow \bigoplus_i\widehat{CF}_{\mathbb{Z}_p}^{loc}(H^{(p)},y_i)\otimes \Lambda_0 \\
\sum_{i,a}\alpha_{i,a}X_i^a+\sum_{i,b}\beta_{i,b}Y_i^b+\sum_{i,c}\gamma_{i,c}Z_i^c\mapsto \sum_{i,a}\alpha_{i,a}\widetilde{X}_{i,\theta}^a+\sum_{i,b}\beta_{i,c}\widetilde{Y}_{i,\theta}^b+\sum_{i,c}\gamma_{i,c}\widetilde{Z}_{i,\theta}^c \\ (\alpha_{i,\alpha}, \beta_{i,b}, \gamma_{i,c}\in \Lambda_0).
\end{gather*}
Then for every ${\{\alpha_{i,a}\}\subset \Lambda_0,}$
\begin{gather}
\widehat{d}_{eq}\big(\sum_{i,a}\alpha_{i,a}\widetilde{X}_i^a\big)=\iota\big(d_F(\sum_{i,a}\alpha_{i,a}X_i^a)\big)\otimes 1+\sum_{i,a}\beta_{i,a}\widetilde{X}_{i,\theta}^a+\sum_{l\ge 1}\sum_{i,a}\big(\delta_{i,a}\widetilde{X}_i^a+\beta_{i,a}\widetilde{X}_{i,\theta}^a\big)u^l  \\
\widehat{d}_{eq}\big(\sum_{i,a}\alpha_{i,a}\widetilde{X}_{i,\theta}^a\big)=\iota_{\theta}\big(d_F(\sum_{i,a}\alpha_{i,a}X_{i}^a)\big)+\sum_{l\ge 1}\sum_{i,a}\big(\delta_{i,a}^{\theta}\widetilde{X}_i^a+\beta_{i,a}^{\theta}\widetilde{X}_{i,\theta}^a\big)u^l
\end{gather}
holds for some ${\{\beta_{i,a}\}}$, ${\{\delta_{i,a}\}}$, ${\{\beta_{i,a}^{\theta}\}}$ ${\{\delta_{i,a}^{\theta}\}\subset \Lambda_0}$. In this sense, ${\widehat{d}_{eq}}$ is written as ${d_F+\textrm{higher \ terms}}$.
\end{Cor}
\vspace{5mm}
\textbf{proof}:
We can define ${\Pi_{eq}}$ and ${\Theta_{eq}}$ as follows:
\begin{gather*}
\Pi_{eq},\Theta_{eq}:\widehat{CF}_{\mathbb{Z}_p}(H^{(p)},\gamma)\longrightarrow  \widehat{CF}_{\mathbb{Z}_p}(H^{(p)},\gamma)  \\
\Pi_{i,eq}\big(\sum_a\alpha_a\widetilde{X}_i^a+\sum_b\beta_b\widetilde{Y}_i^b+\sum_c\gamma_c
\widetilde{Z}_i^c+\sum_a\alpha_{a,\theta}\widetilde{X}_{i,\theta}^a+\sum_b\beta_{b,\theta}\widetilde{Y}_{i,\theta}^b+\sum_c\gamma_{c,\theta}
\widetilde{Z}_{i,\theta}^c\big) \\
=\sum_a\alpha_a\widetilde{X}_i^a+\sum_a\alpha_{a,\theta}\widetilde{X}_{i,\theta}^a  \\
\Theta_{i,eq}\big(\sum_a\alpha_a\widetilde{X}_i^a+\sum_b\beta_b\widetilde{Y}_i^b+\sum_c\gamma_c
\widetilde{Z}_i^c+\sum_a\alpha_{a,\theta}\widetilde{X}_{i,\theta}^a+\sum_b\beta_{b,\theta}\widetilde{Y}_{i,\theta}^b+\sum_c\gamma_{c,\theta}
\widetilde{Z}_{i,\theta}^c\big) \\
=\sum_b\beta_b\widetilde{Z}_i^c+\sum_b\beta_{b,\theta}
\widetilde{Z}_{i,\theta}^c  \\
\Pi_{eq}=\sum_{i}\Pi_{i,eq}  \\
\Theta_{eq}=\sum_{i}\Theta_{i,eq}.
\end{gather*}
We abbreviate ${\{X_i^a,Y_i^b,Z_i^c\}}$ to ${\{W_j\}}$ and ${\{\widetilde{X}_i^a,\widetilde{Y}_i^b,\widetilde{Z}_i^c,\widetilde{X}_{i,\theta}^a,\widetilde{Y}_{i,\theta}^b,\widetilde{Z}_{i,\theta}^c\}}$ to ${\{\widetilde{W}_j,\widetilde{W}_{j,\theta}\}}$.
Let $\Gamma$ be one of the maps ${\{\Pi, D, \Theta \}}$ in ${(4.4)}$. Let ${\Gamma_{eq}}$ be the corresponding map in ${\{\Pi_{eq},D_{eq},\Theta_{eq}\}}$ in ${(4.5)}$. Recall that we assumed that ${\{W_i\}}$ and ${\{\widetilde{W}_{j},\widetilde{W}_{j,\theta}\}}$ satisfy the properties ensured by Lemma 13. So ${\Gamma_{eq}}$ has the following form:

\begin{gather*}
\Gamma_{eq}\big(\sum_{j}\alpha_{j}\widetilde{W}_j\big)
=\iota\big(\Gamma(\sum_{j}\alpha_{j}W_j)\big)+\sum_{j}\beta_{j}\widetilde{W}_{j,\theta}+\sum_{l\ge 1}\sum_{j}\big(\delta_{j,l}\widetilde{W}_j+\beta_{j,l}\widetilde{W}_{j,\theta}\big)u^l \\
\Gamma_{eq}\big(\sum_{j}\alpha_{j}\widetilde{W}_{j,\theta}\big)=\iota_{\theta}\big(\Gamma(\sum_{j}\alpha_{j}W_{j,\theta})\big)+\sum_{l\ge 1}\sum_{j}\big(\delta_{i,l}\widetilde{W}_j+\beta_{j,l}\widetilde{W}_{j,\theta}\big)u^l.
\end{gather*}

So, ${\widehat{d}_{eq}}$ can be written as ${(5.1)}$ and ${(5.2)}$ because ${d_F}$ is generated by ${\{\Pi, D, \Theta\}}$ as in ${(4.4)}$ and ${\widehat{d}_{eq}}$ is generated by ${\{\Pi_{eq},D_{eq},\Theta_{eq}\}}$ as in ${(4.5)}$.
\begin{flushright}
    $\Box$
\end{flushright}

Our next purpose is to calculate ${\widehat{HF}_{\mathbb{Z}_p}(H^{(p)},p\gamma)}$. For this purpose, we introduce the ${\mathbb{Z}_p}$-equivariant (Tate) cohomology for ${(H,\gamma)}$ (see \cite{SZ}). The ${\mathbb{Z}_p}$-equivariant  cochain complex is defined as follows:

\begin{equation*}
C(\mathbb{Z}_p,CF(H,\gamma:\Lambda_0)^{\otimes p})=CF(H,\gamma:\Lambda_0)^{\otimes p}\otimes \Lambda_0[[u]]\langle \theta \rangle.
\end{equation*}

The Floer differential ${d_F}$ on ${CF(H,\gamma:\Lambda_0)}$ naturally extends to a differential ${d_F^{(p)}}$ on ${CF(H,\gamma:\Lambda_0)^{\otimes p}}$. There is a natural ${\mathbb{Z}_p}$ action $\tau$ on ${CF(H,\gamma:\Lambda_0)^{\otimes p}}$:

\begin{equation*}
\tau(x_0\otimes x_1\otimes \cdots \otimes x_{p-1})=(-1)^{|x_{p-1}|(|x_{0}|+\cdots+|x_{p-2}|)}x_{p-1}\otimes x_0\otimes \cdots \otimes x_{p-2}.
\end{equation*}

Let ${N}$ be the sum ${N=1+\tau+\tau^2\cdots +\tau^{p-1}}$. Then the ${\mathbb{Z}_p}$-equivariant coboundary operator ${d_{\mathbb{Z}_p}}$ is a ${\Lambda_0[[u]]}$-linear map defined as follows:

\begin{gather*}
d_{\mathbb{Z}_p}(x\otimes 1)=d_F^{(p)}(x)\otimes 1+(1-\tau)(x)\otimes \theta  \\
d_{\mathbb{Z}_p}(x\otimes \theta)=-d_F^{(p)}(x)\otimes \theta+N(x)\otimes u.
\end{gather*}
The ${\mathbb{Z}_p}$-equivariant cohomology is defined by the cohomology of this complex.

\begin{equation*}
H(\mathbb{Z}_p,CF(H,\gamma:\Lambda_0)^{\otimes p})=H(C(\mathbb{Z}_p,CF(H,\gamma:\Lambda_0)^{\otimes p}),d_{\mathbb{Z}_p})
\end{equation*}
The ${\mathbb{Z}_p}$-equivariant Tate cochain complex is a coefficient extension of the ${\mathbb{Z}_p}$-equivariant cochain complex as follows:

\begin{gather*}
\widehat{C}(\mathbb{Z}_p,CF(H,\gamma:\Lambda_0)^{\otimes p})=CF(H,\gamma:\Lambda_0)^{\otimes p}\otimes \Lambda_0[u^{-1},u]]\langle \theta \rangle  \\
\widehat{d}_{\mathbb{Z}_p}:\widehat{C}(\mathbb{Z}_p,CF(H,\gamma:\Lambda_0)^{\otimes p}) \longrightarrow \widehat{C}(\mathbb{Z}_p,CF(H,\gamma:\Lambda_0)^{\otimes p}).
\end{gather*}
Here, ${\widehat{d}_{\mathbb{Z}_p}}$ is the natural extension of ${d_{\mathbb{Z}_p}}$. The ${\mathbb{Z}_p}$-equivariant Tate cohomology is the cohomology of this cochain complex as follows:

\begin{gather*}
\widehat{H}(\mathbb{Z}_p,CF(H,\gamma:\Lambda_0)^{\otimes p})=H(\widehat{C}(\mathbb{Z}_p,CF(H,\gamma:\Lambda_0)^{\otimes p}),\widehat{d}_{\mathbb{Z}_p}).
\end{gather*}
The module ${\widehat{H}(\mathbb{Z}_p,CF(H,\gamma:\Lambda_0)^{\otimes p})}$ is determined by ${HF(H,\gamma:\Lambda_0)}$ from the following lemma.

\begin{Lem}[\cite{Sh}]
There is a so-called quasi-Frobenius isomorphism as follows:
\begin{equation*}
r_p^*HF(H,\gamma:\Lambda_0)\otimes \Lambda_0[u^{-1},u]]\langle \theta \rangle \cong \widehat{H}(\mathbb{Z}_p,CF(H,\gamma:\Lambda_0)^{\otimes p}).
\end{equation*}
Here, ${r_p:\Lambda_0\rightarrow \Lambda_0}$ is the homomorphism defined by ${T\rightarrow T^{\frac{1}{p}}}$.
\end{Lem}

Assume that there is an isomorphism
\begin{equation*}
HF(H,\gamma:\Lambda_0)\cong \bigoplus_{i=1}^m\Lambda_0/T^{\beta_j}\Lambda_0 \ \ \ (\beta_j>0).
\end{equation*}
Then, Lemma 15 implies that there is the following isomorphism:
\begin{equation*}
\widehat{H}(\mathbb{Z}_p,CF(H,\gamma:\Lambda_0)^{\otimes p})\cong \Big(\bigoplus_{i=1}^m\Lambda_0/T^{p\beta_j}\Lambda_0\Big)\otimes \Lambda_0[u^{-1},u]]\langle \theta \rangle.
\end{equation*}
So, the module structure of the ${\mathbb{Z}_p}$-equivariant Tate cohomology is completely determined by ${HF(H,\gamma:\Lambda_0)}$.

There is another important operation we consider, the so-called ${\mathbb{Z}_p}$-equivariant pair of pants product:

\begin{gather*}
\mathcal{P}:H(\mathbb{Z}_p,CF(H,\gamma:\Lambda_0)^{\otimes p})\longrightarrow HF_{\mathbb{Z}_p}(H^{(p)},p\gamma)  \\
\widehat{\mathcal{P}}:\widehat{H}(\mathbb{Z}_p,CF(H,\gamma:\Lambda_0)^{\otimes p})\longrightarrow \widehat{HF}_{\mathbb{Z}_p}(H^{(p)},p\gamma).
\end{gather*}
The construction of ${\mathcal{P}}$ (and ${\widehat{\mathcal{P}}}$) is given by counting the solutions of Floer equations on a $p$-branched cover of the cylinder ${\mathbb{R}\times S^1}$ ($p$-legged pants) parametrized by ${S^{\infty}}$ (see section $8$ in \cite{SZ}, see also \cite{Se} for the ${\mathbb{Z}_2}$-equivariant case). We give a detailed construction of $\mathcal{P}$ and ${\widehat{\mathcal{P}}}$ in the last section for weakly monotone case. The important point is that ${\widehat{\mathcal{P}}}$ gives a local isomorphism between the corresponding local cohomologies in the following sense (\cite{SZ}).

We define an action filtration on ${\widehat{C}(\mathbb{Z}_p,CF(H,\gamma:\Lambda_0)^{\otimes p})}$ and ${\widehat{CF}_{\mathbb{Z}_p}(H^{(p)},p\gamma)}$. We fix a sufficiently small positive real number ${\epsilon>0}$. We define the filtration ${F^q\widehat{C}(\mathbb{Z}_p,CF(H,\gamma:\Lambda_0)^{\otimes p})}$ and ${F^q\widehat{CF}_{\mathbb{Z}_p}(H^{(p)},p\gamma)}$ ${(q\in \mathbb{Z}_{\ge 0})}$ as follows:

\begin{gather*}
F^q\widehat{C}(\mathbb{Z}_p,CF(H,\gamma:\Lambda_0)^{\otimes p})=T^{q\epsilon}\widehat{C}(\mathbb{Z}_p,CF(H,\gamma:\Lambda_0)^{\otimes p}) \\
F^q\widehat{CF}_{\mathbb{Z}_p}(H^{(p)},p\gamma)=T^{q\epsilon}\widehat{CF}_{\mathbb{Z}_p}(H^{(p)},p\gamma).
\end{gather*}
Let ${\widehat{d}_{\mathbb{Z}_p}^{loc}}$ be a differential map defined as follows:
\begin{gather*}
\widehat{d}_{\mathbb{Z}_p}^{loc}:\widehat{CF}_{\mathbb{Z}_p}(H^{(p)},\gamma)\longrightarrow   \widehat{CF}_{\mathbb{Z}_p}(H^{(p)},\gamma)  \\
\widehat{d}_{\mathbb{Z}_p}^{loc}(x\otimes 1)=(d_F^{loc})^{(p)}(x)\otimes 1+(1-\tau)\otimes \theta  \\
\widehat{d}_{\mathbb{Z}_p}^{loc}(x\otimes \theta)=-(d_F^{loc})^{(p)}(x)\otimes \theta+N(x)\otimes u.
\end{gather*}
If we divide ${\widehat{d}_{\mathbb{Z}_p}}$ into ${\widehat{d}_{\mathbb{Z}_p}^{loc}+D_{\mathbb{Z}_p}}$, 
\begin{equation*}
D_{\mathbb{Z}_p}(F^q\widehat{C}(\mathbb{Z}_p,CF(H,\gamma:\Lambda_0)^{\otimes p}))\subset F^{q+1}\widehat{C}(\mathbb{Z}_p,CF(H,\gamma:\Lambda_0)^{\otimes p})
\end{equation*}
holds for all ${q\in \mathbb{Z}_{\ge 0}}$. So the ${E_1}$-terms of the associated spectral sequences of ${F^q\widehat{C}}$ and ${F^q\widehat{CF}_{\mathbb{Z}_p}}$ are given by local cohomologies. The $E_1$-term of the spectral sequence associated with the filtration ${F^q\widehat{C}(\mathbb{Z}_p,CF(H,\gamma:\Lambda_0)^{\otimes p})}$ is isomorphic to a ${\Lambda/T^{\epsilon}\Lambda_0}$-module
\begin{gather*}
\bigoplus_{x\in P(H,\gamma)}\widehat{H}(\mathbb{Z}_p,HF^{loc}(H,x)^{\otimes p})\otimes \Lambda_0/T^{\epsilon}\Lambda_0
\end{gather*}
and the $E_1$-term of the spectral sequence associated with the filtration ${F^q\widehat{CF}_{\mathbb{Z}_p}(H^{(p)},p\gamma)}$ is isomorphic to a ${\Lambda_0/T^{\epsilon}\Lambda_0}$-module
\begin{gather*}
\bigoplus_{y\in P(H^{(p)},p\gamma)/\mathbb{Z}_p}\widehat{HF}_{\mathbb{Z}_p}^{loc}(H^{(p)},\mathbb{Z}_py)\otimes \Lambda_0/T^{\epsilon}\Lambda_0 .
\end{gather*}
This fact and the local isomorphism theorem of ${\widehat{\mathcal{P}}}$ imply that $\widehat{\mathcal{P}}$ gives an isomorphism between the ${E_1}$-pages of their spectral sequences. So, to prove that ${\widehat{\mathcal{P}}}$ is an isomorphism, it suffices to prove the following strong convergences of each spectral sequence (Theorem 3.2 in Chapter 15 in \cite{CE}).

\begin{Lem}[strong convergence]
Let ${(A,d)}$ be a filtered differential module such that ${(A,d)=(\widehat{C}(\mathbb{Z}_p,CF(H,\gamma:\Lambda_0)^{\otimes p},\widehat{d}_{\mathbb{Z}_p})}$ or ${(A,d)=(\widehat{CF}_{\mathbb{Z}_p}(H^{(p)},p\gamma:\Lambda_0),\widehat{d}_{eq})}$ holds. 
\begin{enumerate}
\item (weak convergence) For each ${q\in \mathbb{Z}_{\ge 0}}$, the intersection of the images of the homomorphisms
\begin{gather*}
H(F^qA/F^{q+r}A)\longrightarrow H(F^{q+1}A) \ \ \ r\ge 1  \\
[z]\mapsto [dz]
\end{gather*}
is zero.
\item The natural map
\begin{equation*}
u:H(A)\longrightarrow \varprojlim H(A)/F^qH(A)
\end{equation*}
is an isomorphism. Here, ${F^qH(A)}$ is the image of 
\begin{equation*}
H(F^qA)\longrightarrow H(A).
\end{equation*}
\end{enumerate}
\end{Lem}
\vspace{5mm}
\textbf{proof}(Lemma 16):
${(1)}$ Let $K$ be another Hamiltonian function and let ${(B,d)}$ be either the Tate chain complex ${(\widehat{C}(\mathbb{Z}_p,CF(K,\gamma:\Lambda_0)^{\otimes p}),\widehat{d}_{\mathbb{Z}_p})}$ or ${(\widehat{CF}_{\mathbb{Z}_p}(K^{(p)},p\gamma:\Lambda_0),\widehat{d}_{eq})}$ respectively. Then, there are  continuation chain maps over ${\Lambda}$
\begin{gather*}
F:(A\otimes \Lambda,d)\longrightarrow (B\otimes \Lambda, d)  \\
G:(B\otimes \Lambda,d)\longrightarrow (A\otimes \Lambda, d)
\end{gather*}
and a chain homotopy ${\mathcal{H}}$ between ${\textrm{Id}}$ and ${GF}$ as follows. Let ${D}$ be the constant 
\begin{equation*}
D=p\times ||H-K||=p\times \int_0^1\Big\{ \max (H_t-K_t)-\min (H_t-K_t)\Big\}dt .
\end{equation*}
Then, ${\mathcal{H}(A^{\alpha})\subset A^{\alpha -D}}$ holds for all ${\alpha \in \mathbb{R}}$ (${A^{\alpha}=T^{\alpha}A}$). 
If $K$ is a ${C^{\infty}}$-small function, then ${P(K,\gamma)=P(K^{(p)},p\gamma)=\emptyset}$ and ${(B,d)=(0,0)}$. This implies that there is a map 
\begin{equation*}
\mathcal{L}:A\otimes \Lambda \rightarrow A\otimes \Lambda
\end{equation*}
such that 
\begin{gather*}
Id_{A\otimes \Lambda}=d\mathcal{L}+\mathcal{L}d  \\
\mathcal{L}(A^{\alpha})\subset A^{\alpha-p||H||} \ \ \  (\forall \alpha \in \mathbb{R})
\end{gather*}
holds.
This implies that for ${r>\frac{p||H||}{\epsilon}+1}$ and any cocycle ${C\in F^{q+r}A}$, ${C'=\mathcal{L}(C)}$ satisfies ${C'\in F^{q+1}A}$ and ${d(C')=C}$. So the map
\begin{gather*}
H(F^qA/F^{q+r}A)\longrightarrow H(F^{q+1}A) \\
[z]\mapsto [dz]
\end{gather*}
is zero for ${r>\frac{p||H||}{\epsilon}+1}$ because $dz$ is a cocycle in ${F^{q+r}A}$ and hence we can find ${C'\in F^{q+1}A}$ such that ${d(C')=dz}$ holds. So we proved ${(1)}$. \\

${(2)}$ Let ${z\in F^qA}$ be a cocycle for ${q>\frac{p||H||}{\epsilon}}$. Then, the above argument implies that we can choose ${z'\in A}$ so that ${d(z')=z}$ holds. So, ${F^qH(A)}$ is zero for ${q>\frac{p||H||}{\epsilon}}$ and hence
\begin{gather*}
H(A)/F^{q}H(A)=H(A)
\end{gather*}
holds. In particular, the natural map 
\begin{gather*}
u:H(A)\longrightarrow \varprojlim H(A)/F^qH(A)
\end{gather*}
is an isomorphism.
\begin{flushright}
    $\Box$
\end{flushright}

So we proved that 
\begin{equation*}
r_p^*HF(H,\gamma:\Lambda_0)\otimes \Lambda_0[u^{-1},u]]\langle \theta\rangle \cong \widehat{H}(\mathbb{Z}_p,CF(H,\gamma:\Lambda_0)^{\otimes p})\cong \widehat{HF}_{\mathbb{Z}_p}(H^{(p)},p\gamma)
\end{equation*}
holds.

Next, we prove the following lemma (see also section 3.2 in \cite{CGG3}).
\begin{Lem}
Assume that ${P(H,\gamma)=\{x_1,\cdots,x_k\}}$, ${P(H^{(p)},p\gamma)=\{y_i,\cdots,y_k\}}$ and ${y_i=x_i^{(p)}}$ hold. We also assume that $p$ is an admissible prime number.
Assume that there is an isomorphism
\begin{gather*}
HF(H,\gamma:\Lambda_0)\cong \bigoplus_{j=1}^m\Lambda_0/T^{\beta_j}\Lambda_0  \\
0<\beta_1\le \beta_2 \cdots \le \beta_m
\end{gather*}
and there is an isomorphism 
\begin{gather*}
HF(H^{(p)},p\gamma:\Lambda_0)\cong \bigoplus_{j=1}^{m}\Lambda_0/T^{\delta_j}\Lambda_0  \\
0<\delta_1\le \delta_2 \cdots \le \delta_{m} .
\end{gather*}
Then, ${\delta_1\ge p\beta_1}$ holds.
\end{Lem}
\vspace{5mm}
\textbf{proof}(Lemma 17):
The above isomorphism implies that the following isomorphism holds.
\begin{equation*}
\widehat{HF}_{\mathbb{Z}_p}(H^{(p)},p\gamma)\cong (\bigoplus _{j=1}^m\Lambda_0/T^{p\beta_j}\Lambda_0)\otimes \Lambda_0[u^{-1},u]]\langle \theta \rangle
\end{equation*}
We define spectral numbers ${\sigma(z)}$ and ${\tau(z)}$ for ${z\in CF(H^{(p)},p\gamma:\Lambda_0)}$ as follows:

\begin{gather*}
\sigma(z)=\sup \{\alpha \in \mathbb{R} \ | \ z\in T^{\alpha}CF(H^{(p)},p\gamma:\Lambda_0) \}  \\
\tau(z)=\sigma(d_F(z))-\sigma(z).
\end{gather*}
We also define ${\sigma_{eq}(z)}$ and ${\tau_{eq}(z)}$ for ${z\in \widehat{CF}_{\mathbb{Z}_p}(H^{(p)},p\gamma)}$ as follows:
\begin{gather*}
\sigma_{eq}(z)=\sup \{\alpha \in \mathbb{R} \ | \ z\in T^{\alpha}\widehat{CF}_{\mathbb{Z}_p}(H^{(p)},p\gamma) \}  \\
\tau_{eq}(z)=\sigma_{eq}(d_{eq}(z))-\sigma_{eq}(z).
\end{gather*}

First, we prove the following claim.
\begin{Claim}
Let ${\beta_1}$ and ${\delta_1}$ be positive real numbers defined in the statement of Lemma 17. Then the following equalities hold:
\begin{gather*}
\delta_1=\inf \{\tau(z) \ | \ z\in CF(H^{(p)},p\gamma:\Lambda_0)\} \\
p\beta_1=\inf \{\tau_{eq}(z) \ | \ \widehat{CF}_{\mathbb{Z}_p}(H^{(p)},p\gamma)\}.
\end{gather*}
Moreover, the infimum is achieved for ${z}$ with ${\sigma(z)=0}$ (or ${\sigma_{eq}(z)=0}$).
\end{Claim}

We prove only the first equality (the proof for the second equality is the same). We can choose a cocycle ${z\in  CF(H^{(p)},p\gamma:\Lambda_0)}$ such that ${\sigma(z)=0}$ holds and ${T^{\delta_1}z}$ is a coboundary. So we can choose ${w\in  CF(H^{(p)},p\gamma:\Lambda_0)}$ so that ${d_F(w)=T^{\delta_1}z}$ holds. This implies that ${\tau(w)\le \delta_1}$ and 
\begin{equation*}
\delta_1\ge \inf \{\tau(z) \ | \ z\in CF(H^{(p)},p\gamma:\Lambda_0)\} 
\end{equation*}
holds. Assume that ${\delta_1>\textrm{RHS}}$ hold. Then there is ${z\in CF(H^{(p)},p\gamma:\Lambda_0) }$ such that 
\begin{gather*}
\sigma(z)=0, \ \ \alpha=\tau(z)=\sigma(d_F(z))<\delta_1
\end{gather*}
holds. So, ${w=T^{-\alpha}d_F(z)}$ satisfies ${\sigma(w)=0}$ and ${T^{\alpha}w=d_F(z)}$ is a coboundary. Let ${\{C_1,\cdots,C_{m}\}}$ be a set of cocycles such that the cohomology class ${[C_i]}$ corresponds to the generator of the i-th factor of the right hand side in the following isomorphism:
\begin{equation*}
HF(H^{(p)},p\gamma:\Lambda_0) \cong \bigoplus_{j=1}^{m}\Lambda_0/T^{\delta_j}\Lambda_0 .
\end{equation*}
Now ${[T^{\alpha}w]=0\in HF(H^{(p)},p\gamma:\Lambda_0)  }$ implies that ${[w]\in HF(H^{(p)},p\gamma:\Lambda_0) }$ is written in the following form.

\begin{gather*}
[w]=\sum_{\delta_1-\alpha\le \lambda_{1,j}<\delta_1}a_{1,\lambda_{1,j}}T^{\lambda_{1,j}}[C_1]+\cdots +\sum_{\delta_{m}-\alpha\le \lambda_{m,j}<\delta_{m}}a_{m,\lambda_{m,j}}T^{\lambda_{m,j}}[C_{m}]  \ \ \  (a_{i,j}\in \mathbb{F}_p)
\end{gather*}
Note that ${[w]\neq 0}$ holds because ${\sigma(w)=0}$ and as our cochain complex is strict and any coboundary ${b\in CF(H^{(p)},p\gamma:\Lambda_0) }$ satisfies ${\sigma(b)>0}$. We choose a cochain ${v\in CF(H^{(p)},p\gamma:\Lambda_0)}$ so that 
\begin{gather*}
C=\sum_{\delta_1-\alpha\le \lambda_{1,j}<\delta_1}a_{1,\lambda_{1,j}}T^{\lambda_{1,j}}C_1+\cdots +\sum_{\delta_{m}-\alpha\le \lambda_{m,j}<\delta_{m}}a_{m,\lambda_{m,j}}T^{\lambda_{m,j}}C_{m}  \\
C=w+d_F(v)
\end{gather*}
holds. However ${\sigma(d_F(v))>0}$ implies that 
\begin{equation*}
0<\delta_1-\alpha \le \sigma(C)=\sigma(w+d_F(v))=\sigma(w)=0
\end{equation*}
holds. This is a contradiction and we proved Claim 18.

We fix a basis ${\{X_i^a,Y_i^b,Z_i^c\}}$ of ${CF^{loc}(H^{(p)},x_i)}$ and a basis ${\{\widetilde{X}_i^a,\widetilde{Y}_i^b,\widetilde{Z}_i^c,\widetilde{X}_{i,\theta}^a,\widetilde{Y}_{i,\theta}^b,\widetilde{Z}_{i,\theta}^c\}}$ of ${\widehat{CF}_{\mathbb{Z}_p}^{loc}(H^{(p)},y_i)}$ as in Lemma 13 and Corollary 14. In Corollary 14 we proved that ${\widehat{d}_{eq}}$ has the following form:
\begin{gather*}
\widehat{d}_{eq}\big(\sum_{i,a}\alpha_{i,a}\widetilde{X}_i^a\big)=\iota\big(d_F(\sum_{i,a}\alpha_{i,a}X_i^a)\big)\otimes 1+\sum_{i,a}\beta_{i,a}\widetilde{X}_{i,\theta}^a+\sum_{l\ge 1}\sum_{i,a}\big(\delta_{i,a}\widetilde{X}_i^a+\beta_{i,a}\widetilde{X}_{i,\theta}^a\big)u^l.
\end{gather*}
This implies that 
\begin{gather*}
\sigma\big(d_F(\sum_{i,a}\alpha_{i,a}X_i^a)\big)\ge \sigma_{eq}\big(\widehat{d}_{eq}\big(\sum_{i,a}\alpha_{i,a}\widetilde{X}_i^a\big) \big) \\
\tau(\sum_{i,a}\alpha_{i,a}X_i^a)\ge \tau_{eq}\big(\sum_{i,a}\alpha_{i,a}\widetilde{X}_i^a\big)
\end{gather*}
holds. So 
\begin{gather*}
\delta_1=\inf\{\tau(x) \ | \ x\in CF(H^{(p)},p\gamma:\Lambda_0)\}\ge \inf\{\tau_{eq}(z) \ | \ z\in CF_{\mathbb{Z}_p}(H^{(p)},p\gamma)\}=p\beta_1
\end{gather*} 
holds and we finished the proof of  Lemma 17.
\begin{flushright}
    $\Box$
\end{flushright}

Next we apply Lemma 17 to prove Theorem 2. Recall that we assumed that
\begin{gather*}
P(H,\gamma)=\{x_1,\cdots, x_k\} \\
P(H^{(p)},p\gamma)=\{x_1^{(p)},\cdots, x_k^{(p)}\}
\end{gather*}
holds and ${p}$ is admissible in the beginning of this section. Our purpose is to prove that there is a simple ${p'}$ periodic orbit in ${p\gamma}$. This is equivalent to proving that ${P(H^{(p')},p\gamma)\neq \emptyset}$ because any periodic orbit in ${P(H^{(p')},p\gamma)}$ is simple if $p$ (hence $p'$) is sufficiently large. This follows from the following argument. Let ${\{e_1,\cdots,e_m\}}$ be a basis of ${H_1(M:\mathbb{Z})/\textrm{Tor}}$. Then $\gamma$ is a sum as follows:
\begin{gather*}
\gamma=\alpha_1e_1+\cdots +\alpha_me_m \ \ \ (\alpha_i\in \mathbb{Z}).
\end{gather*}
If there is a non-simple periodic orbit in ${P(H^{(p')},p\gamma)}$, there is ${\gamma'\in H_1(M:\mathbb{Z})/\textrm{Tor}}$ such that
\begin{gather*}
p'\gamma'=p\gamma
\end{gather*}
holds. This implies that $p'$ is a common divisor of ${\{\alpha_1,\cdots,\alpha_m\}}$. In particular, any periodic orbit in ${P(H^{(p')},p\gamma)}$ is simple if ${p>\max \{\alpha_1,\cdots,\alpha_m\}}$ holds. We fix ${C>0}$ independent of $p$ so that 
\begin{gather*}
HF(H,\gamma:\Lambda_0^{\mathbb{F}_p})\cong \bigoplus_{j=1}^m\Lambda_0^{\mathbb{F}_p}/T^{\beta_{p,j}}\Lambda_0^{\mathbb{F}_p} \ \ \ (0<\beta_{p,1}\le \cdots \le \beta_{p,m})  \\
C<\frac{1}{2}\beta_{p,1}
\end{gather*} 
holds for any prime number $p$. Here ${\Lambda^{\mathbb{F}_p}}$ is the universal Novikov ring of the ground field ${\mathbb{F}_p}$. This is always possible because ${\beta_{p,1}}$ is greater than the minimum energy of a solution to the Floer equation, which is independent of $p$ (Lemma 5). Lemma 17 implies that 
\begin{equation*}
\tau(z)\ge p\beta_{p,1}\ge 2pC
\end{equation*}
holds for any ${z\in CF(H^{(p)},p\gamma:\Lambda_0^{\mathbb{F}_p})}$. This implies that any element ${x\neq 0\in HF^{loc}(H^{(p)},x_i^{(p)})}$ determines a cocycle in ${CF^{[-pC,\epsilon)}(H^{(p)},p\gamma:\Lambda^{\mathbb{F}_p})}$ and ${CF^{[-\epsilon,pC)}(H^{(p)},p\gamma:\Lambda^{\mathbb{F}_p})}$ and they are not coboundaries (${\epsilon>0}$ is sufficiently small) because ${\tau(z)\ge 2pC}$ holds for any $z$. This implies that the natural homomorphism
\begin{equation*}
\iota:HF^{[-\epsilon,pC)}(H^{(p)},p\gamma:\Lambda^{\mathbb{F}_p})\longrightarrow HF^{[-pC,\epsilon)}(H^{(p)},p\gamma:\Lambda^{\mathbb{F}_p})
\end{equation*}
is not zero (${\iota([x])\neq 0}$). Let ${p'}$ be the first prime number greater than ${p}$. We assume $p$ is a sufficiently large prime number such that 
\begin{equation*}
2(p'-p)||H||<pC
\end{equation*}
holds. This is possible because ${p'-p=o(p)}$ holds (see \cite{BHP}). We have two continuation homomorphism as follows:
\begin{gather*}
F:HF^{[-\epsilon,pC)}(H^{(p)},p\gamma:\Lambda^{\mathbb{F}_p})\longrightarrow HF^{[-\epsilon-(p'-p)||H||,pC-(p'-p)||H||)}(H^{(p')},p\gamma:\Lambda^{\mathbb{F}_p})  \\
G:HF^{[-\epsilon-(p'-p)||H||,pC-(p'-p)||H||)}(H^{(p')},p\gamma:\Lambda^{\mathbb{F}_p})\longrightarrow HF^{[-pC,\epsilon)}(H^{(p)},p\gamma:\Lambda^{\mathbb{F}_p}).
\end{gather*}
The composition of $F$ and $G$ satisfies ${GF=\iota\neq 0}$. So, ${P(H^{(p')},p\gamma)\neq \emptyset}$ holds and we proved the theorem.

\section{Weakly monotone case}

\subsection{Preliminary}
 In this section, we prove Theorem 2 for weakly monotone symplectic manifolds. The difference between toroidally monotone cases and weakly monotone cases is that we have not constructed  ${\mathbb{Z}_p}$-equivariant Floer cohomology and ${\mathbb{Z}_p}$-equivariant pair of pants product in the weakly monotone case. So to prove Theorem 2 for the weakly monotone case, it suffices to construct these theories. The rest of the proof is the same as in the toroidally monotone case.

 In the toroidally monotone case, we can exclude sphere bubbles easily because the Maslov index of a holomorphic sphere is greater than or equal to $2$. In the weakly monotone case, we have to exclude sphere bubbles much more carefully. Recall the construction of Floer cohomology theory for weakly monotone symplectic manifolds \cite{HS,On}. For a generic choice of an almost complex structure $J$, there are no holomorphic spheres with negative Chern numbers. Let ${H}$ be a Hamiltonian function. The pair ${(H,J)}$ is not necessary a Floer regular pair. However, we can perturb $H$ to $\widetilde{H}$ so that ${(\widetilde{H},J)}$ is Floer regular and sphere bubbles do not appear in the definition of the Floer coboundary operator ${d_F}$.
 
Recall that in the definition of ${\mathbb{Z}_p}$-equivariant Floer cohomology for toroidally monotone symplectic manifolds, we considered a family of almost complex structures ${\{J_{w,t}\}}$ parametrized by ${S^{\infty}}$ and ${S^1}$ while fixing a Hamiltonian function ${H^{(p)}}$. So one possible modification for the weakly monotone case is to consider a family of Hamiltonian functions parametrized by ${S^{\infty}}$ and ${S^1}$ while fixing an almost complex structure ${J}$. Let ${K\in C^{\infty}(S^1\times M)}$ be a perturbation of ${H^{(p)}}$ such that ${(K,J)}$ is a Floer regular pair. Note that ${K}$ is not necessarily a ${\frac{1}{p}}$-periodic Hamiltonian function. We also consider a family of Hamiltonian functions ${\mathcal{K}_{w,t}}$ parametrized by ${(w,t)\in S^{\infty}\times S^1}$ which satisfies the following conditions (compare it to the definition of ${J_{w,t}}$ in section 3).

\begin{itemize}
\item (locally constant at critical points) For all ${w}$ in a small neighborhood of ${Z_i^m\in S^{\infty}}$, 
\begin{equation*}
\mathcal{K}_{w,t}=K_{t-\frac{m}{p}}.
\end{equation*}
\item ($\mathbb{Z}_p$-equivariance) ${\mathcal{K}_{m\cdot w,t}=\mathcal{K}_{w,t-\frac{m}{p}}}$ holds for any ${m\in \mathbb{Z}_p}$ and ${w\in S^{\infty}}$.
\item (invariance under the shift $\tau$) ${\mathcal{K}_{\tau(w),t}=\mathcal{K}_{w,t}}$ holds.
\end{itemize}

We consider the following equation for ${x,y\in P(K)}$, ${m\in \mathbb{Z}_p}$ and ${i\in \mathbb{Z}}$:

\begin{gather*}
(u,v)\in C^{\infty}(\mathbb{R}\times S^1,M)\times C^{\infty}(\mathbb{R},S^{\infty})  \\
\partial_su(s,t)+J(u(s,t))(\partial_tu(s,t)-X_{\mathcal{K}_{v(s),t}}(u(s,t)))=0  \\
\frac{d}{ds}v(s)-\textrm{grad}\widetilde{F}=0  \\
\lim_{s\to -\infty}v(s)=Z_{\alpha}^0, \lim_{s\to +\infty}v(s)=Z_i^m, \lim_{s\to -\infty}u(s,t)=x(t), \lim_{s\to +\infty}u(s,t)=y(t-\frac{m}{p}) .
\end{gather*}
One might try to define 
\begin{equation*}
d_{\alpha}^{i,m}:CF(K,\gamma:\Lambda_0)\longrightarrow CF(K,\gamma:\Lambda_0)
\end{equation*}
by counting above solutions and define $d_{eq}$ in terms of the ${d_{\alpha}^{i,m}}$. However, this attempt contains the following difficulty. The action gap of the solution ${(u,v)}$ of the above equation
\begin{equation*}
\int_{\mathbb{R}\times S^1}u^*\omega+\int_0^1K(t,x(t))-K(t,y(t))dt
\end{equation*}
is not necessarily non-negative. This problem happens when ${i\in \mathbb{Z}}$ becomes sufficiently large. As $i$ becomes bigger and bigger, the effect of the perturbation ${H_t^{(p)}\to \mathcal{K}_{w,t}}$ becomes bigger and we cannot define a coboundary operator ${d_{eq}}$ over ${\Lambda_0}$. What we can do is to fix some ${N\in \mathbb{Z}}$ and define a finite operators ${\{d_{\alpha}^{i,m}\}}$ for ${i\le N}$.

This difficulty is very similar to the difficulty in \cite{FOOO}. In \cite{FOOO}, Fukaya, Oh, Ohta and Ono constructed an ${A_{\infty}}$-algebra associated to a Lagrangian submanifold in a symplectic manifold. They constructed ${A_{\infty}}$-operators ${\{m_k\}_{k\in \mathbb{Z}_{\ge0}}}$ by using moduli spaces of holomorphic discs bounding the Lagrangian submanifold. To determine ${\{m_k\}_{k\in \mathbb{Z}_{\ge0}}}$, they had to achieve transversality of infinitely many interrelated moduli spaces by perturbing multisections of Kuranishi structures. Unfortunately, this is impossible because perturbations of lower-order operators influence perturbations of higher-order operators, and the higher perturbations become bigger and bigger. So what they could do is to construct finite operators ${\{m_{0,0},\cdots,m_{n,K}\}}$ (${A_{n,K}}$-algebra) by one perturbation. They constructed ${A_{\infty}}$-operators ${\{m_k\}_{k\in \mathbb{Z}_{\ge0}}}$ by ``gluing" infinitely many ${A_{n,K}}$-algebras as ${(n,K)\to (+\infty,+\infty)}$ by applying homological algebra developed in \cite{FOOO}. Our situation is much simpler because we do not have to consider higher algebraic operators (we only need a coboundary operator) and we do not have to consider the space of infinitely many singular chains. So we can mimic the construction of Lagrangian ${A_{\infty}}$-algebra by applying the algebraic machinery developed in \cite{FOOO}. In the rest of this section, we explain how we can apply \cite{FOOO} in our cases.

\subsection{${X_K}$-module, ${X_K}$-morphism and ${X_K}$-homotopy}
We consider the following situation. Let $C$ be a $\mathbb{Z}_2$-graded ${\Lambda_0}$-module. We want to construct a family of degree $1$ maps ${\{d_i:C\rightarrow C\}_{i\in \mathbb{Z}_{\ge 0}}}$ so that ${\sum_{i+j=k}d_id_j=0}$ holds for any ${k\in \mathbb{Z}_{\ge 0}}$. Then the infinite sum
\begin{equation*}
d=d_0+d_1u+d_2u^2+d_3u^3+\cdots
\end{equation*}
becomes a coboudary operator on ${C\otimes \Lambda_0[[u]]}$ as follows:
\begin{gather*}
d:C\otimes \Lambda_0[[u]]\longrightarrow C\otimes \Lambda_0[[u]] \\
x\otimes u^m \mapsto \sum_{i\ge 0}d_i(x)\otimes u^{m+i}.
\end{gather*}

\begin{Def}
\begin{enumerate}
\item Let ${\mathfrak{f}=\{f_i:C\rightarrow D\}_{i=0}^K}$ be a family of maps between ${\Lambda_0}$-modules. We define ${\mathfrak{f}^{(L)}\in Hom(C,D)\otimes \Lambda_0[[u]]}$ for some ${L\in \mathbb{Z}_{\ge 0}}$ as follows:
\begin{gather*}
\mathfrak{f}^{(L)}=\begin{cases}f_0+f_1\otimes u+\cdots +f_L\otimes u^L & L\le K  \\ 
f_0+f_1\otimes u+\cdots +f_K\otimes u^K & L\ge K.
\end{cases}
\end{gather*}
\item Let ${\mathfrak{d}=\{d_i:C\rightarrow C\}}_{i=0}^K$ be a family of degree $1$ maps (${0\le K\le \infty}$). We call ${(C,\mathfrak{d})}$ an ${X_K}$-module if ${\mathfrak{d}}$ satisfies 
\begin{equation*}
\mathfrak{d}^{(K)}\circ \mathfrak{d}^{(K)}\equiv 0 \ \ \ \textrm{mod}(u^{K+1}). 
\end{equation*}
\item Let ${(C,\mathfrak{d})}$ and ${(D,\mathfrak{l})}$ be ${X_K}$-modules. An ${X_K}$-morphism between ${C}$ and ${D}$ is a family of degree $0$ maps ${\mathfrak{f}=\{f_i\}_{i=0}^K}$ which satisfies the following equality:
\begin{equation*}
\mathfrak{f}^{(K)}\circ \mathfrak{d}^{(K)}\equiv \mathfrak{l}^{(K)}\circ \mathfrak{f}^{(K)} \ \ \textrm{mod}(u^{K+1}).
\end{equation*}
\end{enumerate}
\end{Def}

We define a ${\Lambda_0}$-algebra ${\Lambda_0^{(K)}}$ as follows:
\begin{gather*}
\Lambda_0^{(K)}=\begin{cases} \Lambda_0[[u]]/(u^{K+1}) & K<+\infty \\ \Lambda_0[[u]] & K=+\infty .\end{cases}
\end{gather*}
For any ${\Lambda_0}$-module $C$, we define ${\Lambda_0^{(K)}}$-module ${C^{(K)}}$ by
\begin{gather*}
C^{(K)}=C\otimes_{\Lambda_0}\Lambda_0^{(K)} .
\end{gather*}
We define ${\textrm{deg}(u)=2}$, so ${C^{(K)}}$ is also ${\mathbb{Z}_2}$-graded. Note that ${(C,\mathfrak{d})}$ is an ${X_K}$-module if and only if ${\mathfrak{d}^{(K)}\circ\mathfrak{d}^{(K)}=0 }$ holds on ${C^{(K)}}$ and ${\mathfrak{f}:(C,\mathfrak{d})\rightarrow (D,\mathfrak{l})}$ is an ${X_K}$-morphism if and only if 
\begin{gather*}
\mathfrak{f}^{(K)}:(C^{(K)},\mathfrak{d}^{(K)})\longrightarrow (D^{(K)},\mathfrak{l}^{(K)})
\end{gather*}
is a cochain map.

Let ${\mathfrak{f}_0}$ and ${\mathfrak{f_1}}$ be ${X_K}$-morphisms from ${(C,\mathfrak{d})}$ to ${(D,\mathfrak{l})}$. We say ${\mathfrak{f}_0}$ and ${\mathfrak{f}_1}$ are ${X_K}$-homotopic if there is a family of maps ${\mathfrak{h}=\{h_i:C\rightarrow D\}_{i=0}^K}$ which satisfies the following relations:

\begin{equation*}
\mathfrak{f}_0^{(K)}-\mathfrak{f}_1^{(K)}\equiv \mathfrak{h}^{(K)}\circ \mathfrak{d}^{(K)}+\mathfrak{l}^{(K)}\circ \mathfrak{h}^{(K)} \ \ \ \textrm{mod}(u^{K+1}).
\end{equation*}
The composition of two ${X_K}$-morphisms ${\{f_i:C\rightarrow D\}}$ and ${\{g_j:D\rightarrow E\}}$ is defined by ${(\mathfrak{f}\circ \mathfrak{g})_i=\sum_{j+l=i}f_jg_l}$. Note that 

\begin{equation*}
(\mathfrak{f}\circ \mathfrak{g})^{(K)}\equiv \mathfrak{f}^{(K)}\circ \mathfrak{g}^{(K)} \ \ \ \textrm{mod}(u^{K+1})
\end{equation*}
holds. We call an ${X_K}$-morphism ${\mathfrak{f}:(C,\mathfrak{d})\rightarrow (D,\mathfrak{l})}$ an ${X_K}$-homotopy equivalence if there is an ${X_K}$-morphism ${\mathfrak{g}:(D,\mathfrak{l})\rightarrow (C,\mathfrak{d})}$ such that ${\mathfrak{f}\mathfrak{g}}$ and ${\mathfrak{g}\mathfrak{f}}$ are ${X_K}$-homotopic to the identity. 

We also give an equivalent definition of ${X_K}$-homotopy. For this purpose, we introduce a new ${X_K}$-module ${(C\times [0,1],\widetilde{\mathfrak{d}}=\{\widetilde{d}_i\})}$ for an ${X_K}$-module ${(C,\mathfrak{d})}$. We define ${\widetilde{d}^{(K)}=\sum_{i=0}^{K}\widetilde{d}_iu^i}$ as follows:
\begin{gather*}
C\times [0,1]=C\oplus C[-1]\oplus C \\
\begin{cases}
\widetilde{d}^{(K)}(x,0,0)=(d^{(K)}(x),(-1)^{\textrm{deg}x}x,0)  \ \ \ \ \  (x\in C^{(K)})\\
\widetilde{d}^{(K)}(0,y,0)=(0,d^{(K)}(y),0) \ \ \ \ \  (y\in(C[-1])^{(K)}) \\
\widetilde{d}^{(K)}(0,0,z)=(0,-(-1)^{\textrm{deg}z}z,d^{(K)}(z)) \ \ \ \ \ (z\in C^{(K)}).
\end{cases}
\end{gather*}
Note that ${C[-1]}$ is a copy of $C$ with degree shift $1$:
\begin{equation*}
C[-1]^m=C^{m-1}.
\end{equation*}
It is straightforward to see that ${\widetilde{d}^{(K)}\widetilde{d}^{(K)}=0}$ holds on ${(C\times [0,1])^{(K)}}$. So ${(C\times [0,1],\widetilde{\mathfrak{d}})}$ is an ${X_K}$-module. We also define ${X_K}$-morphisms ${\textrm{Incl}=\{(\textrm{Incl})_i\}_{i\le K}}$, ${\textrm{Eval}_{s=0}=\{(\textrm{Eval}_{s=0})_i\}_{i\le K}}$ and ${\textrm{Eval}_{s=1}=\{(\textrm{Eval}_{s=1})_i\}_{i\le K}}$ between ${(C,\mathfrak{d})}$ and ${(C\times [0,1],\widetilde{\mathfrak{d}})}$ as follows:

\begin{gather*}
\textrm{Incl}:(C,\mathfrak{d})\longrightarrow (C\times [0,1],\widetilde{\mathfrak{d}})  \\
\textrm{(Incl)}_i(x)=\begin{cases}(x,0,x)  & i=0 \\ (0,0,0) & i\ge 1  \end{cases}  
\end{gather*}
\begin{gather*}
\textrm{Eval}_{s=0}:(C\times [0,1],\widetilde{\mathfrak{d}})\longrightarrow (C,\mathfrak{d}) \\
(\textrm{Eval}_{s=0})_i(x,y,z)=\begin{cases}x  & i=0  \\  0  & i\ge 1  \end{cases}
\end{gather*}
\begin{gather*}
\textrm{Eval}_{s=1}:(C\times [0,1],\widetilde{\mathfrak{d}})\longrightarrow (C,\mathfrak{d}) \\
(\textrm{Eval}_{s=1})_i(x,y,z)=\begin{cases}z  & i=0  \\ 0  & i\ge 1.  \end{cases}
\end{gather*}

\begin{Def}
Let ${\mathfrak{f}_0}$ and ${\mathfrak{f}_1}$ be ${X_K}$-morphisms between ${(C,\mathfrak{d})}$ and ${(D,\mathfrak{l})}$. We say ${\mathfrak{f}_0}$ and ${\mathfrak{f}_1}$ are ${X_K}$-homotopic if there is an ${X_K}$-morphism ${\mathfrak{h}:C\rightarrow D\times [0,1]}$ such that ${\textrm{Eval}_{s=0}\circ \mathfrak{h}=\mathfrak{f}_0}$ and ${\textrm{Eval}_{s=1}\circ \mathfrak{h}=\mathfrak{f}_1}$ hold. In this case $\mathfrak{h}$ is called an ${X_K}$-homotopy between ${\mathfrak{f}_0}$ and ${\mathfrak{f}_1}$.
\end{Def}
The next lemma implies that the above two definitions of $X_K$-homotopy are equivalent. 
\begin{Lem}
Let ${\mathfrak{f}_0}$ and ${\mathfrak{f}_1}$ be ${X_K}$-morphisms between ${(C,\mathfrak{d})}$ and ${(D,\mathfrak{l})}$. A family of degree $-1$ maps ${\mathfrak{h}=\{h_i\}_{i=0}^K\subset \textrm{Hom}(C,D)}$ satisfies
\begin{gather*}
\mathfrak{f}_0^{(K)}-\mathfrak{f}_1^{(K)}\equiv \mathfrak{h}^{(K)}\circ \mathfrak{d}^{(K)}+\mathfrak{l}^{(K)}\circ \mathfrak{h}^{(K)} \ \ \ \textrm{mod}(u^{K+1})
\end{gather*}
if and only if the family of maps
\begin{gather*}
\widetilde{\mathfrak{h}}=\{\widetilde{h}_i\};(C,\mathfrak{d})\longrightarrow (D\times [0,1], \widetilde{l})  \\
\widetilde{h}_i(x)=((f_0)_i(x), -(-1)^{\textrm{deg}x}h_i(x),(f_1)_i(x))
\end{gather*}
is an ${X_K}$-homotopy between ${\mathfrak{f}_0}$ and ${\mathfrak{f}_1}$.
\end{Lem}
\vspace{5mm}
\textbf{proof}:
We calculate ${\widetilde{\mathfrak{l}}^{(K)}\circ \widetilde{\mathfrak{h}}^{(K)}}$ and ${\widetilde{\mathfrak{h}}^{(K)}\circ \mathfrak{d}^{(K)}}$. Note that ${\widetilde{\mathfrak{h}}}$ is an ${X_K}$-homotopy between ${\mathfrak{f}_0}$ and ${\mathfrak{f}_1}$ if and only if they coincide modulo ${(u^{K+1})}$.
\begin{gather*}
\widetilde{\mathfrak{l}}^{(K)}\circ \widetilde{\mathfrak{h}}^{(K)}(x)\\ =\big(\mathfrak{l}^{(K)}\mathfrak{f}_0^{(K)}(x),-(-1)^{\textrm{deg}x}\mathfrak{l}^{(K)}\mathfrak{h}^{(K)}(x)+(-1)^{\textrm{deg}x}\mathfrak{f}_0^{(K)}(x)-(-1)^{\textrm{deg}x}\mathfrak{f}_1^{(K)}(x),\mathfrak{l}^{(K)}\mathfrak{f}_1^{(K)}(x)\big)
\end{gather*}
\begin{gather*}
\widetilde{\mathfrak{h}}^{(K)}\circ \mathfrak{d}^{(K)}(x)=(\mathfrak{f}_0^{(K)}\mathfrak{d}^{(K)}(x), (-1)^{\textrm{deg}x}\mathfrak{h}^{(K)}\mathfrak{d}^{(K)}(x),\mathfrak{f}_1^{(K)}\mathfrak{d}^{(K)}(x))
\end{gather*}
By comparing the middle factors of the right hand sides, we can see that ${\widetilde{\mathfrak{h}}}$ is an ${X_K}$-homotopy between ${\mathfrak{f}_0}$ and ${\mathfrak{f}_1}$ if and only if
\begin{gather*}
\mathfrak{f}_0^{(K)}-\mathfrak{f}_1^{(K)}\equiv \mathfrak{h}^{(K)}\circ \mathfrak{d}^{(K)}+\mathfrak{l}^{(K)}\circ \mathfrak{h}^{(K)} \ \ \ \textrm{mod}(u^{K+1})
\end{gather*}
holds.
\begin{flushright}
    $\Box$
\end{flushright}

\begin{Lem}
The following two maps
\begin{gather*}
(\textrm{Incl})_0:(C,d_0)\longrightarrow (C\times [0,1],\widetilde{d}_0) \\
(\textrm{Eval}_{s=s_0})_0:(C\times [0,1],\widetilde{d}_0)\longrightarrow (C,d_0)  \ \ \  (s_0=0,1)
\end{gather*}
are cochain homotopy equivalences.
\end{Lem}
\vspace{5mm}
\textbf{proof}:
It suffices to prove that ${(\textrm{Eval}_{s=s_0})_0(\textrm{Incl})_0}$ and ${(\textrm{Incl})_0(\textrm{Eval}_{s=s_0})_0}$ are homotopic to the identity. The former case is trivial because it is equal to the identity.
\begin{gather*}
{(\textrm{Eval}_{s=s_0})_0(\textrm{Incl})_0}=\textrm{Id}
\end{gather*}
\begin{gather*}
\big\{\textrm{Id}-(\textrm{Incl})_0\circ (\textrm{Eval}_{s=0})_0\big\}(x,y,z)=(0,y,z-x)
\end{gather*}
\begin{gather*}
\big\{\textrm{Id}-(\textrm{Incl})_0\circ (\textrm{Eval}_{s=1})_0\big\}(x,y,z)=(x-z,y,0)
\end{gather*}
For the sake of simplicity, we do not use the shifted degree of ${C\times [0,1]}$ here. We can define cochain homotopies ${h_{s=0}}$ and ${h_{s=1}}$ as follows:
\begin{gather*}
h_{s=0}:C\oplus C\oplus C\longrightarrow C\oplus C\oplus C  \\
(x,y,z)\mapsto (0,0,-(-1)^{\textrm{deg}y}y)
\end{gather*}
\begin{gather*}
h_{s=1}:C\oplus C\oplus C\longrightarrow C\oplus C\oplus C  \\
(x,y,z)\mapsto ((-1)^{\textrm{deg}y}y,0,0).
\end{gather*}
The following calculations imply that ${h_{s=s_0}}$ is a homotopy between ${\textrm{Id}}$ and ${(\textrm{Incl})_0\circ (\textrm{Eval}_{s=s_0})_0}$:
\begin{gather*}
\big\{h_{s=0}\circ \widetilde{d}_0+\widetilde{d}_0\circ h_{s=0}\big\}(x,y,z)\\
=h_{s=0}\big(d_0(x),(-1)^{\textrm{deg}x}x+d_0(y)-(-1)^{\textrm{deg}z}z,d_0(z)\big)
+\widetilde{d}_0\big(0,0,-(-1)^{\textrm{deg}y}y\big)  \\
=\big(0,0,-x+(-1)^{\textrm{deg}y}d_0(y)+z\big)+\big(0,y,-(-1)^{\textrm{deg}y}d_0(y)\big) \\
=(0,y,z-x)=\big\{\textrm{Id}-(\textrm{Incl})_0\circ(\textrm{Eval}_{s=0})_0\big\}(x,y,z)
\end{gather*}
\begin{gather*}
\big\{h_{s=1}\circ \widetilde{d}_0+\widetilde{d}_0\circ h_{s=1}\big\}(x,y,z) \\
=h_{s=1}\big(d_0(x),(-1)^{\textrm{deg}x}x+d_0(y)-(-1)^{\textrm{deg}z}z,d_0(z),d_0(z)\big)
+\widetilde{d}_0\big((-1)^{\textrm{deg}y}y,0,0\big) \\
=\big(x-(-1)^{\textrm{deg}y}d_0(y)-z,0,0\big)+\big((-1)^{\textrm{deg}y}d_0(y),y,0\big) \\
=(x-z,y,0)=\big\{\textrm{Id}-(\textrm{Incl})_0\circ(\textrm{Eval}_{s=1})_0\big\}(x,y,z).
\end{gather*}
\begin{flushright}
    $\Box$
\end{flushright}

We need a cohomology theory on the space of homomorphisms between two cochain complexes  ${\textrm{Hom}((A,d_A),(B,d_B))}$.

\begin{Def}
Let ${(A,d_A)}$ and ${(B,d_B)}$ be ${\mathbb{Z}_2}$-graded differential modules. We define a coboundary operator on ${\textrm{Hom}(A,B)}$ as follows:
\begin{gather*}
\partial_{\textrm{Hom}((A,d_A),(B,d_B))}:\textrm{Hom}(A,B)\longrightarrow \textrm{Hom}(A,B) \\
\phi \mapsto d_B\circ \phi-(-1)^{\textrm{deg}\phi}\phi\circ d_A.
\end{gather*}
We denote the cohomology of this cochain complex by ${H\big(\textrm{Hom}(A,B),\partial_{\textrm{Hom}((A,d_A),(B,d_B))}\big)}$.
\end{Def}

\begin{Lem}
Let ${(A,d_A)}$, ${(A',d_{A'})}$, ${(B,d_B)}$ and ${(B',d_{B'})}$ be ${\mathbb{Z}_2}$-graded complexes. Assume that there are cochain maps
\begin{gather*}
f:(A',d_{A'})\longrightarrow (A,d_A)\\
g:(B,d_{B})\longrightarrow (B',d_{B'}) .
\end{gather*}
Then, 
\begin{gather*}
(g,f)_*:\textrm{Hom}(A,B)\longrightarrow \textrm{Hom}(A',B')  \\
\phi \mapsto g\circ \phi \circ f 
\end{gather*}
is a cochain map. Moreover, if ${f'}$ is cochain homotopic to ${f}$ and ${g'}$ is a cochain homotopic to $g$, ${(g',f')_*}$ is also cochain homotopic to ${(g,f)_*}$.
\end{Lem}
\vspace{5mm}
\textbf{proof}:
An equality
\begin{gather*}
(g,f)_*(\partial_{\textrm{Hom}(A,B)}(\phi))=g\circ (d_B\phi-(-1)^{\textrm{deg}\phi}\phi d_A)\circ f \\
=d_{B'}g\phi f-(-1)^{\textrm{deg}\phi}g\phi f d_{A'} =\partial_{\textrm{Hom}(A',B')}(g\phi f)  \\
=\partial_{\textrm{Hom}(A',B')}\big((g,f)_*(\phi)\big)
\end{gather*}
holds. So, ${(g,f)_*}$ is a cochain map. Assume that ${f}$ is cochain homotopic to $f$ and ${g'}$ is cochain homotopic to ${g}$. It suffices to prove that ${(g,\textrm{Id})_*}$ and ${(g',\textrm{Id})_*}$ are cochain homotopic and ${(\textrm{Id},f)_*}$ and ${(\textrm{Id},f)_*}$ are cochain homotopic. Let ${h_{(f,f')}}$ be a homotopy between ${f}$ and ${f'}$ and let ${h_{(g,g')}}$ be a homotopy between ${g}$ and ${g'}$. Then
\begin{gather*}
f-f'=d_A\circ h_{(f,f')}+h_{(f,f')}\circ d_{A'}  \\
g-g'=d_{B'}\circ h_{(g,g')}+h_{(g,g')}\circ d_{B}
\end{gather*}
holds. Next, we define cochain homotopies ${H_{(f,f')}}$ and ${H_{(g,g')}}$ as follows:
\begin{gather*}
H_{(f,f')}:\textrm{Hom}(A,B)\longrightarrow \textrm{Hom}(A',B) \\
\phi \mapsto (-1)^{\textrm{deg}\phi}\phi \circ h_{(f,f')}
\end{gather*}
\begin{gather*}
H_{(g,g')}:\textrm{Hom}(A,B)\longrightarrow \textrm{Hom}(A,B') \\
\phi \mapsto  h_{(g,g')}\circ \phi.
\end{gather*}
The following calculations imply that ${H_{(f,f')}}$ is a cochain homotopy between ${(\textrm{Id},f)_*}$ and ${(\textrm{Id},f')_*}$, and ${H_{(g,g')}}$ is a cochain homotopy between ${(g,\textrm{Id})_*}$ and ${(g',\textrm{Id})_*}$:
\begin{gather*}
\big\{\partial_{\textrm{Hom(A',B)}}\circ H_{(f,f')}+H_{(f,f')}\circ \partial_{\textrm{Hom}(A,B)}\big\}(\phi) \\
=\partial_{\textrm{Hom}(A',B)}\big((-1)^{\textrm{deg}\phi}\phi\circ h_{(f,f')}\big)+H_{(f,f')}\big(d_B\phi-(-1)^{\textrm{deg}\phi}\phi d_A\big) \\
=(-1)^{\textrm{deg}\phi}d_B\phi h_{(f,f')}+\phi h_{(f,f')}d_{A'}+(-1)^{\textrm{deg}\phi +1}(d_B\phi h_{(f,f')}-(-1)^{\textrm{deg}\phi}\phi d_Ah_{(f,f')}) \\
=\phi(h_{(f,f')}d_{A'}+d_Ah_{(f,f')})=\phi \circ f-\phi \circ f' \\
=\big\{(\textrm{Id},f)_*-(\textrm{Id},f')_*\big\}(\phi)
\end{gather*}
\begin{gather*}
\big\{\partial_{\textrm{Hom}(A,B')}\circ H_{(g,g')}+H_{(g,g')}\circ \partial_{\textrm{Hom}(A,B)}\big\}(\phi) \\
=\partial_{\textrm{Hom}(A,B')}(h_{(g,g')}\phi)+H_{(g,g')}\big(d_B\phi-(-1)^{\textrm{deg}\phi}\phi d_A\big) \\
=d_{B'}h_{(g,g')}\phi-(-1)^{\textrm{deg}\phi+1}(h_{(g,g')}\phi d_A)+h_{(g,g')}d_B\phi-(-1)^{\textrm{deg}\phi}h_{(g,g')}\phi d_A \\
=(d_{B'}h_{(g,g')}+h_{(g,g')}d_B)\phi =g\circ \phi-g'\circ \phi \\
=\big\{(g,\textrm{Id})_*-(g',\textrm{Id})\big\}(\phi).
\end{gather*}
\begin{flushright}
    $\Box$
\end{flushright}

The next proposition describes an obstruction for extensions of ${X_K}$-morphisms.

\begin{Prop}
Let ${(C,\mathfrak{d})}$ and ${(D,\mathfrak{l})}$ be ${X_{K+1}}$-modules. Assume that ${\mathfrak{f}:C\rightarrow D}$ is an ${X_K}$-morphism between the restrictions of $C$ and $D$ to $X_K$-modules. Then there is degree $1$ element 
\begin{equation*}
\mathfrak{o}_{K+1}(\mathfrak{f})\in \textrm{Hom}(C,D)
\end{equation*}
which satisfies the following properties:
\begin{enumerate}
\item $\partial_{\textrm{Hom}((C,d_0),(D,l_0))}(\mathfrak{o}_{K+1}(\mathfrak{f}))=0$ holds. So ${\mathfrak{o}_{K+1}(\mathfrak{f})}$ is a cocycle.
\item We can extend ${\mathfrak{f}}$ to an ${X_{K+1}}$-morphism if and only if ${\mathfrak{o}_{K+1}(\mathfrak{f})}$ is a coboundary  (${[\mathfrak{o}_{K+1}(\mathfrak{f})]=0}$).
\item If ${\mathfrak{f}}$ and ${\mathfrak{f}'}$ are ${}X_K$-homotopic,
\begin{gather*}
[\mathfrak{o}_{K+1}(\mathfrak{f})]=[\mathfrak{o}_{K+1}(\mathfrak{f}')]\in H\big(\textrm{Hom}(C,D),\partial_{\textrm{Hom}((C,d_0),(D,l_0))}\big)
\end{gather*}
holds.
\item Let ${\mathfrak{g}:C'\rightarrow C}$ and ${\mathfrak{g}':D\rightarrow D'}$ be ${X_{K+1}}$-morphisms. Then,
\begin{gather*}
[\mathfrak{o}_{K+1}(\mathfrak{g}'\circ \mathfrak{f}\circ \mathfrak{g})]=(g_{0}',g_0)_*[\mathfrak{o}_{K+1}(\mathfrak{f})]
\end{gather*}
holds.
\end{enumerate}
\end{Prop}
\vspace{5mm}
\textbf{proof}:
Let ${\mathfrak{f}:(C,\mathfrak{d})\rightarrow (D,\mathfrak{l})}$ be an ${X_K}$-morphism. Then, 
\begin{gather*}
\mathfrak{l}^{(K+1)}\circ \mathfrak{f}^{(K)}-\mathfrak{f}^{(K)}\circ \mathfrak{d}^{(K+1)}\equiv \mathfrak{o}_{K+1}(\mathfrak{f})\otimes u^{K+1} \ \ \ \textrm{mod}(u^{K+2})
\end{gather*}
holds for some ${\mathfrak{o}_{K+1}(\mathfrak{f})\in \textrm{Hom}(C,D)}$. This is the definition of ${\mathfrak{o}_{K+1}(\mathfrak{f})}$. Note that 
\begin{gather*}
\mathfrak{l}^{(K+1)}\circ \mathfrak{f}^{(K)}-\mathfrak{f}^{(K)}\circ \mathfrak{d}^{(K+1)}=\partial_{\textrm{Hom}((C^{(K+1)},\mathfrak{d}^{(K+1)})(D^{(K+1)},\mathfrak{l}^{(K+1)}))}(\mathfrak{f}^{(K)})
\end{gather*}
holds. This implies that 
\begin{gather*}
0=\big\{\partial_{\textrm{Hom}((C^{(K+1)},\mathfrak{d}^{(K+1)}),(D^{(K+1)},\mathfrak{l}^{(K+1)}))}\big\}^2(\mathfrak{f}^{(K)}) \\
=\partial_{\textrm{Hom}((C^{(K+1)},\mathfrak{d}^{(K+1)}),(D^{(K+1)},\mathfrak{l}^{(K+1)}))}(\mathfrak{o}_{K+1}(\mathfrak{f})\otimes u^{K+1})  \\
=\partial_{\textrm{Hom}((C,d_0),(D,l_0))}(\mathfrak{o}_{K+1}(\mathfrak{f}))\otimes u^{K+1}
\end{gather*}
holds on ${\textrm{Hom}((C^{(K+1)},\mathfrak{d}^{(K+1)})(D^{(K+1)},\mathfrak{l}^{(K+1)}))}$. So 
\begin{gather*}
\partial_{\textrm{Hom}((C,d_0),(D,l_0))}(\mathfrak{o}_{K+1}(\mathfrak{f}))=0
\end{gather*}
holds and ${\mathfrak{o}_{K+1}(\mathfrak{f})\in \textrm{Hom}((C,d_0),(D,l_0))}$ is a cocycle.
Now ${\mathfrak{f}=\{f_0,\cdots,f_K\}}$ and ${f_{K+1}}$ determines an ${X_{K+1}}$-morphism if and only if 
\begin{gather*}
0\equiv \mathfrak{l}^{(K+1)}\circ (\mathfrak{f}^{(K)}+f_{K+1}\otimes u^{K+1})- (\mathfrak{f}^{(K)}+f_{K+1}\otimes u^{K+1})\circ \mathfrak{d}^{(K+1)} \\
\equiv \big(\mathfrak{o}_{K+1}(\mathfrak{f})+\partial_{\textrm{Hom}((C,d_0),(D,l_0))}(f_{K+1})\big)\otimes u^{K+1} \ \ \ \textrm{mod}(u^{K+2})
\end{gather*}
holds. So, ${\mathfrak{f}}$ can be extended to an ${X_{K+1}}$-morphism if and only if ${\mathfrak{o}_{K+1}(\mathfrak{f})}$ is a coboundary. 
Next, we prove $(\mathrm{iv})$. Let ${\mathfrak{g}:C'\rightarrow C}$ and ${\mathfrak{g}':D\rightarrow D'}$ be ${X_{K+1}}$-morphisms. Then there is ${e:C'\rightarrow D'}$
such that
\begin{gather*}
(\mathfrak{g}'\circ \mathfrak{f}\circ \mathfrak{g})^{(K)}\equiv \mathfrak{g}'^{(K)}\circ \mathfrak{f}^{(K)}\circ \mathfrak{g}^{(K)}+e\otimes u^{K+1} \ \ \textrm{mod}(u^{K+2})
\end{gather*}
holds. Then, 
\begin{gather*}
\mathfrak{o}_{K+1}(\mathfrak{g}'\circ \mathfrak{f}\circ \mathfrak{g})\otimes u^{K+1}\equiv \mathfrak{l}'^{(K+1)}\circ (\mathfrak{g}'\circ \mathfrak{f}\circ \mathfrak{g})^{(K)}-
(\mathfrak{g}'\circ \mathfrak{f}\circ \mathfrak{g})^{(K)}\circ \mathfrak{d}'^{(K+1)}  \\
\equiv \mathfrak{l}'^{(K+1)}\circ (\mathfrak{g}'^{(K)}\circ \mathfrak{f}^{(K)}\circ \mathfrak{g}^{(K)})-
(\mathfrak{g}'^{(K)}\circ \mathfrak{f}^{(K)}\circ \mathfrak{g}^{(K)})\circ \mathfrak{d}'^{(K+1)}\\
+\partial_{\textrm{Hom}((C',d_0'),(D',l_0'))}(e)\otimes u^{K+1} \ \ \ \textrm{mod}(u^{K+2})
\end{gather*}
holds. Moreover, ${\mathfrak{l}'^{(K+1)}\circ (\mathfrak{g}'^{(K)}\circ \mathfrak{f}^{(K)}\circ \mathfrak{g}^{(K)})-
(\mathfrak{g}'^{(K)}\circ \mathfrak{f}^{(K)}\circ \mathfrak{g}^{(K)})\circ \mathfrak{d}'^{(K+1)} \ \ \textrm{mod}(u^{K+2})}$ is divided into a sum of the following three terms:
\begin{enumerate}
\item \begin{gather*}
(\mathfrak{l}'^{(K+1)}\circ \mathfrak{g}'^{(K)}-\mathfrak{g}'^{(K)}\circ \mathfrak{l}^{(K+1)})\circ \mathfrak{f}^{(K)}\circ \mathfrak{g}^{(K)} \\
\equiv (\mathfrak{o}_{K+1}(\mathfrak{g}')\otimes u^{K+1})\circ \mathfrak{f}^{(K)}\circ \mathfrak{g}^{(K)}\equiv (\mathfrak{o}_{K+1}(\mathfrak{g}')f_0g_0)\otimes u^{K+1} \ \ \textrm{mod}(u^{K+2})
\end{gather*}
\item \begin{gather*}
\mathfrak{g}'^{(K)}\circ (\mathfrak{l}^{(K+1)}\circ\mathfrak{f}^{(K)}-\mathfrak{f}^{(K)}\circ \mathfrak{d}^{(K+1)})\circ \mathfrak{g}^{(K)} \\
\equiv \mathfrak{g}'^{(K)}\circ (\mathfrak{o}_{K+1}(\mathfrak{f})\otimes u^{K+1})\circ \mathfrak{g}^{(K)}\equiv (g'_0\mathfrak{o}_{K+1}(\mathfrak{f})g_0)\otimes u^{K+1} \ \ \textrm{mod}(u^{K+2})
\end{gather*}
\item \begin{gather*}
\mathfrak{g}'^{(K)}\circ \mathfrak{f}^{(K)}\circ (\mathfrak{d}^{(K+1)}\circ \mathfrak{g}^{(K)}-\mathfrak{g}^{(K)}\circ \mathfrak{d}'^{(K+1)}) \\
\equiv\mathfrak{g}'^{(K)}\circ \mathfrak{f}^{(K)}\circ (\mathfrak{o}_{K+1}(\mathfrak{g})\otimes u^{K+1})\equiv (g'_0f_0\mathfrak{o}_{K+1}(\mathfrak{g}))\otimes u^{K+1} \ \ \textrm{mod}(u^{K+2}).
\end{gather*}
\end{enumerate}
So, an equality
\begin{gather*}
\mathfrak{o}_{K+1}(\mathfrak{g}'\circ \mathfrak{f}\circ \mathfrak{g})=\mathfrak{o}_{K+1}(\mathfrak{g}')f_0g_0+g'_0\mathfrak{o}_{K+1}(\mathfrak{f})g_0+g'_0f_0\mathfrak{o}_{K+1}(\mathfrak{g})+\partial_{\textrm{Hom}((C',d'_0),(D',l'_0))}(e)
\end{gather*}
holds. Note that ${[\mathfrak{o}_{K+1}(\mathfrak{g})]=[\mathfrak{o}_{K+1}(\mathfrak{g}')]=0}$ because ${\mathfrak{g}}$ and ${\mathfrak{g}'}$ are ${X_{K+1}}$-morphisms.
This implies that 
\begin{gather*}
[\mathfrak{o}_{K+1}(\mathfrak{g}'\circ \mathfrak{f}\circ \mathfrak{g})]=[g'_0\mathfrak{o}_{K+1}(\mathfrak{f})g_0]=(g'_0,g_0)_*[\mathfrak{o}_{K+1}(\mathfrak{f})]
\end{gather*}
holds. Next, we prove ${(\mathrm{iii})}$. Let ${\mathfrak{h}:C\rightarrow D\times [0,1]}$ be an ${X_K}$-homotopy between ${\mathfrak{f}}$ and ${\mathfrak{f}'}$. Then,
\begin{gather*}
[\mathfrak{o}_{K+1}(\mathfrak{f})]=[\mathfrak{o}_{K+1}(\textrm{Eval}_{s=0}\circ \mathfrak{h})]=\big((\textrm{Eval}_{s=0})_0,\textrm{Id}\big)_*[\mathfrak{o}_{K+1}(\mathfrak{h})] \\ 
=\big((\textrm{Eval}_{s=1})_0,\textrm{Id}\big)_*[\mathfrak{o}_{K+1}(\mathfrak{h})]=[\mathfrak{o}_{K+1}(\textrm{Eval}_{s=1}\circ \mathfrak{h})]=[\mathfrak{o}_{K+1}(\mathfrak{f}')]
\end{gather*}
holds.
\begin{flushright}
    $\Box$
\end{flushright}

We have the following corollary.
\begin{Cor}
Let ${\mathfrak{f}:C\rightarrow D}$ be an ${X_{K+1}}$-morphism. Assume that there are ${X_K}$-morphisms ${\mathfrak{g}:C\rightarrow D}$ and ${\mathfrak{h}:C\rightarrow D\times [0,1]}$ such that ${\mathfrak{h}}$ is an ${X_K}$-homotopy between the restriction of ${\mathfrak{f}}$ to an $X_K$-morphism and ${\mathfrak{g}}$. Then, we can extend ${\mathfrak{g}}$ and ${\mathfrak{h}}$ to ${X_{K+1}}$-morphisms so that ${\mathfrak{h}}$ is an ${X_{K+1}}$-homotopy between ${\mathfrak{f}}$ and ${\mathfrak{g}}$.
\end{Cor}
\vspace{5mm}
\textbf{proof}:
We define ${h'_{K+1}:C\rightarrow D\times [0,1]}$ by ${h'_{K+1}=(\textrm{Incl})_0\circ f_{K+1}}$. Then, ${h'_{K+1}}$ satisfies the following equalities:
\begin{gather*}
((\textrm{Eval}_{s=0})_0,\textrm{Id})_*(\mathfrak{o}_{K+1}(\mathfrak{h})+\partial_{\textrm{Hom}((C,d_0),(D\times [0,1],\widetilde{l}_0))}(h'_{K+1})) \\
=\mathfrak{o}_{K+1}(\mathfrak{f})+\partial_{\textrm{Hom}(C,D)}(f_{K+1})=0 
\end{gather*}
\begin{gather*}
(\textrm{Eval}_{s=0})_0\circ h'_{K+1}=f_{K+1}.
\end{gather*}
Let ${N}$ be a subcomplex of ${(D\times [0,1],\widetilde{l}_0)}$ as follows:
\begin{gather*}
N=\textrm{Ker}((\textrm{Eval}_{s=0})_0)=\{(0,y,z)\in D\times [0,1]\}.
\end{gather*}
Then ${\mathfrak{o}_{K+1}(\mathfrak{h})+\partial_{\textrm{Hom}((C,d_0),(D\times [0,1],\widetilde{l}_0))}(h'_{K+1})}$ is a cocycle in ${\textrm{Hom}(C,N)}$. Note that ${N}$ is ${\widetilde{l}_0}$-acyclic. Indeed,
\begin{gather*}
H:N=0\oplus D\oplus D \longrightarrow 0\oplus D\oplus D \\
(0,y,z)\mapsto (0,0,-(-1)^{\textrm{deg}y}y)
\end{gather*}
satisfies ${\textrm{Id}_{N}-0=\widetilde{l}_0H+H\widetilde{l}_0}$. So, 
\begin{gather*}
H\big(\textrm{Hom}(C,N),\partial_{\textrm{Hom}((C,d_0),(N,\widetilde{l}_0))}\big)=0
\end{gather*}
holds and we can choose ${\Delta h_{K+1}:C\rightarrow D\times [0,1]}$ such that 
\begin{gather*}
(\textrm{Eval}_{s=0})_0\circ \Delta h_{K+1}=0
\end{gather*}
\begin{gather*}
\mathfrak{o}_{K+1}(\mathfrak{h})+\partial_{\textrm{Hom}(C,D\times [0,1])}(h'_{K+1})=-\partial_{\textrm{Hom}(C,D\times [0,1])}(\Delta h_{K+1})
\end{gather*}
holds. ${h_{K+1}=h'_{K+1}+\Delta h_{K+1}}$ satisfies the following equalities:
\begin{gather*}
\mathfrak{o}_{K+1}(\mathfrak{h})+\partial_{\textrm{Hom}(C,D\times [0,1])}(h_{K+1})=0
\end{gather*}
\begin{gather*}
(\textrm{Eval}_{s=0})_0\circ h_{K+1}=f_{K+1}.
\end{gather*}
The first equality implies that ${h_0,\cdots,h_K,h_{K+1}}$ determines an ${X_{K+1}}$-morphism. So this is an ${X_{K+1}}$-extension of the original ${X_K}$-morphism ${\mathfrak{h}}$. The second equality implies that the extended ${\mathfrak{h}}$ is an ${X_{K+1}}$-homotopy between ${\mathfrak{f}}$ and ${\textrm{Eval}_{s=1}\circ \mathfrak{h}}$. Note that ${\textrm{Eval}_{s=1}\circ \mathfrak{h}}$ is an ${X_{K+1}}$-extension of ${}\mathfrak{g}$. So it suffices to define an ${X_{K+1}}$-extension of ${\mathfrak{g}}$ by ${\textrm{Eval}_{s=1}\circ \mathfrak{h}}$.
\begin{flushright}
    $\Box$
\end{flushright}
Next, we prove a very important proposition.
\begin{Prop}
Let ${\mathfrak{f}:(C,\mathfrak{d})\rightarrow (D,\mathfrak{l})}$ be an ${X_K}$-morphism and assume that 
\begin{gather*}
f_0:(C,d_0)\longrightarrow (D,l_0)
\end{gather*}
is a cochain homotopy equivalence. Then, ${\mathfrak{f}}$ is an ${X_K}$-homotopy equivalence.
\end{Prop}
\vspace{5mm}
\textbf{proof}:
Assume that ${\mathfrak{g}:D\rightarrow C}$ is an ${X_L}$-morphism (${L<K}$) such that there is an ${X_L}$-homotopy ${\mathfrak{h}:C\rightarrow C\times [0,1]}$ between ${\textrm{Id}}$ and ${\mathfrak{g}\circ \mathfrak{f}}$ (${\mathfrak{g}\circ \mathfrak{f}}$ is a composition of an $X_L$-morphism $\mathfrak{f}$ and the restriction of $\mathfrak{g}$ to an $X_L$-morphism.). Our first purpose is to extend ${\mathfrak{g}}$ and ${\mathfrak{h}}$ to ${X_{L+1}}$-morphisms. Corollary 26 implies that we can choose ${h'_{L+1}:C\rightarrow C\times [0,1]}$ such that ${\mathfrak{h}'=\{h_0,\cdots,h_L,h'_{L+1}\}}$ is an ${X_{L+1}}$-morphism which satisfies the following properties:
\begin{itemize}
\item ${\textrm{Eval}_{s=0}\circ \mathfrak{h}'=\textrm{Id}}$
\item ${\textrm{Eval}_{s=1}\circ \mathfrak{h}'}$ is an ${X_{L+1}}$-extension of ${\mathfrak{g}\circ \mathfrak{f}}$.
\end{itemize}
Note that 
\begin{gather*}
0=[\mathfrak{o}_{L+1}(\mathfrak{g}\circ \mathfrak{f})]=(\textrm{Id},f_0)_*[\mathfrak{o}_{L+1}(\mathfrak{g})]
\end{gather*}
implies that ${[\mathfrak{o}_{L+1}(\mathfrak{g})]=0}$ holds because ${f_0}$ is a cochain homotopy equivalence and hence ${(\textrm{Id},f_0)_*}$ is an isomorphism. So, we can choose ${g'_{L+1}:D\rightarrow C}$ such that ${\mathfrak{g}'=\{g_0,\cdots,g_L,g'_{L+1}\}}$ is an ${X_{L+1}}$-morphism. Then, ${\mathfrak{g}'\circ \mathfrak{f}}$ and ${\textrm{Eval}_{s=1}\circ \mathfrak{h}'}$ are $X_{L+1}$-extensions of ${\mathfrak{g}\circ \mathfrak{f}}$. This implies that the difference of these two extensions
\begin{gather*}
\Theta=(\mathfrak{g}'\circ \mathfrak{f})_{L+1}-(\textrm{Eval}_{s=1})_0\circ h'_{L+1}
\end{gather*}
is a cocycle in ${(\textrm{Hom}(C,C),\partial_{\textrm{Hom}((C,d_0),(C,d_0))})}$ because 
\begin{gather*}
\partial\big((\mathfrak{g}'\circ \mathfrak{f})_{L+1}\big)=\partial((\textrm{Eval}_{s=1})_0\circ h'_{L+1})=-\mathfrak{o}_{L+1}(\mathfrak{g}\circ \mathfrak{f})
\end{gather*}
holds. Recall that
\begin{gather*}
(\textrm{Id},f_0)_*:(\textrm{Hom}(D,C),\partial_{\textrm{Hom}(D,C)})\longrightarrow (\textrm{Hom}(C,C),\partial_{\textrm{Hom}(C,C)})
\end{gather*}
is a cochain homotopy equivalence. So we can choose a cocycle ${\Delta g'_{L+1}\in \textrm{Hom}(D,C)}$ such that 
\begin{gather*}
[\Delta g'_{L+1}\circ f_0+\Theta]=0
\end{gather*}
holds. ${\mathfrak{g}=\{g_0,\cdots,g_L,g'_{L+1}+\Delta g'_{L+1}\}}$ is another ${X_{L+1}}$-extension of ${\mathfrak{g}}$ because
\begin{gather*}
\partial(g'_{L+1}+\Delta g'_{L+1})=\partial(g'_{L+1})=-\mathfrak{o}_{L+1}(\mathfrak{g})
\end{gather*}
holds. We fix ${\Delta_1h_{L+1}\in \textrm{Hom}(C,C)}$ such that 
\begin{gather*}
\partial(\Delta_1h_{L+1})=\Delta g'_{L+1}\circ f_0+\Theta=(\mathfrak{g}\circ \mathfrak{f}-\textrm{Eval}_{s=1}\circ \mathfrak{h}')_{L+1}
\end{gather*}
holds. We also fix ${\Delta h_{L+1}\in \textrm{Hom}(C,C\times [0,1])}$ such that 
\begin{gather*}
(\textrm{Eval}_{s=0})_0\circ \Delta h_{L+1}=0 \\
(\textrm{Eval}_{s=1})_0\circ \Delta h_{L+1}=\Delta_1h_{L+1}
\end{gather*}
holds. Note that ${\mathfrak{h}=\{h_0,\cdots,h_L,h'_{L+1}+\partial(\Delta h_{L+1})\}}$ is another ${X_{L+1}}$-extension of ${\mathfrak{h}}$ because
\begin{gather*}
\partial(h'_{L+1}+\partial(\Delta h_{L+1}))=\partial(h'_{L+1})=-\mathfrak{o}_{L+1}(\mathfrak{h})
\end{gather*}
holds. The equalities
\begin{gather*}
\textrm{Eval}_{s=0}\circ \mathfrak{h}=\textrm{Eval}_{s=0}\circ \mathfrak{h}'=\textrm{Id}
\end{gather*}
\begin{gather*}
(\textrm{Eval}_{s=1}\circ \mathfrak{h})_{L+1}=(\textrm{Eval}_{s=1}\circ \mathfrak{h}')_{L+1}+\partial (\Delta_1h_{L+1}) \\
(\textrm{Eval}_{s=1}\circ \mathfrak{h}')_{L+1}+(\mathfrak{g}\circ \mathfrak{f}-\textrm{Eval}_{s=1}\circ \mathfrak{h}')_{L+1} 
=(\mathfrak{g}\circ \mathfrak{f})_{L+1}
\end{gather*}
imply that ${\mathfrak{h}}$ is an ${X_{L+1}}$-homotopy between ${\textrm{Id}}$ and ${\mathfrak{g}\circ \mathfrak{f}}$. Inductively, we can construct an ${X_K}$-morphism ${\mathfrak{g}:D\rightarrow C}$ so that ${\mathfrak{g}\circ \mathfrak{f}}$ is ${X_K}$-homotopic to the identity. By applying these arguments to this ${\mathfrak{g}}$, we can construct an ${X_{K}}$-morphism ${\mathfrak{f}':C\rightarrow D}$ such that ${\mathfrak{f}'\circ \mathfrak{g}}$ is ${X_K}$-homotopic to the identity. ${\mathfrak{f}}$ is ${X_K}$-homotopic to ${\mathfrak{f}'}$ because  ${\mathfrak{f}}$ is homotopic to ${\mathfrak{f}'\circ \mathfrak{g}\circ \mathfrak{f}}$ and ${\mathfrak{f}'\circ \mathfrak{g}\circ \mathfrak{f}}$ is homotopic to ${\mathfrak{f}'}$. So ${\mathfrak{f}}$ is an ${X_K}$-homotopy equivalence and ${\mathfrak{g}}$ is a homotopy inverse of ${\mathfrak{f}}$.
\begin{flushright}
    $\Box$
\end{flushright}

Next, we explain how to extend ${X_K}$-modules to ${X_{K+1}}$-modules.

\begin{Prop}
Let ${(C,\mathfrak{d})}$ be an ${X_K}$-module. There is an obstruction class
\begin{gather*}
\mathfrak{p}_{K+1}(\mathfrak{d})\in (\textrm{Hom}(C,C), \partial_{\textrm{Hom}((C,d_0),(C,d_0))})
\end{gather*}
such that we can extend ${(C,\mathfrak{d})}$ to an ${X_{K+1}}$-module if and only if ${[\mathfrak{p}_{K+1}(\mathfrak{d})]=0}$ holds. Moreover, if ${\mathfrak{f}:C\rightarrow D}$ is an ${X_K}$-homotopy equivalence, ${[\mathfrak{p}_{K+1}(\mathfrak{d})]=0}$ holds if and only if ${[\mathfrak{p}_{K+1}(\mathfrak{l})]=0}$ holds.
\end{Prop}
\vspace{5mm}
\textbf{proof}:
We define ${\mathfrak{p}_{K+1}(\mathfrak{d})\in \textrm{Hom}(C,C)}$ by 
\begin{gather*}
\mathfrak{d}^{(K)}\circ \mathfrak{d}^{(K)}\equiv \mathfrak{p}_{K+1}(\mathfrak{d})\otimes u^{K+1} \ \ \ \textrm{mod}(u^{K+2}).
\end{gather*}
Then
\begin{gather*}
\partial_{\textrm{Hom}(C,C)}(\mathfrak{p}_{K+1}(\mathfrak{d}))\otimes u^{K+1}
\equiv (d_0\circ \mathfrak{p}_{K+1}(\mathfrak{d})-\mathfrak{p}_{K+1}(\mathfrak{d})\circ d_0)\otimes u^{K+1} \\
\equiv (\mathfrak{d}^{(K)}\circ \mathfrak{p}_{K+1}(\mathfrak{d})-\mathfrak{p}_{K+1}(\mathfrak{d})\circ \mathfrak{d}^{(K)})u^{K+1} \\
\equiv \mathfrak{d}^{(K)}\circ \mathfrak{d}^{(K)}\circ \mathfrak{d}^{(K)}-\mathfrak{d}^{(K)}
\circ \mathfrak{d}^{(K)}\circ \mathfrak{d}^{(K)}\equiv 0 \ \ \ \textrm{mod}(u^{K+2})
\end{gather*}
holds. So, ${\mathfrak{p}_{K+1}(\mathfrak{d})}$ is a cocycle. Note that ${d_{K+1}\in \textrm{Hom}(C,C)}$ satisfies
\begin{gather*}
\partial_{\textrm{Hom}(C,C)}(d_{K+1})=-\mathfrak{p}_{K+1}(\mathfrak{d})
\end{gather*}
if and only if ${\{d_0,\cdots,d_k,d_{K+1}\}}$ determines an ${X_{K+1}}$-module structure on $C$. So, we can extend ${(C,\mathfrak{d})}$ to an ${X_{K+1}}$-module if and only if ${[\mathfrak{p}_{K+1}(\mathfrak{d})]=0}$ holds.
Assume that ${\mathfrak{f}:C\rightarrow D}$ and ${\mathfrak{g}:D\rightarrow C}$ are ${X_K}$-morphisms such that ${\mathfrak{f}\circ \mathfrak{g}}$ and ${\mathfrak{g}\circ \mathfrak{f}}$ are ${X_K}$-homotopic to the identity. It suffices to prove that 
\begin{gather*}
[\mathfrak{p}_{K+1}(\mathfrak{l})]=(f_0,g_0)_*[\mathfrak{p}_{K+1}(\mathfrak{d})]
\end{gather*}
holds. We define ${e:C\rightarrow D}$ by
\begin{gather*}
\mathfrak{f}^{(K)}\circ \mathfrak{d}^{(K)}-\mathfrak{l}^{(K)}\circ \mathfrak{f}^{(K)}\equiv e\otimes u^{K+1}  \ \ \ \textrm{mod}(u^{K+2}).
\end{gather*}
Then
\begin{gather*}
\mathfrak{f}^{(K)}\circ \mathfrak{d}^{(K)} \circ \mathfrak{d}^{(K)}\circ \mathfrak{g}^{(K)} \\
\equiv (eu^{K+1}+\mathfrak{l}^{(K)}\mathfrak{f}^{(K)})  \mathfrak{d}^{(K)} \mathfrak{g}^{(K)}
\equiv ed_0g_0\otimes u^{K+1}+\mathfrak{l}^{(K)}\mathfrak{f}^{(K)}  \mathfrak{d}^{(K)} \mathfrak{g}^{(K)} \\
\equiv  ed_0g_0\otimes u^{K+1}+\mathfrak{l}^{(K)}(eu^{K+1}+\mathfrak{l}^{(K)}\mathfrak{f}^{(K)})
\mathfrak{g}^{(K)}  \\
\equiv (ed_0g_0+l_0eg_0)\otimes u^{K+1}+\mathfrak{p}_{K+1}(\mathfrak{l})f_0g_0\otimes u^{K+1} \\
\equiv (\partial_{\textrm{Hom}(C,D)}(e))g_0\otimes u^{K+1}+\mathfrak{p}_{K+1}(\mathfrak{l})f_0g_0\otimes u^{K+1} \ \ \ \textrm{mod}(u^{K+2})
\end{gather*}
holds. On the other hand,
\begin{gather*}
\mathfrak{f}^{(K)}\circ \mathfrak{d}^{(K)} \circ \mathfrak{d}^{(K)}\circ \mathfrak{g}^{(K)} \\
\equiv \mathfrak{f}^{(K)}\circ (\mathfrak{p}_{K+1}(\mathfrak{d})\otimes u^{K+1}) \circ \mathfrak{g}^{(K)}\equiv (f_0\circ \mathfrak{p}_{K+1}(\mathfrak{d})\circ g_0)\otimes u^{K+1} \ \ \ \textrm{mod}(u^{K+2})
\end{gather*}
holds. This implies that
\begin{gather*}
(f_0,g_0)_*[\mathfrak{p}_{K+1}(\mathfrak{d})]=[f_0\circ \mathfrak{p}_{K+1}(\mathfrak{d})\circ g_0]
=[\mathfrak{p}_{K+1}(\mathfrak{l})f_0g_0]=[\mathfrak{p}_{K+1}(\mathfrak{l})]
\end{gather*}
holds.
\begin{flushright}
    $\Box$
\end{flushright}

We apply the above proposition to the next proposition.
\begin{Prop}
Let ${(C,\mathfrak{d})}$ be an ${X_K}$-module and let ${(D,\mathfrak{l})}$ be an ${X_{K+1}}$-module. Assume that ${\mathfrak{f}:C\rightarrow D}$ is an ${X_K}$-homotopy equivalence from $C$ to the restriction of $D$ to an ${X_K}$-module. Then we can extend ${(C,\mathfrak{d})}$ to an ${X_{K+1}}$-module and we can also extend ${\mathfrak{f}}$ to an ${X_{K+1}}$-homotopy equivalence.
\end{Prop}
\vspace{5mm}
\textbf{proof}:
Proposition 28 implies that ${[\mathfrak{p}_{K+1}(\mathfrak{d})]=0}$ holds. So we can choose ${d'_{K+1}\in \textrm{Hom}(C,C)}$ such that ${\{d_0,\cdots,d_K,d'_{K+1}\}}$ determines an ${X_{K+1}}$-module structure on ${C}$. Our purpose is to extend ${\mathfrak{f}}$ to an ${X_{K+1}}$-morphism. However, ${[\mathfrak{o}_{K+1}(\mathfrak{f})]=0}$ does not hold in general. So, we add a cocycle ${\Delta d_{K+1}\in \textrm{Hom}(C,C)}$ to ${d'_{K+1}}$ so that ${[\mathfrak{o}^{\textrm{new}}_{K+1}(\mathfrak{f})]=0}$ holds. Note that ${\{d_0,\cdots,d_k,d'_{K+1}+\Delta d_{K+1}\}}$ is a new ${X_{K+1}}$-extension of ${\mathfrak{d}}$ and
\begin{gather*}
\mathfrak{o}^{\textrm{new}}_{K+1}(\mathfrak{f})=\mathfrak{o}_{K+1}(\mathfrak{f})-f_0\circ \Delta d_{K+1}
\end{gather*}
holds. We can choose ${\Delta d_{K+1}}$ such that ${[\mathfrak{o}_{K+1}^{new}(\mathfrak{f})]=0}$ holds  because ${f_0:(C,d_0)\rightarrow (D,l_0)}$ is a cochain homotopy equivalence and hence
\begin{gather*}
(f_0,\textrm{Id})_*:H(\textrm{Hom}(C,C),\partial_{\textrm{Hom}(C,C)})\longrightarrow H(\textrm{Hom}(C,D),\partial_{\textrm{Hom}(C,D)})
\end{gather*}
is an isomorphism. Then, ${[\mathfrak{o}^{\textrm{new}}_{K+1}(\mathfrak{f})]=0}$ implies we can extend ${\mathfrak{f}}$ to an ${X_{K+1}}$-morphisms. Proposition 27 implies that this ${\mathfrak{f}}$ is also an ${X_{K+1}}$-homotopy equivalence.
\begin{flushright}
    $\Box$
\end{flushright}

We also consider ${X_K}$-modules and ${X_K}$-morphisms over ${\mathbb{F}_p}$.

\begin{Def}[local ${X_K}$-modules, local ${X_K}$-morphisms]
\begin{enumerate}
\item Let ${\overline{C}}$ be a ${\mathbb{Z}}$-graded finite dimensional ${\mathbb{F}_p}$-vector space. Let ${\overline{\mathfrak{d}}_{loc}=\{\overline{d}_{i,loc}\}_{i=0}^K}$ be a family of maps ${(0\le K \le \infty)}$ such that ${\textrm{deg}(\overline{d}_{i,loc})=1-2i}$ holds. We call ${(\overline{C},\overline{\mathfrak{d}}_{loc})}$ a local $X_K$-module if ${\overline{\mathfrak{d}}_{loc}}$ satisfies
\begin{gather*}
\overline{\mathfrak{d}}_{loc}^{(K)}\circ \overline{\mathfrak{d}}_{loc}^{(K)}\equiv 0 \ \ \ \textrm{mod}(u^{K+1}).
\end{gather*}
If we assume that ${\textrm{deg}(u)=2}$, ${\overline{\mathfrak{d}}^{(K)}=\sum_i \overline{d}_{i,loc}\otimes u^i}$ is a degree ${1}$ map from ${\overline{C}^{(K)}}$ to itself.
\item Let ${(\overline{C},\overline{\mathfrak{d}}_{loc})}$ and ${(\overline{D},\overline{\mathfrak{l}}_{loc})}$ be local ${X_K}$-modules. A local ${X_K}$-morphism between $\overline{C}$ and ${\overline{D}}$ is a family of maps ${\overline{\mathfrak{f}}_{loc}=\{\overline{f}_{i,loc}\}_{i=0}^K}$ which satisfies the following conditions:
\begin{gather*}
\textrm{deg}(\overline{f}_{i,loc})=-2i  \\
\overline{\mathfrak{f}}_{loc}^{(K)}\circ \overline{\mathfrak{d}}_{loc}^{(K)}\equiv \overline{\mathfrak{l}}_{loc}^{(K)}
\circ \overline{\mathfrak{f}}_{loc}^{(K)} \ \ \ \textrm{mod}(u^{K+1}).
\end{gather*}
If we assume that ${\textrm{deg}(u)=2}$, ${\overline{\mathfrak{f}}_{loc}^{(K)}=\sum_i\overline{f}_{i,loc}\otimes u^i}$ is a degree $0$ map from ${\overline{C}^{(K)}}$ to ${\overline{D}^{(K)}}$.
\end{enumerate}
\end{Def}

Next we introduce the ``${\epsilon}$-gapped condition" to ${X_K}$-modules and ${X_K}$-morphisms. We fix ${\epsilon >0}$.

\begin{Def}[$\epsilon$-gapped condition]
\begin{enumerate}
\item Let ${\overline{C}}$ be a ${\mathbb{Z}}$-graded finite dimensional ${\mathbb{F}_p}$-vector space and let ${C}$ be a ${\Lambda_0}$ module defined by
\begin{equation*}
C=\overline{C} \otimes_{\mathbb{F}_p}\Lambda_0 .
\end{equation*}
Assume that ${\overline{\mathfrak{d}}_{loc}=\{\overline{d}_{i,loc}:\overline{C}\rightarrow \overline{C}\}_{i\le K}}$ determines a local ${X_K}$-module structure on ${\overline{C}}$. Let ${d_{i,loc}\in \textrm{Hom}(C,C)}$ be the natural extension of ${\overline{d}_{i,loc}\in \textrm{Hom}(\overline{C},\overline{C})}$. An ${X_K}$-module ${(C,\mathfrak{d})}$ is called ${\epsilon}$-gapped if 
\begin{gather*}
d_i=d_{i,loc}+T^{\epsilon}d_{i,\epsilon} \ \ \ \big(d_{i,\epsilon}\in \textrm{Hom}(C,C)\big)
\end{gather*}
holds.
\item Let ${(\overline{C},\overline{\mathfrak{d}}_{loc})}$ and ${(\overline{D},\overline{\mathfrak{l}}_{loc})}$ be local ${X_K}$-modules. Assume that ${(C,\mathfrak{d})}$ and ${(D,\mathfrak{l})}$ are ${\epsilon}$-gapped ${X_K}$-modules. An ${X_K}$-morphism ${\mathfrak{f}:(C,\mathfrak{d})\rightarrow (D,\mathfrak{l})}$ is called ${\epsilon}$-gapped if there is a local ${X_K}$-morphism ${\overline{\mathfrak{f}}_{loc}:(\overline{C},\overline{\mathfrak{d}}_{loc})\rightarrow (\overline{D},\overline{\mathfrak{l}}_{loc})}$ such that 
\begin{gather*}
f_i=f_{i,loc}+T^{\epsilon}f_{i,\epsilon} \ \ \ \big(f_{i,loc}\in \textrm{Hom}(C,D)\big)
\end{gather*}
holds. Note that ${f_{i,loc}\in \textrm{Hom}(C,D)}$ is the natural extension of ${\overline{f}_{i,loc}\in \textrm{Hom}(\overline{C},\overline{D})}$.
\item Let ${\mathfrak{f}}$ and ${\mathfrak{g}}$ be ${\epsilon}$-gapped ${X_K}$-morphisms between ${(C,\mathfrak{d})}$ and ${(D,\mathfrak{l})}$. We say that ${\mathfrak{f}}$ and ${\mathfrak{g}}$ are ${\epsilon}$-gapped ${X_K}$-homotopic if there is an ${\epsilon}$-gapped ${X_K}$-morphism
\begin{gather*}
\mathfrak{h}:(C,\mathfrak{d})\longrightarrow (D\times [0,1],\widetilde{\mathfrak{l}})
\end{gather*}
such that
\begin{gather*}
\textrm{Eval}_{s=0}\circ \mathfrak{h}=\mathfrak{f}  \\
\textrm{Eval}_{s=1}\circ \mathfrak{h}=\mathfrak{g}
\end{gather*}
holds.
\end{enumerate}
\end{Def}

Note that ${\mathfrak{f}_{loc}:(C,\mathfrak{d}_{loc})\rightarrow (D,\mathfrak{l}_{loc})}$ is an $X_K$-morphism if ${\overline{\mathfrak{f}}_{loc}}$ is an ${X_K}$-morphism. The next proposition is obstruction theory corresponding to the local-to-global extension.

\begin{Prop}
Let ${(C,\mathfrak{d})}$ and ${(D,\mathfrak{l})}$ be ${\epsilon}$-gapped ${X_K}$-modules. Let ${\overline{\mathfrak{f}}_{loc}:(\overline{C},\overline{\mathfrak{d}}_{loc})\rightarrow (\overline{D},\overline{\mathfrak{l}}_{loc})}$ be a local ${X_K}$-morphism.
\begin{enumerate}
\item  There is an obstruction class
\begin{gather*}
\mathfrak{o}_{\epsilon}(\mathfrak{f}_{loc})\in \textrm{Hom}((C^{(K)},\mathfrak{d}^{(K)}),(C^{(K)},\mathfrak{d}^{(K)}))
\end{gather*}
such that we can extend ${\mathfrak{f}_{loc}}$ to ${\epsilon}$-gapped ${X_K}$-morphism if and only if ${[\mathfrak{o}_{\epsilon}(\mathfrak{f}_{loc})]=0}$ holds.
\item Assume that ${\overline{\mathfrak{f}}}_{loc}$ and ${\overline{\mathfrak{f}}}'_{loc}$ are ${X_K}$-homotopic. Then,
\begin{gather*}
[\mathfrak{o}_{\epsilon}(\mathfrak{f}_{loc})]=[\mathfrak{o}_{\epsilon}(\mathfrak{f}'_{loc})]
\end{gather*}
holds.
\item Assume that ${\mathfrak{g}':C'\rightarrow C}$ and ${\mathfrak{g}:D\rightarrow D'}$ are $\epsilon$-gapped ${X_K}$-morphism. Then,
\begin{gather*}
[\mathfrak{o}_{\epsilon}(\mathfrak{g}'_{loc}\circ \mathfrak{f}_{loc}\circ \mathfrak{g}_{loc})]=(\mathfrak{g}^{(K)}_{loc},\mathfrak{g}'^{(K)}_{loc})_*[\mathfrak{o}_{\epsilon}(\mathfrak{f}_{loc})]
\end{gather*}
holds.
\end{enumerate}
\end{Prop}
\vspace{5mm}
\textbf{proof}:
We define ${\mathfrak{o}_{\epsilon}(\mathfrak{f}_{loc})}$ as follows:
\begin{gather*}
\mathfrak{l}^{(K)}\circ \mathfrak{f}^{(K)}_{loc}-\mathfrak{f}^{(K)}_{loc}\circ \mathfrak{d}^{(K)}\equiv T^{\epsilon}\mathfrak{o}_{\epsilon}(\mathfrak{f}_{loc}) \ \ \ \textrm{mod}(u^{K+1}).
\end{gather*}
Note that ${\mathfrak{o}_{\epsilon}(\mathfrak{f}_{loc})\in \textrm{Hom}(C^{(K)},D^{(K)})}$ holds because ${\mathfrak{f}_{loc}:(C,\mathfrak{d}_{loc})\rightarrow (D,\mathfrak{l}_{loc})}$ is an ${}X_K$-morphism and since ${(C,\mathfrak{d})}$ and ${(D,\mathfrak{l})}$ are ${\epsilon}$-gapped. Now, an equality
\begin{gather*}
0=\big\{\partial_{\textrm{Hom}((C^{(K)},\mathfrak{d}^{(K)}),(D^{(K)},\mathfrak{l}^{(K)}))}\big\}^2(\mathfrak{f}^{(K)}_{loc}) \\
=\big\{\partial_{\textrm{Hom}((C^{(K)},\mathfrak{d}^{(K)}),(D^{(K)},\mathfrak{l}^{(K)}))}\big\}(\mathfrak{l}^{(K)}\circ \mathfrak{f}^{(K)}_{loc}-\mathfrak{f}_{loc}^{(K)}\circ \mathfrak{d}^{(K)})
\\
=\big\{\partial_{\textrm{Hom}((C^{(K)},\mathfrak{d}^{(K)}),(D^{(K)},\mathfrak{l}^{(K)}))}\big\}(T^{\epsilon}\mathfrak{o}_{\epsilon}(\mathfrak{f}^{(K)})) \\
=T^{\epsilon} \partial_{\textrm{Hom}((C^{(K)},\mathfrak{d}^{(K)}),(D^{(K)},\mathfrak{l}^{(K)}))}(\mathfrak{o}_{\epsilon}(\mathfrak{f}_{loc}))
\end{gather*}
holds. So, ${\mathfrak{o}_{\epsilon}(\mathfrak{f}_{loc})}$ is a cocycle.
Let ${\mathfrak{f}_{\epsilon}=\{f_{i,\epsilon}\}_{i\le K}}$ be a family of maps in ${\textrm{Hom}(C,D)}$. Then, ${\mathfrak{f}=\{f_{i,loc}+T^{\epsilon}f_{i,\epsilon}\}}$ is an ${\epsilon}$-gapped ${X_K}$-morphism if and only if 
\begin{gather*}
0\equiv \mathfrak{l}^{(K)}\circ \mathfrak{f}^{(K)}-\mathfrak{f}^{(K)}\circ \mathfrak{d}^{(K)}\\
\equiv \mathfrak{l}^{(K)}\circ (\mathfrak{f}_{loc}^{(K)}+T^{\epsilon}\mathfrak{f}_{\epsilon}^{(K)})-(\mathfrak{f}_{loc}^{(K)}+T^{\epsilon}\mathfrak{f}_{\epsilon}^{(K)})\circ \mathfrak{d}^{(K)} \\
\equiv T^{\epsilon}\big\{\mathfrak{o}_{\epsilon}(\mathfrak{f}_{loc})+\partial_{\textrm{Hom}((C^{(K)},\mathfrak{d}^{(K)}),(D^{(K)},\mathfrak{l}^{(K)}))}(\mathfrak{f}_{\epsilon})\big\}
\ \ \ \textrm{mod}(u^{K+1})
\end{gather*}
holds. So, we can extend ${\mathfrak{f}_{loc}}$ to an ${\epsilon}$-gapped ${X_K}$-morphism if and only if ${[\mathfrak{o}_{\epsilon}(\mathfrak{f}_{loc})]=0}$ holds.
Next, we prove ${(\mathrm{iii})}$. By definition,
\begin{gather*}
T^{\epsilon}\cdot \mathfrak{o}_{\epsilon}(\mathfrak{g}'_{loc}\circ \mathfrak{f}_{loc}\circ \mathfrak{g}_{loc})\equiv \mathfrak{l}^{(K)}\circ (\mathfrak{g}'^{(K)}_{loc}\circ \mathfrak{f}^{(K)}_{loc}\circ \mathfrak{g}^{(K)}_{loc})-(\mathfrak{g}'^{(K)}_{loc}\circ \mathfrak{f}^{(K)}_{loc}\circ \mathfrak{g}^{(K)}_{loc})\circ \mathfrak{d}^{(K)} \ \ \textrm{mod}(u^{K+1})
\end{gather*}
holds. The right hand side is a sum of the following three terms:
\begin{enumerate}
\item 
\begin{gather*}
(\mathfrak{l}'^{(K)}\mathfrak{g}^{(K)}_{loc}-\mathfrak{g}^{(K)}_{loc}\mathfrak{l}^{(K)})\mathfrak{f}^{(K)}_{loc}
\mathfrak{g}'^{(K)}_{loc}\equiv (T^{\epsilon}\mathfrak{o}_{\epsilon}(\mathfrak{g}_{loc}))\mathfrak{f}^{(K)}_{loc}
\mathfrak{g}'^{(K)}_{loc}
\end{gather*}
\item
\begin{gather*}
\mathfrak{g}^{(K)}_{loc}(\mathfrak{l}^{(K)}\mathfrak{f}^{(K)}_{loc}-\mathfrak{f}_{loc}^{(K)}\mathfrak{d}^{(K)})\mathfrak{g}'^{(K)}_{loc}\equiv
\mathfrak{g}^{(K)}_{loc}(T^{\epsilon}\mathfrak{o}_{\epsilon}(\mathfrak{f}^{(K)}_{loc})) \mathfrak{g}'^{(K)}_{loc}
\end{gather*}
\item 
\begin{gather*}
\mathfrak{g}_{loc}^{(K)}\mathfrak{f}^{(K)}_{loc}(\mathfrak{d}^{(K)}\mathfrak{g}'^{(K)}_{loc}-\mathfrak{g}'^{(K)}_{loc}\mathfrak{d}'^{(K)})\equiv \mathfrak{g}^{(K)}_{loc}
\mathfrak{f}^{(K)}_{loc}(T^{\epsilon}\mathfrak{o}_{\epsilon}(\mathfrak{g}'_{loc})).
\end{gather*}
\end{enumerate}
Note that ${[\mathfrak{o}_{\epsilon}(\mathfrak{g}_{loc})]=0}$ and ${[\mathfrak{o}_{\epsilon}(\mathfrak{g}'_{loc})]=0}$ hold because ${\mathfrak{g}}$ and ${\mathfrak{g}'}$ are ${\epsilon}$-gapped ${X_K}$-morphisms which extends ${\mathfrak{g}_{loc}}$ and ${\mathfrak{g}'_{loc}}$. This implies that
\begin{gather*}
[\mathfrak{o}_{\epsilon}(\mathfrak{g}'_{loc}\circ \mathfrak{f}_{loc}\circ \mathfrak{g}_{loc})]
=[\mathfrak{g}^{(K)}_{loc}\circ \mathfrak{o}_{\epsilon}(\mathfrak{f}_{loc})\circ \mathfrak{g}'^{(K)}_{loc}]=(\mathfrak{g}^{(K)}_{loc},\mathfrak{g}'^{(K)}_{loc})_*[\mathfrak{o}_{\epsilon}(\mathfrak{f}_{loc})]
\end{gather*}
holds. Next, we prove ${(\mathrm{ii})}$. Let ${\overline{\mathfrak{h}}_{loc}:(\overline{C},\overline{\mathfrak{d}}_{loc})\rightarrow (\overline{D}\times [0,1],\widetilde{\overline{\mathfrak{l}}}_{loc})}$ be an ${X_K}$-homotopy between ${\overline{\mathfrak{f}}_{loc}}$ and ${\overline{\mathfrak{f}}'_{loc}}$. Then, 
\begin{gather*}
[\mathfrak{o}_{\epsilon}(\mathfrak{f}_{loc})]=(\textrm{Eval}_{s=0},\textrm{Id})_*[\mathfrak{h}_{loc}]
=(\textrm{Eval}_{s=1},\textrm{Id})_*[\mathfrak{h}_{loc}]=[\mathfrak{o}_{\epsilon}(\mathfrak{f}'_{loc})]
\end{gather*}
holds.
\begin{flushright}
    $\Box$
\end{flushright}

Next, we prove the following proposition.
\begin{Prop}
Let ${\mathfrak{f}:(C,\mathfrak{d})\rightarrow (D,\mathfrak{l})}$ be an ${\epsilon}$-gapped ${X_K}$-morphism. Assume that ${f_0:(C,d_0)\rightarrow (D,l_0)}$ is a cochain homotopy equivalence. Then we can construct an ${\epsilon}$-gapped ${X_K}$-morphism ${\mathfrak{g}:(D,\mathfrak{l})\rightarrow (C,\mathfrak{d})}$ such that ${\mathfrak{f}\circ \mathfrak{g}}$ and ${\mathfrak{g}\circ \mathfrak{f}}$ are ${\epsilon}$-gapped ${X_K}$-homotopic to the identity.
\end{Prop}
\vspace{5mm}
\textbf{proof}:
Observe that ${f_{0,loc}:(\overline{C},\overline{d}_{0,loc})\rightarrow (\overline{D},\overline{l}_{0,loc})}$ is also a cochain homotopy equivalence. So Proposition 27 implies that we can construct an ${X_K}$-morphism
\begin{gather*}
\mathfrak{g}_{loc}:(\overline{D},\overline{\mathfrak{l}}_{loc})\longrightarrow (\overline{C},\overline{\mathfrak{d}}_{loc})
\end{gather*}
and an ${X_K}$-homotopy 
\begin{gather*}
\overline{\mathfrak{h}}_{loc}:(\overline{C},\overline{\mathfrak{d}}_{loc})\rightarrow (\overline{C}\times [0,1],\widetilde{\overline{\mathfrak{d}}}_{loc})
\end{gather*}
between ${\textrm{Id}}$ and ${\mathfrak{g}_{loc}\circ \mathfrak{f}_{loc}}$. Our purpose is to extend ${\mathfrak{g}_{loc}}$ and ${\mathfrak{h}_{loc}}$ to ${\epsilon}$-gapped ${X_K}$-morphisms so that $\mathfrak{h}$ becomes an ${\epsilon}$-gapped ${X_K}$-homotopy between ${\textrm{Id}}$ and ${\mathfrak{g}\circ \mathfrak{f}}$. The equality
\begin{gather*}
0=[\mathfrak{o}_{\epsilon}(\textrm{Id})]=(\textrm{Eval}_{s=0},\textrm{Id})_*[\mathfrak{o}_{\epsilon}(\mathfrak{h}_{loc})]
\end{gather*}
implies that ${[\mathfrak{o}_{\epsilon}(\mathfrak{h}_{loc})]=0}$ holds. Moreover, note that the image of 
\begin{gather*}
\mathfrak{o}_{\epsilon}(\mathfrak{h}_{loc}):C^{(K)}\longrightarrow (C\times [0,1])^{(K)}=C^{(K)}\oplus C[1]^{(K)}\oplus C^{(K)}
\end{gather*}
is included in ${\textrm{Ker}(\textrm{Eval}_{s=0})=\{(0,y,z)\in  (C\times [0,1])^{(K)} \}}$. Because ${\textrm{Ker}(\textrm{Eval}_{s=0})}$ is acyclic, we can choose ${\mathfrak{h}'_{\epsilon}=\{h'_{i,\epsilon}\}_{0\le i \le K}}$ such that 
\begin{gather*}
h'_{i,\epsilon}:C\longrightarrow C\times [0,1]  \\
\textrm{Eval}_{s=0}(h'_{i,\epsilon})=0 \\
\partial_{\textrm{Hom}(C^{(K)},(C\times [0,1])^{(K)})}(\mathfrak{h}'^{(K)}_{\epsilon})=-\mathfrak{o}_{\epsilon}(\mathfrak{h}_{loc})
\end{gather*}
holds. Then, ${\mathfrak{h}_{loc}+T^{\epsilon}\mathfrak{h}'_{\epsilon}}$ is an ${\epsilon}$-gapped ${X_K}$-morphism which satisfies ${\textrm{Eval}_{s=0}\circ (\mathfrak{h}_{loc}+T^{\epsilon}\mathfrak{h}'_{\epsilon})=\textrm{Id}}$. On the other hand, an equality
\begin{gather*}
0=[\mathfrak{o}_{\epsilon}(\textrm{Id})]=(\textrm{Id},\mathfrak{f}^{(K)}_{loc})_*[\mathfrak{o}_{\epsilon}(\mathfrak{g}_{loc})]
\end{gather*}
implies that ${\mathfrak{o}_{\epsilon}([\mathfrak{g}_{loc}])=0}$ holds. So we can choose ${\mathfrak{g}'_{\epsilon}=\{g'_{i,\epsilon}\}_{0\le i\le K}}$ which satisfies
\begin{gather*}
g'_{i,\epsilon}:D\longrightarrow C \\
\partial_{\textrm{Hom}(D^{(K)},C^{(K)})}(\mathfrak{g}'^{(K)}_{\epsilon})=-\mathfrak{o}_{\epsilon}(\mathfrak{g}_{loc}) .
\end{gather*}
Then ${\mathfrak{g}_{loc}+T^{\epsilon}\mathfrak{g}'_{\epsilon}}$ is an ${\epsilon}$-gapped ${X_K}$-morphism. So ${(\mathfrak{g}_{loc}+T^{\epsilon}\mathfrak{g}'_{\epsilon})\circ \mathfrak{f}}$ and ${\textrm{Eval}_{s=1}\circ (\mathfrak{h}_{loc}+T^{\epsilon}\mathfrak{h}'_{\epsilon})}$ are ${\epsilon}$-gapped ${}X_K$-morphisms which extends ${\mathfrak{g}_{loc}\circ \mathfrak{f}_{loc}}$. Let ${\Theta}$ be the difference of these two ${X_K}$-morphisms. Define
\begin{gather*}
\Theta:C^{(K)}\longrightarrow C^{(K)}  \\
T^{\epsilon}\Theta \stackrel{\textrm{def}}{=}(\mathfrak{g}_{loc}+T^{\epsilon}\mathfrak{g}'_{\epsilon})^{(K)}\circ \mathfrak{f}^{(K)}-
\textrm{Eval}_{s=1}^{(K)}\circ (\mathfrak{h}_{loc}+T^{\epsilon}\mathfrak{h}'_{\epsilon})^{(K)}\\
=T^{\epsilon}\cdot \big\{(\mathfrak{g}^{(K)}_{loc}\circ \mathfrak{f}^{(K)}_{\epsilon}+\mathfrak{g}'^{(K)}_{\epsilon}\circ \mathfrak{f}^{(K)}_{loc})-\textrm{Eval}_{s=1}\circ \mathfrak{h}'^{(K)}_{loc} \big\}.
\end{gather*}
Then ${\Theta}$ is a cocycle because
\begin{gather*}
\partial_{\textrm{Hom}(C^{(K)},C^{(K)})} (\Theta)=\partial (\mathfrak{g}^{(K)}_{loc}\circ \mathfrak{f}^{(K)}_{\epsilon}+\mathfrak{g}'^{(K)}_{\epsilon}\circ \mathfrak{f}^{(K)}_{loc})-\partial(\textrm{Eval}_{s=1}\circ \mathfrak{h}'^{(K)}_{loc} ) \\
=\mathfrak{o}_{\epsilon}(\mathfrak{g}_{loc}\circ \mathfrak{f}_{loc})-\mathfrak{o}_{\epsilon}(\mathfrak{g}_{loc}\circ \mathfrak{f}_{loc})=0
\end{gather*}
holds. Proposition 27 implies that 
\begin{gather*}
\mathfrak{f}^{(K)}: (C^{(K)},\mathfrak{d}^{(K)})\longrightarrow (D^{(K)},\mathfrak{l}^{(K)})
\end{gather*}
is a cochain homotopy equivalence. This implies that
\begin{gather*}
(\textrm{Id},\mathfrak{f}^{(K)})_*:\textrm{Hom}\big((D^{(K)},\mathfrak{l}^{(K)}),(C^{(K)},\mathfrak{d}^{(K)})\big)\longrightarrow \textrm{Hom}\big((C^{(K)},\mathfrak{d}^{(K)}),(C^{(K)},\mathfrak{d}^{(K)})\big)   \\
\phi \mapsto \phi\circ \mathfrak{f}^{(K)}
\end{gather*}
is also a cochain homotopy equivalence. So we can choose ${\Delta \mathfrak{g}'_{\epsilon}=\{\Delta g'_{i,\epsilon}\}_{0\le i\le K}}$ such that
\begin{gather*}
g'_{i,\epsilon}:D\longrightarrow C \\
\partial (\Delta \mathfrak{g}'^{(K)}_{\epsilon})=0 \\
[\Delta \mathfrak{g}'^{(K)}_{\epsilon}\circ \mathfrak{f}^{(K)}]=-[\Theta]
\end{gather*}
holds. Then ${\mathfrak{g}=\mathfrak{g}_{loc}+T^{\epsilon}(\mathfrak{g}'_{\epsilon}+\Delta \mathfrak{g}'_{\epsilon})}$ is an ${\epsilon}$-gapped ${X_K}$-morphism which extends ${\mathfrak{g}_{loc}}$. We choose ${\Delta_1\mathfrak{h}=\{\Delta_1h_{i,\epsilon}\}_{0\le i\le K}}$ such that
\begin{gather*}
\Delta_1h_{i,\epsilon}:C\longrightarrow C \\
\partial_{\textrm{Hom}((C^{(K)},\mathfrak{d}^{(K)}),(C^{(K)},\mathfrak{d}^{(K)}))}(\Delta_1\mathfrak{h}^{(K)}_{\epsilon})=\Theta+\Delta g'^{(K)}_{\epsilon}\circ \mathfrak{f}^{(K)}
\end{gather*}
holds. Next we choose ${\Delta \mathfrak{h}_{\epsilon}=\{\Delta h_{i,\epsilon}\}_{0\le i\le K}}$ such that
\begin{gather*}
\Delta h_{i,\epsilon}:C\longrightarrow C\times [0,1] \\
\textrm{Eval}_{s=0}(\Delta h_{i,\epsilon})=0 \\
\textrm{Eval}_{s=1}(\Delta h_{i,\epsilon})=\Delta_1 h_{i,\epsilon}
\end{gather*}
holds. We define ${\mathfrak{h}=\{h_i\}_{0\le i\le K}}$ as follows:
\begin{gather*}
h_i:C\longrightarrow C\times [0,1] \\
\mathfrak{h}^{(K)}=\mathfrak{h}^{(K)}_{loc}+T^{\epsilon}\big(\mathfrak{h}'^{(K)}_{\epsilon}+\partial (\Delta \mathfrak{h}_{\epsilon}^{(K)})\big).
\end{gather*}
${\mathfrak{h}}$ is an ${\epsilon}$-gapped ${X_K}$-morphism because
\begin{gather*}
\partial \big(\mathfrak{h}'^{(K)}_{\epsilon}+\partial (\Delta \mathfrak{h}_{\epsilon}^{(K)})\big)=\partial (\mathfrak{h}'^{(K)}_{\epsilon})=-\mathfrak{o}_{\epsilon}(\mathfrak{h}_{loc})
\end{gather*}
holds. Moreover, ${\mathfrak{g}\circ \mathfrak{f}=\textrm{Eval}_{s=1}\circ \mathfrak{h}}$ is satisfied because
\begin{gather*}
\mathfrak{g}^{(K)}\circ \mathfrak{f}^{(K)}-\textrm{Eval}^{(K)}_{s=1}\circ \mathfrak{h}^{(K)}
\equiv T^{\epsilon}\{\Theta+\Delta \mathfrak{g}'^{(K)}_{\epsilon}\circ \mathfrak{f}^{(K)}-\partial (\Delta_1\mathfrak{h}^{(K)}_{\epsilon})\}\equiv 0 \ \ \textrm{mod}(u^{K+1})
\end{gather*}
holds. So ${\mathfrak{h}}$ is an ${\epsilon}$-gapped ${X_K}$-homotopy between ${\textrm{Id}}$ and ${\mathfrak{g}\circ \mathfrak{f}}$. By applying these arguments to ${\mathfrak{g}}$, we can construct an ${\epsilon}$-gapped ${X_K}$-morphism ${\mathfrak{f}':C\rightarrow D}$ so that ${\mathfrak{f}'\circ \mathfrak{g}}$ is homotopic to the identity. Then ${\mathfrak{f}}$ is homotopic to ${\mathfrak{f}'}$ because ${\mathfrak{f}}$ is homotopic to ${\mathfrak{f}'\circ \mathfrak{g}\circ \mathfrak{f}}$ and ${\mathfrak{f}'\circ \mathfrak{g}\circ \mathfrak{f}}$ is homotopic to ${\mathfrak{f}'}$. So we proved the proposition.
\begin{flushright}
    $\Box$
\end{flushright}

Next, we explain how to extend ${\epsilon}$-gapped ${X_K}$-modules to ${\epsilon}$-gapped ${X_{K+1}}$-modules.

\begin{Prop}
\begin{enumerate}
\item Let ${(C,\mathfrak{d})}$ be an ${\epsilon}$-gapped $X_K$-module. Assume that there is a map ${\overline{d}_{K+1,loc}:\overline{C}\rightarrow \overline{C}}$ such that ${\{\overline{d}_{0,loc},\cdots, \overline{d}_{K+1,loc}\}}$ determines a local ${X_{K+1}}$-module structure on ${\overline{C}}$. There is an obstruction class
\begin{gather*}
\mathfrak{p}_{K+1,\epsilon}(\mathfrak{d},d_{K+1,loc})\in \textrm{Hom}(C,C)
\end{gather*}
such that there is an ${\epsilon}$-gapped ${X_K}$-module structure ${\{d_0,\cdots,d_{K+1,loc}+T^{\epsilon}d_{K+1,\epsilon}\}}$ on $C$ if and only if ${[\mathfrak{p}_{K+1,\epsilon}(\mathfrak{d},d_{K+1,loc})]=0}$ holds.
\item Let ${\mathfrak{f}:(C,\mathfrak{d})\rightarrow (D,\mathfrak{l})}$ be an ${\epsilon}$-gapped ${X_K}$-morphism which is also an ${X_K}$-homotopy equivalence. Assume that there are maps
\begin{gather*}
\overline{f}_{K+1,loc}:\overline{C}\longrightarrow \overline{D} \\
\overline{d}_{K+1,loc}:\overline{C}\longrightarrow \overline{C} \\
\overline{l}_{K+1,loc}:\overline{D}\rightarrow \overline{D} 
\end{gather*}
such that ${\{\overline{d}_{0,loc},\cdots,\overline{d}_{K+1,loc}\}}$ and ${\{\overline{l}_{0,loc},\cdots,\overline{l}_{K+1,loc}\}}$ are ${X_{K+1}}$-module structures on ${\overline{C}}$ and ${\overline{D}}$, and ${\{\overline{f}_{0,loc},\cdots,\overline{f}_{K+1,loc}\}}$ is an ${X_{K+1}}$-homotopy equivalence between them. Then ${[\mathfrak{p}_{K+1,\epsilon}(\mathfrak{d},d_{K+1,\epsilon})]}=0$ holds if and only if ${[\mathfrak{p}_{K+1,\epsilon}(\mathfrak{l},l_{K+1,\epsilon})]}=0$ holds.
\end{enumerate}
\end{Prop}
\vspace{5mm}
\textbf{proof}:
We define ${\widetilde{d}:C^{(K+1)}\rightarrow C^{(K+1)}}$ as follows:
\begin{gather*}
\widetilde{d}:C^{(K+1)}\longrightarrow C^{(K+1)}  \\
\widetilde{d}(x\otimes 1)=d_0(x)\otimes 1+\cdots+d_K(x)\otimes u^K+d_{K+1,loc}(x)\otimes u^{K+1}.
\end{gather*}
Then 
\begin{gather*}
\widetilde{d}\circ \widetilde{d}=T^{\epsilon}\cdot \mathfrak{p}_{K+1,\epsilon}(\mathfrak{d},d_{K+1,loc})\otimes u^{K+1} \ \ \textrm{mod}(u^{K+2})
\end{gather*}
holds because ${\mathfrak{d}}$ is an ${X_K}$-module structure and ${\overline{\mathfrak{d}}_{loc}}$ is an ${X_{K+1}}$-module structure. ${ \mathfrak{p}_{K+1,\epsilon}(\mathfrak{d},d_{K+1,loc})}$ can be regarded as an element of ${\textrm{Hom}(C,C)}$ because it vanishes on the subset ${u\cdot C^{(K+1)}\subset C^{(K+1)}}$ and its image is included in ${C\subset C^{(K+1)}}$. We calculate
\begin{gather*}
T^{\epsilon}\cdot \partial_{\textrm{Hom}(C,C)}\big(\mathfrak{p}_{K+1,\epsilon}(\mathfrak{d},d_{K+1,\epsilon})\big)\otimes u^{K+1} 
\equiv d_0\circ \widetilde{d}\circ \widetilde{d}-\widetilde{d}\circ \widetilde{d}\circ d_0 \\
\equiv 
\widetilde{d}\circ \widetilde{d}\circ \widetilde{d}-\widetilde{d}\circ \widetilde{d}\circ \widetilde{d}\equiv 0 \ \ \ \textrm{mod}(u^{K+2}).
\end{gather*}
So ${\mathfrak{p}_{K+1,\epsilon}(\mathfrak{d},d_{K+1,loc})}$ is a cocycle. Now ${\{d_0,\cdots,d_{K+1,loc}+T^{\epsilon}d_{K+1,\epsilon}\}}$ is an ${\epsilon}$-gapped ${X_{K+1}}$-module structure on ${C}$ if and only if 
\begin{gather*}
(\widetilde{d}+T^{\epsilon}d_{K+1,\epsilon}\otimes u^{K+1})^2\equiv T^{\epsilon}\big(\mathfrak{p}_{K+1,\epsilon}(\mathfrak{d},d_{K+1,loc})+\partial(d_{K+1,\epsilon})\big)\otimes u^{K+1}\equiv 0 \ \ \ \textrm{mod}(u^{K+2})
\end{gather*}
holds. In other words, ${\{d_0,\cdots,d_{K+1,loc}+T^{\epsilon}d_{K+1,\epsilon}\}}$ is an ${\epsilon}$-gapped ${X_{K+1}}$-module structure on ${C}$ if and only if  ${[\mathfrak{p}_{K+1,\epsilon}(\mathfrak{d},d_{K+1,loc})]=0}$ holds. So we proved ${(\mathrm{i})}$. Next we prove ${(\mathrm{ii})}$. By assumption and Proposition 27 and Proposition 33, we can choose 
\begin{gather*}
\overline{\mathfrak{g}}=\{\overline{g}_{i,loc}\}_{i\le K+1}:(\overline{D},\overline{l})\longrightarrow (\overline{C},\overline{d}) \\
\mathfrak{g}:(D,\mathfrak{l})\longrightarrow (C,\mathfrak{d})
\end{gather*}
such that ${\overline{\mathfrak{g}}}$ is an ${X_{K+1}}$-morphism and ${\mathfrak{g}}$ is an ${\epsilon}$-gapped ${X_K}$-morphism which extends ${\overline{\mathfrak{g}}}$. Moreover, ${\overline{\mathfrak{g}}}$ is a homotopy inverse of ${\overline{\mathfrak{f}}}$ and ${\mathfrak{g}}$ is a homotopy inverse of ${\mathfrak{f}}$. We define ${\widetilde{f}}$, ${\widetilde{g}}$, ${\widetilde{d}}$ and ${\widetilde{l}}$ as follows:
\begin{gather*}
\widetilde{f}=f_0\otimes 1+\cdots +f_K\otimes u^K+f_{K+1,loc}\otimes u^{K+1} \\
\widetilde{g}=g_0\otimes 1+\cdots +g_K\otimes u^K+g_{K+1,loc}\otimes u^{K+1} \\
\widetilde{d}=d_0\otimes 1+\cdots +d_K\otimes u^K+d_{K+1,loc}\otimes u^{K+1} \\
\widetilde{l}=l_0\otimes 1+\cdots +l_K\otimes u^K+l_{K+1,loc}\otimes u^{K+1}.
\end{gather*}
Then 
\begin{gather*}
\widetilde{f}\circ \widetilde{d}-\widetilde{l}\circ \widetilde{f}\equiv T^{\epsilon}\cdot q\otimes u^{K+1} \ \ \ \textrm{mod}(u^{K+2})
\end{gather*}
holds for some ${q\in \textrm{Hom}(C,D)}$ because ${\mathfrak{f}}$ is an ${\epsilon}$-gapped ${X_K}$-morphism and ${\overline{\mathfrak{f}}}$ is an ${X_{K+1}}$-morphism. By the above definition of ${\mathfrak{p}_{K+1,\epsilon}(\cdot)}$ we have the following equalities:
\begin{gather*}
T^{\epsilon}(f_0\circ\mathfrak{p}_{K+1,\epsilon}(\mathfrak{d},d_{K+1,loc})\circ g_0)\otimes u^{K+1}\equiv \widetilde{f}\circ (\widetilde{d})^2\circ \widetilde{g} \\
\equiv (T^{\epsilon}q\otimes u^{K+1}+\widetilde{l}\circ \widetilde{f})\circ \widetilde{d}\circ \widetilde{g} \\
\equiv T^{\epsilon}qd_0g_0\otimes u^{K+1}+\widetilde{l}\circ (T^{\epsilon}q\otimes u^{K+1}+\widetilde{l}\circ \widetilde{f})\circ \widetilde{g}  \\
\equiv T^{\epsilon}(qg_0d_0+l_0qg_0)+T^{\epsilon}\mathfrak{p}_{K+1,\epsilon}(\mathfrak{l},l_{K+1,loc})f_0g_0\otimes u^{K+1} \ \ \ \textrm{mod}(u^{K+2}).
\end{gather*}
So, an equality
\begin{gather*}
f_0\circ \mathfrak{p}_{K+1,\epsilon}(\mathfrak{d},d_{K+1,loc})\circ g_0=\partial (qg_0)+\mathfrak{p}_{K+1,\epsilon}(\mathfrak{l},l_{K+1,loc})f_0g_0
\end{gather*}
holds. This implies that
\begin{gather*}
(f_0,g_0)_*[\mathfrak{p}_{K+1,\epsilon}(\mathfrak{d},d_{K+1,loc})]=[\mathfrak{p}_{K+1,\epsilon}(\mathfrak{l},l_{K+1,loc})f_0g_0]=[\mathfrak{p}_{K+1,\epsilon}(\mathfrak{l},l_{K+1,loc})]
\end{gather*}
holds and ${[\mathfrak{p}_{K+1,\epsilon}(\mathfrak{d},d_{K+1,loc})]=0}$ holds if and only if ${[\mathfrak{p}_{K+1,\epsilon}(\mathfrak{l},l_{K+1,loc})]=0}$ holds.
\begin{flushright}
    $\Box$
\end{flushright}
We apply this to the following proposition.

\begin{Prop}
Let ${(C,\mathfrak{d})}$ be an ${\epsilon}$-gapped ${X_K}$-module and let ${(D,\mathfrak{l})}$ be an ${\epsilon}$-gapped ${X_{K+1}}$ module. Let ${\mathfrak{f}:C\rightarrow D}$ be an ${\epsilon}$-gapped ${X_K}$-homotopy equivalence from $C$ to the restriction of $D$ to an $\epsilon$-gapped $X_K$-module. Assume that we have maps
\begin{gather*}
\overline{d}_{K+1,loc}:\overline{C}\longrightarrow \overline{C} \\
\overline{f}_{K+1,loc}:\overline{C}\longrightarrow \overline{D}
\end{gather*} 
such that ${\{\overline{d}_{0,loc},\cdots,\overline{d}_{K+1,loc}\}}$ is an ${X_{K+1}}$-module structure on ${\overline{C}}$ and ${\{\overline{f}_{0,loc},\cdots,\overline{f}_{K+1,loc}\}}$ is an ${X_{K+1}}$-homotopy equivalence between ${\overline{C}}$ and ${\overline{D}}$. Then we can extend ${(\mathfrak{d},d_{K+1,loc})}$ to an ${\epsilon}$-gapped ${X_{K+1}}$-module structure on ${C}$ and we can extend ${(\mathfrak{f},f_{K+1,loc})}$ to an ${\epsilon}$-gapped ${X_{K+1}}$-homotopy equivalence between ${C}$ and ${D}$.
\end{Prop}
\vspace{5mm}
\textbf{proof}:
By Proposition 34, ${[\mathfrak{p}_{K+1,\epsilon}(\mathfrak{d},d_{K+1,loc})]=0}$ holds and we can choose ${d'_{K+1,\epsilon}:C\rightarrow C}$ so that ${\{d_0,\cdots,d_{K+1,loc}+T^{\epsilon}d'_{K+1,\epsilon}\}}$ becomes an ${\epsilon}$-gapped ${X_{K+1}}$-module structure on $C$. Our next purpose is to extend ${(\mathfrak{f},f_{K+1,loc})}$ to ${\epsilon}$-gapped ${X_{K+1}}$-morphism. We define ${\widetilde{f}}$, $\widetilde{d}$ and ${\widetilde{l}}$ as follows:
\begin{gather*}
\widetilde{f}=f_0\otimes 1+\cdots +f_K\otimes u^K+f_{K+1,loc}\otimes u^{K+1}  \\
\widetilde{d}=d_0\otimes 1+\cdots +d_K\otimes u^K+(d_{K+1,loc}+T^{\epsilon}d'_{K+1,\epsilon})\otimes u^{K+1} \\
\widetilde{l}=l_0\otimes 1+\cdots +l_{K+1}\otimes u^{K+1}.
\end{gather*}
We define ${\mathfrak{o}_{K+1,\epsilon}(\mathfrak{f},f_{K+1,loc}):C\rightarrow D}$ as follows:
\begin{gather*}
\widetilde{l}\circ \widetilde{f}-\widetilde{f}\circ \widetilde{d}\equiv T^{\epsilon}\mathfrak{o}_{K+1,\epsilon}(\mathfrak{f},f_{K+1,loc})\otimes u^{K+1} \ \ \ \textrm{mod}(u^{K+2}).
\end{gather*}
Note that
\begin{gather*}
T^{\epsilon}(\partial(\mathfrak{o}_{K+1,\epsilon}(\mathfrak{f},f_{K+1,loc}))) \equiv \partial_{\textrm{Hom}((C^{(K+1)},\widetilde{d}),(D^{(K+1)},\widetilde{l}))}(\widetilde{l}\circ \widetilde{f}-\widetilde{f}\circ \widetilde{d}) \\
\equiv (\partial_{\textrm{Hom}((C^{(K+1)},\widetilde{d}),(D^{(K+1)},\widetilde{l}))})^2(\widetilde{f})\equiv 0 \ \ \ \textrm{mod}(u^{K+2})
\end{gather*}
implies ${\mathfrak{o}_{K+1,\epsilon}(\mathfrak{f},f_{K+1,loc})}$ is a cocycle. Now ${\{f_0,\cdots,f_K,f_{K+1,loc}+T^{\epsilon}f_{K+1,\epsilon}\}}$ is an ${\epsilon}$-gapped ${X_{K+1}}$-morphism if and only if
\begin{gather*}
0\equiv \widetilde{l}\circ (\widetilde{f}+T^{\epsilon}f_{K+1,\epsilon})- (\widetilde{f}+T^{\epsilon}f_{K+1,\epsilon})\circ \widetilde{d} \\
\equiv T^{\epsilon}(\mathfrak{o}_{K+1,\epsilon}(\mathfrak{f},f_{K+1,loc})+\partial (f_{K+1,\epsilon})) \ \ \ \textrm{mod}(u^{K+2})
\end{gather*}
holds. However, ${[\mathfrak{o}_{K+1,\epsilon}(\mathfrak{f},f_{K+1,loc})]=0}$ does not hold in general. So we add a cocycle ${\Delta d'_{K+1,\epsilon}\in \textrm{Hom}(C,C)}$ to ${d'_{K+1}}$ such that ${[\mathfrak{o}^{\textrm{new}}_{K+1,\epsilon}(\mathfrak{f},f_{K+1,loc})]=0}$ holds. Note that 
\begin{gather*}
\mathfrak{o}^{\textrm{new}}_{K+1,\epsilon}(\mathfrak{f},f_{K+1,loc})=\mathfrak{o}_{K+1,\epsilon}(\mathfrak{f},f_{K+1,loc})-f_0\circ \Delta d'_{K+1,\epsilon}
\end{gather*}
holds. Now
\begin{gather*}
(f_0,\textrm{Id}): (\textrm{Hom}(C,C),\partial_{(\textrm{Hom}(C,C)})\longrightarrow (\textrm{Hom}(C,D),\partial_{(\textrm{Hom}(C,D)})
\end{gather*}
is a cochain homotopy equivalence and we can choose a cocycle ${\Delta d'_{K+1}}$ such that ${[\mathfrak{o}^{\textrm{new}}_{K+1,\epsilon}(\mathfrak{f},f_{K+1,loc})]=0}$ holds. Then, ${\{d_0,\cdots,d_{K,loc}+T^{\epsilon}(d'_{K+1,\epsilon}+\Delta d'_{K+1,\epsilon})\}}$ is an ${\epsilon}$-gapped ${X_K}$-module structure and we can extend ${(\mathfrak{f},f_{K+1,loc})}$ to an ${\epsilon}$-gapped ${X_{K+1}}$-morphism. Proposition 33 implies that this is an ${\epsilon}$-gapped ${X_{K+1}}$-homotopy equivalence.
\begin{flushright}
    $\Box$
\end{flushright}

The following corollary is the ${\epsilon}$-gapped version of Corollary 26.
\begin{Cor}
Let ${(C,\mathfrak{d})}$ and ${(D,\mathfrak{l})}$ be $\epsilon$-gapped ${X_{K+1}}$ modules. Assume that ${\mathfrak{f}:C\rightarrow D}$ is an ${\epsilon}$-gapped ${X_{K+1}}$ morphism and ${\mathfrak{g}:C\rightarrow D}$ is an ${\epsilon}$-gapped ${X_K}$ morphism such that the restriction of ${\mathfrak{f}}$ to an $\epsilon$-gapped $X_K$-morphism and ${\mathfrak{g}}$ are ${\epsilon}$-gapped $X_K$-homotopic. Then we can extend ${\mathfrak{g}}$ to an ${\epsilon}$-gapped ${X_{K+1}}$-morphism such that it is ${\epsilon}$-gapped ${X_{K+1}}$ homotopic to ${\mathfrak{f}}$.
\end{Cor}
\vspace{5mm}
\textbf{proof}:
Let ${\mathfrak{h}:C\rightarrow D\times [0,1]}$ be an ${\epsilon}$-gapped ${X_K}$-homotopy between ${\mathfrak{f}}$ and ${\mathfrak{g}}$. Corollary 26 implies that we can choose ${\overline{h}_{K+1,loc}:\overline{C}\rightarrow \overline{D}\times [0,1]}$ and ${\overline{g}_{K+1,loc}:\overline{C}\rightarrow \overline{D}}$ such that ${\{\overline{g}_{i,loc}\}_{i=0}^{K+1}}$ is a local ${X_{K+1}}$-morphism and ${\{\overline{h}_{0,loc}\}_{i=0}^{K+1}}$ is an ${X_{K+1}}$-homotopy between ${\{\overline{f}_{i,loc}\}_{i=0}^{K+1}}$ and ${\{\overline{g}_{i,loc}\}_{i=0}^{K+1}}$. We consider the obstruction ${\mathfrak{o}_{K+1,\epsilon}(\mathfrak{g},g_{K+1,loc})}$ defined in the proof of Proposition 35. Then, ${h'_{K+1,\epsilon}=(\textrm{Incl})_0\circ f_{K+1,\epsilon}}$ satisfies the following equalities:
\begin{gather*}
((\textrm{Eval}_{s=0})_0,\textrm{Id})_*(\mathfrak{o}_{K+1,\epsilon}(\mathfrak{h},h_{K+1,loc})+\partial_{\textrm{Hom}(C,D\times [0,1])}(h'_{K+1,\epsilon})) \\
=\mathfrak{o}_{K+1,\epsilon}(\mathfrak{f},f_{K+1,loc})+\partial_{\textrm{Hom}(C,D)}(f_{K+1,\epsilon})=0   \\
(\textrm{Eval}_{s=0})_0\circ h'_{K+1,\epsilon}=f_{K+1,\epsilon}.
\end{gather*}
The first equality implies that ${\mathfrak{o}_{K+1,\epsilon}(\mathfrak{h},h_{K+1,loc})+\partial(h'_{K+1,\epsilon})}$ is a cocycle in ${\textrm{Hom}(C,N)}$ where ${N}$ is the kernel of ${(\textrm{Eval}_{s=0})_0:D\times [0,1]\rightarrow D}$. Recall that ${(\textrm{Hom}(C,N),\partial_{\textrm{Hom}(C,N)})}$ is acyclic. So we can find ${\Delta h_{K+1,\epsilon}:C\rightarrow D\times [0,1]}$ such that 
\begin{gather*}
(\textrm{Eval}_{s=0})_0\circ \Delta h_{K+1,\epsilon}=0 \\
\mathfrak{o}_{K+1,\epsilon}(\mathfrak{h},h_{K+1,loc})+\partial_{\textrm{Hom}(C,D\times [0,1])}(h'_{K+1,\epsilon})=-\partial_{\textrm{Hom}(C,D\times [0,1])}(\Delta h_{K+1,\epsilon})
\end{gather*}
holds. So ${h_{K+1,\epsilon}=h'_{K+1,\epsilon}+\Delta h_{K+1,\epsilon}}$ satisfies the following equalities:
\begin{gather*}
\mathfrak{o}_{K+1,\epsilon}(\mathfrak{h},h_{K+1,loc})+\partial_{\textrm{Hom}(C,D\times [0,1])}(h_{K+1,\epsilon})=0 \\
(\textrm{Eval}_{s=0})_0\circ h_{K+1,\epsilon}=f_{K+1,\epsilon}.
\end{gather*}
The first equality implies that ${\mathfrak{h}=\{h_0,\cdots,h_{K+1,loc}+T^{\epsilon}h_{K+1,\epsilon}\}}$ is an ${\epsilon}$-gapped ${X_{K+1}}$-morphism. ${\textrm{Eval}_{s=1}\circ \mathfrak{h}}$ is an ${\epsilon}$-gapped ${X_{K+1}}$-morphism which extends ${\mathfrak{g}}$. The second equality implies that ${\mathfrak{h}}$ is an ${\epsilon}$-gapped ${X_{K+1}}$-homotopy between ${\mathfrak{f}}$ and ${\mathfrak{g}}$.
\begin{flushright}
    $\Box$
\end{flushright}

\subsection{Directed family of $X_K$-modules and $X_K$-morphisms}

First, we consider a family of ${X_K}$-modules and ${X_K}$-morphisms.
\begin{Def}
We call ${\mathcal{C}=\{C_{K,i},\iota_{K,i\rightarrow j},\tau_{K\rightarrow K+1,i}\}}$ a directed family of ${X_K}$-modules if it satisfies the following conditions. Here ${C_{K,i}}$ is an $X_K$-module ${(K\in \mathbb{Z}_{\ge 0}, i\in \mathbb{N})}$.
\begin{enumerate}
\item ${\iota_{K,i\rightarrow j}:C_{K,i}\rightarrow C_{K,j}}$ is an ${X_K}$-homotopy equivalence. Moreover, ${\iota_{K,j\rightarrow l}\circ \iota_{K,i\rightarrow j}}$ is ${X_K}$-homotopic to ${\iota_{K,i\rightarrow l}}$ for any ${i,j,l\in \mathbb{N}}$. Moreover, ${\iota_{K,i\rightarrow i}:C_{K,i}\rightarrow C_{K,i}}$ is equal to the identity for any ${i\in \mathbb{N}}$. In particular, ${\iota_{K,i\rightarrow j}}$ is an $X_K$-homotopy inverse of ${\iota_{K,j\rightarrow i}}$.
\item For any ${K\in \mathbb{Z}_{\ge 0}}$, there is a constant ${N(K)\in \mathbb{N}}$ and a family of ${X_K}$-morphisms as follows:
\begin{gather*}
\tau_{K\rightarrow K+1,i}:C_{K,i}\longrightarrow C_{K+1,i} \ \ \ (\forall i\ge N(K)).
\end{gather*}
Here, ${C_{K+1,i}}$ is the restriction of ${C_{K+1,i}}$ to an $X_K$-module.
\item ${\tau_{K\rightarrow K+1,j}\circ \iota_{K,i\rightarrow j}}$ and ${\iota_{K+1,i\rightarrow j}\circ \tau_{K\rightarrow K+1,i}}$ are ${X_K}$-homotopic for any ${i,j\ge N(K)}$. Here, ${\iota_{K+1,i\rightarrow j}}$ is the restriction of ${\iota_{K+1,i\rightarrow j}}$ to an $X_K$-morphism.
\end{enumerate}
\end{Def}

\begin{Def}
Let ${\mathcal{C}=\{C_{K,i},\iota_{K,i\rightarrow j},\tau_{K\rightarrow K+1,i}\}}$ and ${\mathcal{D}=\{D_{K,i},\iota'_{K,i\rightarrow j},\tau'_{K\rightarrow K+1,i}\}}$ be directed families of ${X_K}$-modules. We call ${\mathcal{F}=\{\mathfrak{f}_{K,i}\}}$ a morphism between ${\mathcal{C}}$ and ${\mathcal{D}}$ if it satisfies the following conditions.
\begin{enumerate}
\item There is a constant ${M(K)\in \mathbb{N}}$ and $X_K$-morphisms ${\mathfrak{f}_{K,i}}$ for every ${K}$ as follows: 
\begin{gather*}
\mathfrak{f}_{K,i}:C_{K,i}\longrightarrow D_{K,i} \ \ \ (\forall i\ge M(K)).
\end{gather*}
\item ${\mathfrak{f}_{K,j}\circ \iota_{K,i\rightarrow j}}$ is ${X_K}$-homotopic to ${\iota'_{K,i\rightarrow j}\circ \mathfrak{f}_{K,i}}$ for any ${i,j\ge M(K)}$.
\item ${\mathfrak{f}_{K+1,i}\circ \tau_{K\rightarrow K+1,i}}$ is ${X_K}$-homotopic to ${\tau'_{K\rightarrow K+1,i}\circ \mathfrak{f}_{K,i}}$ for any ${i\ge \max\{M(K),M(K+1)\}}$. Here, ${\mathfrak{f}_{K+1,i}}$ is the restriction of ${\mathfrak{f}_{K+1,i}}$ to an $X_K$-morphism.
\end{enumerate}
\end{Def}

Let ${\mathcal{F}=\{\mathfrak{f}_{K,i}\}:\mathcal{C}\rightarrow \mathcal{D}}$ be a morphism between directed families of $X_K$-modules ${\mathcal{C}=\{C_{K,i},\iota_{K,i\rightarrow j},\tau_{K\rightarrow K+1,i}\}}$ and ${\mathcal{D}=\{D_{K,i},\iota'_{K,i\rightarrow j},\tau'_{K\rightarrow K+1,i}\}}$. Assume that each ${X_K}$-morphisms ${\mathfrak{f}_{K,i}:C_{K,i}\rightarrow D_{K,i}}$ is an $X_K$-homotopy equivalence. Then, we can construct $X_K$-morphisms ${\mathfrak
{g}_{K,i}:D_{K,i}\rightarrow C_{K,i}}$ for any ${i\ge M(K)}$ such that ${\mathfrak{g}_{K,i}\circ\mathfrak{f}_{K,i}}$ and ${\mathfrak{f}_{K,i}\circ\mathfrak{g}_{K,i}}$ are $X_K$-homotopic to the identity. We want to make sure that ${\mathcal{G}=\{\mathfrak{g}_{K,i}\}}$ becomes a morphism from $\mathcal{D}$ to $\mathcal{C}$. We have the following relations. Here, ${f\sim g}$ represents that $f$ and $g$ are $X_K$-homotopic.

\begin{gather*}
   \iota_{K,i\rightarrow j}\circ\mathfrak{g}_{K,i}\sim (\mathfrak{g}_{K,j}\circ \mathfrak{f}_{K,j})\circ (\iota_{K,i\rightarrow j}\circ\mathfrak{g}_{K,i})  \\
   \sim \mathfrak{g}_{K,j}\circ(\iota'_{K,i\rightarrow j}\circ \mathfrak{f}_{K,i})\circ \mathfrak{g}_{K,i}\sim \mathfrak{g}_{K,j}\circ \iota'_{K,i\rightarrow j}
\end{gather*}

\begin{gather*}
    \tau_{K\rightarrow K+1,i}\circ \mathfrak{g}_{K,i}\sim (\mathfrak{g}_{K+1,i}\circ \mathfrak{f}_{K+1,i})\circ (\tau_{K\rightarrow K+1,i}\circ \mathfrak{g}_{K,i}) \\
    \mathfrak{g}_{K+1,i}\circ (\tau'_{K\rightarrow K+1,i}\circ \mathfrak{f}_{K,i})\circ \mathfrak{g}_{K,i}\sim \mathfrak{g}_{K+1,i}\circ \tau'_{K\rightarrow K+1,i}
\end{gather*}

So, ${\mathfrak{g}_{K,j}\circ \iota'_{K,i\rightarrow j}}$ is $X_K$-homotopic to ${\iota_{K,i\rightarrow j}\circ \mathfrak{g}_{K,i}}$ and ${\mathfrak{g}_{K+1,i}\circ \tau'_{K\rightarrow K+1,i}}$ is $X_K$-homotopic to ${\tau_{K\rightarrow K+1,i}\circ \mathfrak{g}_{K,i}}$. In particular, ${\mathcal{G}=\{\mathfrak{g}_{K,i}\}:\mathcal{D}\rightarrow \mathcal{C}}$ is a morphism between directed families of $X_K$-modules such that each $X_K$-morphisms $\mathfrak{g}_{K,i}$ is $X_K$-homotopy equivalence. This implies that we can define the following equivalence relation between directed families of $X_K$-modules.

\begin{Def}
We call two directed families of ${X_K}$-modules ${\mathcal{C}}$ and ${\mathcal{D}}$ equivalent if there is a morphism ${\mathcal{F}:\mathcal{C}\rightarrow \mathcal{D}}$ such that each ${\mathfrak{f}_{K,i}:C_{K,i}\rightarrow D_{K,i}}$ is an ${X_K}$-homotopy equivalence. We call such $\mathcal{F}$ a homotopy equivalence between directed families of $X_K$-modules. When a morphism ${\mathcal{F}=\{\mathfrak{f}_{K,i}\}:\mathcal{C}\rightarrow \mathcal{D}}$ is homotopy equivalence, we can construct a homotopy equivalence ${\mathcal{G}=\{\mathfrak{g}_{K,i}\}:\mathcal{D}\rightarrow \mathcal{C}}$ such that each $\mathfrak{g}_{K,i}$ is $X_K$-homotopy inverse of $\mathfrak{f}_{K,i}$. We call $\mathcal{G}$ a homotopy inverse of $\mathcal{F}$.
\end{Def}

We also define homotopy of morphisms between directed families of $X_K$-modules.
\begin{Def}
    Let $\mathcal{C}$ and ${\mathcal{D}}$ be directed families of $X_K$-modules. We say two morphisms ${\mathcal{F}=\{\mathfrak{f}_{K,i}\}:\mathcal{C}\rightarrow \mathcal{D}}$ and ${\mathcal{G}=\{\mathfrak{g}_{K,i}\}:\mathcal{C}\rightarrow \mathcal{D}}$ are homotopic if each ${\mathfrak{f}_{K,i}:C_{K,i}\rightarrow D_{K,i}}$ and ${\mathfrak{g}_{K,i}:C_{K,i}\rightarrow D_{K,i}}$ are $X_K$-homotopic (when ${\mathfrak{f}_{K,i}}$ and ${\mathfrak{g}_{K,i}}$ are defined). If a morphism $\mathcal{F}:\mathcal{C}\rightarrow \mathcal{D}$ is a homotopy equivalence, its homotopy inverse ${\mathcal{G}:\mathcal{D}\rightarrow \mathcal{C}}$ is unique up to homotopy. 
\end{Def}

\begin{Def}
Let ${C}$ be an ${X_{\infty}}$-module. Then it defines a natural directed family of ${X_K}$-modules by ${C_{K,i}=C}$ (The right hand side is the restriction of $C$ to an $X_K$-module.). We call this directed family of ${X_K}$-modules constant directed family of ${X_K}$-modules. Note that an $X_{\infty}$-morphism ${\mathfrak{f}:C\rightarrow D}$ between $X_{\infty}$-modules also induces a morphism between constant directed families of $X_K$-modules.
\end{Def}

Let ${\mathfrak{f}, \mathfrak{g}:C\rightarrow D}$ be $X_{\infty}$-morphisms. Assume that $\mathfrak{f}$ and ${\mathfrak{g}}$ are homotopic as morphisms between directed family of $X_K$-modules. Then, $\mathfrak{f}$ and ${\mathfrak{g}}$ are $X_K$-homotopic for any ${K<\infty}$. 

\begin{Prop}
\begin{enumerate}
\item Let ${\mathcal{C}=\{C_{K,i},\iota_{K,i\rightarrow j},\tau_{K\rightarrow K+1,i}\}}$ be a directed family of $X_K$-modules. Then, there is a constant directed family of ${X_K}$-modules ${C}$ and a morphism ${\mathcal{I}:C\rightarrow \mathcal{C}}$ which gives a homotopy equivalence between them. In other words, we can construct an ${X_{\infty}}$-module from ${\mathcal{C}}$. Moreover, this ${X_{\infty}}$-module $C$ is unique up to ${X_{\infty}}$-homotopy equivalence.
\item Let ${\mathcal{C}}$ and ${\mathcal{D}}$ be directed families of ${X_K}$-modules and let ${\mathcal{F}:\mathcal{C}\rightarrow \mathcal{D}}$ be a morphism. We assume that ${C}$ and ${D}$ be ${X_{\infty}}$-modules constructed from ${\mathcal{C}}$ and ${\mathcal{D}}$ and ${\mathcal{I}:C\rightarrow \mathcal{C}}$ and ${\mathcal{J}:D\rightarrow \mathcal{D}}$ be homotopy equivalences constructed in ${(\mathrm{i})}$. Then we can construct an ${X_{\infty}}$-morphism ${\mathfrak{g}:C\rightarrow D}$ such that ${\mathfrak{g}}$ is ${X_K}$-homotopic to the composition of ${C\rightarrow C_{K,i}}$ and ${f_{K,i}:C_{K,i}\rightarrow D_{K,i}}$ and an $X_K$-homotopy inverse of ${D_{K,i}\rightarrow D}$ for any ${K<\infty}$. So, the following diagram of directed families of $X_K$-moduels is commutative up to homotopy.
\[
\begin{tikzcd}
C \arrow[dd, "\mathfrak{g}"'] \arrow[rr, "\mathcal{I}"] &  & \mathcal{C} \arrow[dd, "\mathcal{F}"] \\
                                   &  &                   \\
D \arrow[rr, "\mathcal{J}"]                  &  & \mathcal{D}        
\end{tikzcd}
\]

In other words, $\mathfrak{g}$ is homotopic to the comosition of ${\mathcal{F}\mathcal{I}}$ and the homotopy inverse of $\mathcal{J}$ as a morphism between directed families of $X_K$-modules. In particular, ${\mathfrak{g}}:C\rightarrow D$ is unique up to ${X_K}$-homotopy for any ${K<\infty}$.
\end{enumerate}
\end{Prop}
\vspace{5mm}
\textbf{proof}:
${(\mathrm{i})}$ Let ${\{N(K)\}_{K\in \mathbb{N}}}$ be a sequence of integers such that ${\tau_{K\rightarrow K+1,i}:C_{K,i}\rightarrow C_{K+1,i}}$ is defined for ${i\ge N(K)}$. We fix ${(L,i)}$ and assume that ${C=C_{L,i}}$. So ${C}$ is an ${X_L}$-module. We also fix a constant ${M(L+1)=\max\{i,N(L)\}}$. First, we extend $C$ to an ${X_{L+1}}$-module using the following procedure. We define a family of $X_L$-morphisms ${\mathfrak{l}_{L+1,j}:C\rightarrow C_{L+1,j}}$ for ${j\ge M(L+1)}$ as the composition of ${\iota_{L,i\rightarrow j}:C\rightarrow C_{L,j}}$ and ${\tau_{L\rightarrow L+1,j}:C_{L,j}\rightarrow C_{L+1,j}}$. Then, we have the following diagram which is commutative up to $X_L$-homotopy (${\mathfrak{l}_{L+1,j}}$ is ${X_L}$-homotopic to ${\iota_{L+1,j'\rightarrow j}\circ \mathfrak{l}_{L+1,j'}}$). The following diagram extends infinitely to the right.

\[
\begin{tikzcd}
                  & C \arrow[dd, "\mathfrak{l}_{L+1,M(L+1)+1}" description] \arrow[ldd, "\mathfrak{l}_{L+1,M(L+1)}"'] \arrow[rdd] \arrow[rdd] \arrow[rrdd] &                     &    \\
                  &                                                                                        &                     &    \\
C_{L+1,M(L+1)} \arrow[r] & C_{L+1,M(L+1)+1} \arrow[r]                                                                      & C_{L+1,M(L+1)+2} \arrow[r] & {}
\end{tikzcd}
\]
Each ${\mathfrak{l}_{L+1,j}}$ (${j\ge M(L+1)}$) is an $X_L$-homotopy equivalence from $C$ to the restruction of ${C_{L+1,j}}$ to an $X_L$-module. Proposition 29 implies that we can extend $C$ to an $X_{L+1}$-module and ${\mathfrak{l}_{L+1,M(L+1)}}$ to an ${X_{L+1}}$-homotopy equivalence. Then, ${\mathfrak
{l}_{L+1,M(L+1)+1}}$ is an ${X_L}$-morphism which is $X_L$-homotopic to the restriction of $X_{L+1}$-morphism ${\iota_{L+1,M(L+1)\rightarrow M(L+1)+1}\circ \mathfrak{l}_{L+1,M(L+1)}}$ to an $X_L$-morphism. Corollary 26 implies that we can extend ${\mathfrak{l}_{L+1,M(L+1)+1}}$ to $X_{L+1}$-morphism which is $X_{L+1}$-homotopic to  ${\iota_{L+1,M(L+1)\rightarrow M(L+1)+1}\circ \mathfrak{l}_{L+1,M(L+1)}}$. So the following diagram (the leftmost part of the above diagram) commutes up to $X_{L+1}$-homotopy.

\[
\begin{tikzcd}
                     &  & C \arrow[lldd, "\mathfrak{l}_{L+1,M(L+1)}"'] \arrow[rrdd, "\mathfrak{l}_{L+1,M(L+1)+1}"] &  &   \\
                     &  &                                        &  &   \\
C_{L+1,M(L+1)} \arrow[rrrr, "\iota_{L+1,M(L+1)\rightarrow M(L+1)+1}"'] &  &                                        &  & C_{L+1,M(L+1)+1}
\end{tikzcd}
\]
Note that ${\mathfrak{l}_{L+1,M(L+1)+1}}$ is an ${X_{L+1}}$-homotopy equivalence because ${\iota_{L+1,M(L+1)\rightarrow M(L+1)+1}}$ and ${\mathfrak{l}_{L+1,M(L+1)}}$ are ${X_{L+1}}$-homotopy equivalences. Inductively, we can extend each $X_L$-morphisms ${\mathfrak{l}_{K+1,j}:C\rightarrow C_{L+1,j}}$ to an ${X_{L+1}}$-homotopy equivalence such that the first diagram commutes up to $X_{L+1}$-homotopy.

Next, we fix a constant ${M(L+2)=\max\{M(L+1),N(L+1)\}}$. For any ${j\ge M(L+2)}$, we define an ${X_{L+1}}$-morphism ${\mathfrak{l}_{L+2,j}:C\rightarrow C_{L+2,j}}$ as the compostion of ${\mathfrak{l}_{L+1,j}:C\rightarrow C_{L+1,j}}$ and ${\tau_{L+1\rightarrow L+2,j}:C_{L+1,j}\rightarrow C_{L+2,j}}$. Then we have the following diagram which is commutative up to $X_{L+1}$-homotopy.

\[
\begin{tikzcd}
                  & C \arrow[dd, "\mathfrak{l}_{L+2,M(L+2)+1}" description] \arrow[ldd, "\mathfrak{l}_{L+2,M(L+2)}"'] \arrow[rdd] \arrow[rdd] \arrow[rrdd] &                     &    \\
                  &                                                                                        &                     &    \\
C_{L+2,M(L+2)} \arrow[r] & C_{L+2,M(L+2)+1} \arrow[r]                                                                      & C_{L+2,M(L+2)+2} \arrow[r] & {}
\end{tikzcd}
\]
As in the first step, we can extend $C$ to an ${X_{L+2}}$-module and ${\mathfrak{l}_{L+2,j}}$ to an $X_{L+2}$-homotopy equivalence such that the above diagram commutes up to $X_{L+2}$-homotopy (${\mathfrak{l}_{L+2,j}}$ is $X_{L+2}$ homotopic to ${\iota_{L+2,j'\rightarrow j}\circ \mathfrak{l}_{L+2,j'}}$ for any ${j,j'\ge M(L+2)}$.).

Inductively, we can extend $C$ to an ${X_{\infty}}$-module and we can construct a family of $X_K$-homotopy equivalence ${\mathfrak{l}_{K,j}:C\rightarrow C_{K,j}}$ for any ${j\ge M(K)}$ such that 

\begin{itemize}
    \item $\mathfrak{l}_{K,j}$ is $X_K$ homotopic to ${\iota_{K,j'\rightarrow j}\circ 
    \mathfrak{l}_{L,j'}}$.
    \item $\tau_{K\rightarrow K+1,j}\circ \mathfrak{l}_{K,j}$ is $X_K$-homotopic to the restriction of ${\mathfrak{l}_{K+1,j}}$ to an $X_K$-morphism.
\end{itemize}
holds. So, ${\mathcal{I}=\{\mathfrak{l}_{K,j}\}}$ is a morphis from a constant directed family of $X_K$-modules $C$ to $\mathcal{C}$. Moreover, $\mathcal{I}$ is an homotopy equivalence.

Let ${C'}$ be an $X_{\infty}$-module and ${\mathcal{I}':C'\rightarrow \mathcal{C}}$ be a homotopy equivalence from a constant directed family of $X_K$-modules $C'$ to $\mathcal{C}$. Let ${\mathcal{I}'':\mathcal{C}\rightarrow C'}$ be a homotopy inverse of $\mathcal{I}'$. Then, the composition ${\mathcal{I}''\circ \mathcal{I}:C\rightarrow C'}$ is a homotopy equivalence between constant directed families of $X_K$-modules. In other words, ${\mathcal{I}''\circ \mathcal{I}=\{\mathfrak{l}''_{K,j}:C\rightarrow C'\}}$ is a family of $X_K$-homotopy equivalences defined for $j\ge M''(K)$ for some constant $M''(K)$ and they satisfy the suitable commutativity up to homotopy. We define $X_K$-homotopy equivalence ${\widetilde{\mathfrak{l}}_K:C\rightarrow C'}$ by ${\widetilde{\mathfrak{l}}_K=\mathfrak{l}''_{K,M''(K)}}$ for each ${K\in \mathbb{N}}$. Note that each $\widetilde{l}_K$ is $X_K$-homotopic to the restriction of ${\widetilde{\mathfrak{l}}_{K+1}}$ to an $X_K$-morphism. We fix ${L\in \mathbb{N}}$ and we define an $X_L$-homotopy equivalence ${\mathfrak{l}=\widetilde{\mathfrak{l}}_L}$. Then, Corollary 26 implies that we can extend $\mathfrak{l}$ to $X_{L+1}$-morphism so that $\mathfrak{l}$ is $X_{L+1}$-homotopic to $\widetilde{\mathfrak{l}}_{L+1}$. In particular, $\mathfrak{l}$ is an $X_{L+1}$-homotopy equivalence. Inductively, we can extend $\mathfrak{l}$ to $X_{\infty}$-homotopy equivalence ${\mathfrak{l}:C\rightarrow C'}$. So, $X_{\infty}$-module $C$ is unique up to $X_{\infty}$-homotopy equivalence.

${(\mathrm{ii})}$ Let ${\mathcal{J}':\mathcal{D}\rightarrow D}$ be a homotopy inverse of ${\mathcal{J}:D\rightarrow \mathcal{D}}$. Then, the comosition ${\widetilde{\mathfrak{g}}=\mathcal{J}'\circ \mathcal{F}\circ \mathcal{I}:C\rightarrow D}$ is a morphism between constant directed families of $X_K$-modules. So, ${\widetilde{\mathfrak{g}}=\{\mathfrak{f}_{K,i}:C\rightarrow D\}}$ is a family of $X_K$-morphisms which satisfies the following conditions. 
\begin{itemize}
    \item There is a famiy of integers ${\{M(K)\in \mathbb{N}\}_{K\in \mathbb{N}}}$ such that each $X_K$-morphism ${\mathfrak{f}_{K,i}:C\rightarrow D}$ is defined for ${i\ge M(K)}$.
    \item ${\mathfrak{f}_{K,i}}$ and ${\mathfrak{f}_{K,j}}$ are $X_K$-homotopic.
    \item $\mathfrak{f}_{K,i}$ and the restriction of ${\mathfrak{f}_{K+1,j}}$ to an $X_K$-morphism are $X_K$-homotopic for any ${i\ge M(K)}$ and ${j\ge M(K+1)}$.
\end{itemize}
We fix ${L\in \mathbb{N}}$ and define ${\mathfrak{g}=\mathfrak{f}_{L,M(L)}}$. Then, $\mathfrak{g}$ is an $X_L$-morphism which is $X_L$-homotopic to any ${\mathfrak{f}_{L,i}}$ (${i\ge M(L)}$) and any restriction of ${\mathfrak{f}_{L+1,j}}$ to an $X_L$-morphism (${j\ge M(L+1)}$). Corollary 26 implies that we can extend $\mathfrak{g}$ to an $X_{L+1}$-morphism so that $\mathfrak{g}$ is $X_L$-homotopic to each $\mathfrak{f}_{L+1,j}$. Inductively, we can extend $\mathfrak{g}$ to an $X_{\infty}$-morphism so that $\mathfrak{g}$ is $X_K$-homotopic to each ${\mathfrak{f}_{K,i}}$. In particular, ${\mathfrak{g}}$ is unique up to ${X_K}$-homotopy for any ${K<\infty}$.
\begin{flushright}
    $\Box$
\end{flushright}

\begin{Rem}[Uniqueness of ${\mathfrak{g}:C\rightarrow D}$ constructed in Proposition 42(${\mathrm{ii}}$)]
Let $\mathcal{C}$ and ${\mathcal{D}}$ be directed families of $X_K$-modules. Let ${C_1}$, ${C_2}$, ${D_1}$ and ${D_2}$ be $X_{\infty}$-modules and let ${\mathcal{I}_1:C_1\rightarrow \mathcal{C}}$, ${\mathcal{I}_2:C_2\rightarrow \mathcal{C}}$, ${\mathcal{J}_1:D_1\rightarrow \mathcal{D}}$ and ${\mathcal{J}_2:D_2\rightarrow \mathcal{D}}$ be homotopy equivalences between directed families of $X_K$-modules. Proposition 42(${\mathrm{ii}}$) implies that any morphism ${\mathcal{F}:\mathcal{C}\rightarrow \mathcal{D}}$ (and the identity morphisms ${\mathrm{Id}_{\mathcal{C}}:\mathcal{C}\rightarrow \mathcal{C}}$ and ${\mathrm{Id}_{\mathcal{D}}:\mathcal{D}\rightarrow \mathcal{D}}$) induces $X_{\infty}$-morphisms
\begin{gather*}
    \mathfrak{i}_{12}:C_1\longrightarrow C_2, \ \ \mathfrak{i}_{21}:C_2\longrightarrow C_1   \\
    \mathfrak{j}_{12}:D_1\longrightarrow D_2, \ \ \mathfrak{j}_{21}:D_2\longrightarrow D_1   \\
    \mathfrak{g}_{11}:C_1\longrightarrow D_1, \ \ \mathfrak{g}_{12}:C_1\longrightarrow D_2   \\
    \mathfrak{g}_{21}:C_2\longrightarrow D_1, \ \ \mathfrak{g}_{22}:C_2\longrightarrow D_2
\end{gather*}
such that the following diagram commutes up to homotopy as morphisms between directed families of $X_K$-modules.

\[
\begin{tikzcd}
C_1 \arrow[dddddd, bend right=5, "\mathfrak{i}_{12}"', shift right=2] \arrow[rrrddd, bend right=10,"\mathcal{I}_1"'] \arrow[rrrrrrrrr, "\mathfrak{g}_{11}"] \arrow[rrrrrrrrrdddddd, bend left=5,,"\mathfrak{g}_{12}" near start]  &&&&&&&&& D_1 \arrow[dddddd, bend left=5,"\mathfrak{j}_{12}", shift left=2] \arrow[lllddd, bend left=10,"\mathcal{J}_1"] \\  &&&&&&&&&   \\ &&&&&&&&& \\  &&& \mathcal{C} \arrow[rrr, "\mathcal{F}" near start] &&& \mathcal{D} &&&  \\    &&&&&&&&&   \\ &&&&&&&&&   \\
C_2 \arrow[uuuuuu, bend right=5,"\mathfrak{i}_{21}"'] \arrow[rrruuu, "\mathcal{I}_2"] \arrow[rrrrrrrrr, "\mathfrak{g}_{22}"'] \arrow[rrrrrrrrruuuuuu, "\mathfrak{g}_{21}"' near start] &  &  &                    &  &  &   &  &  & D_2 \arrow[uuuuuu, bend left=5, "\mathfrak{j}_{21}"] \arrow[llluuu,bend right=10, "\mathcal{J}_2"']             
\end{tikzcd}
\]
In particular, we have the following diagram. It is commutative up to $X_K$-homotopy for any ${K<\infty}$ as restrictions of $X_K$-morphisms.

\[
\begin{tikzcd}
C_1 \arrow[ddddd, bend right=5, shift right=1,"\mathfrak{i}_{12}"'] \arrow[rrrrrrr, "\mathfrak{g}_{11}"] \arrow[rrrrrrrddddd, "\mathfrak{g}_{12}" near start]   &&&&&&& D_1 \arrow[ddddd, bend left=5, shift left=1,"\mathfrak{j}_{12}"]     \\ &&&&&&&   \\ &&&&&&& \\ &&&&&&& \\ &&&&&&&  \\ C_2 \arrow[uuuuu, bend right=5, shift right=1,"\mathfrak{i}_{21}"'] \arrow[rrrrrrr, "\mathfrak{g}_{22}"'] \arrow[rrrrrrruuuuu, "\mathfrak{g}_{21}" near start] &&&&&&& D_2 \arrow[uuuuu, "\mathfrak{j}_{21}", bend left=5,shift left=1]
\end{tikzcd}
\]
In this sence, $X_{\infty}$-morphism ${\mathfrak{g}:C\rightarrow D}$ constructed in Proposition 42(${\mathrm{ii}}$) has a uniquness property.
\end{Rem}

We can also define directed families of ${\epsilon}$-gapped ${X_K}$-modules, morphisms between directed families of ${\epsilon}$-gapped ${X_K}$-modules and equivalence by replacing ``${X_K}$-" to ``$\epsilon$-gapped ${X_K}$-". We also have an $\epsilon$-gapped version of Proposition 42 as follows.

\begin{Prop}
\begin{enumerate}
\item Let ${\mathcal{C}=\{C_{K,i},\iota_{K,i\rightarrow j},\tau_{K\rightarrow K+1,i}\}}$ be a directed family of ${\epsilon}$-gapped $X_K$-modules. Then, there is a constant directed family of $\epsilon$-gapped ${X_K}$-modules ${C}$ and a morphism ${\mathcal{I}:C\rightarrow \mathcal{C}}$ which gives a homotopy equivalence between them. In other words, we can construct an ${\epsilon}$-gapped ${X_{\infty}}$-modules from ${\mathcal{C}}$. Moreover, this ${\epsilon}$-gapped ${X_{\infty}}$-module $C$ is unique up to ${\epsilon}$-gapped ${X_{\infty}}$-homotopy equivalence.
\item Let ${\mathcal{C}}$ and ${\mathcal{D}}$ be directed families of ${\epsilon}$-gapped ${X_K}$-modules and let ${\mathcal{F}:\mathcal{C}\rightarrow \mathcal{D}}$ be a morphism. We assume that ${C}$ and ${D}$ be ${\epsilon}$-gapped ${X_{\infty}}$-modules constructed from ${\mathcal{C}}$ and ${\mathcal{D}}$ and ${\mathcal{I}:C\rightarrow \mathcal{C}}$ and ${\mathcal{J}:D\rightarrow \mathcal{D}}$ be homotopy equivalences constructed in ${(\mathrm{i})}$. Then we can construct an ${\epsilon}$-gapped ${X_{\infty}}$-morphism ${\mathfrak{g}:C\rightarrow D}$ such that ${\mathfrak{g}}$ is ${X_K}$-homotopic to the composition of ${C\rightarrow C_{K,i}}$ and ${f_{K,i}:C_{K,i}\rightarrow D_{K,i}}$ and an $X_K$-homotopy inverse of ${D\rightarrow D_{K,i}}$ for any ${K\in \mathbb{Z}_{\ge 0}}$. In particular, ${\mathfrak{g}}:C\rightarrow D$ is unique up to ${X_K}$-homotopy for any ${K<\infty}$. Moreover, the underlying local $X_{\infty}$-morphism ${\overline{\mathfrak{g}}_{loc}:\overline{C}\rightarrow \overline{D}}$ is unique up to ${X_{\infty}}$-homotopy.
\end{enumerate}
\end{Prop}
\vspace{5mm}
\textbf{proof}:
${(\mathrm{i})}$ The proof is the same as the proof of Proposition 42 ${(\mathrm{i})}$. We fix ${(L,i)}$ and assume that ${C=C_{L,i}}$. So $C$ is an ${\epsilon}$-gapped $X_L$-module. First, we apply Proposition 29 to extend ${\overline{C}}$ to a local ${X_{L+1}}$-module such that it is ${X_{L+1}}$-homotopy equivalent to each ${\overline{C}_{L+1,j}}$. Next, we apply Proposition 35 to extend ${C}$ to an ${\epsilon}$-gapped ${X_{L+1}}$-module such that it is $\epsilon$-gapped ${X_{L+1}}$-homotopy equivalent to every ${C_{L+1,j}}$. Inductively, we can extend ${C}$ to an ${\epsilon}$-gapped $X_{\infty}$-module and construct an $\epsilon$-gapped ${X_{\infty}}$-equivalence between constant directed family of ${\epsilon}$-gapped ${X_{\infty}}$-module ${C}$ and ${\mathcal{C}}$. Let ${C'}$ be another $\epsilon$-gapped ${X_{\infty}}$-module which is equivalent to ${\mathcal{C}}$. Then we have a family of $\epsilon$-gapped ${X_K}$-equivalences ${\mathfrak{f}_K:C\rightarrow C'}$ such that ${\mathfrak{f}_K}$ and the restriction of ${\mathfrak{f}_{K+1}}$ to an $\epsilon$-gapped $X_K$-morphism are $\epsilon$-gapped ${X_K}$-homotopic. As in the proof of Proposition 42 (${\mathrm{i}}$), we can construct an ${\epsilon}$-gapped ${X_{\infty}}$-homotopy equivalence between ${C}$ and ${C'}$. In particular, $C$ is unique up to $\epsilon$-gapped $X_{\infty}$-homotopy equivalence.

${(\mathrm{ii})}$ The proof is the same as the proof of Proposition 42 ${(\mathrm{ii})}$. Let ${\mathcal{G}:C\rightarrow D}$ be the compositon of ${\mathcal{I}:C\rightarrow \mathcal{C}}$, ${\mathcal{F}:\mathcal{C}\rightarrow \mathcal{D}}$ and the homotopy inverse of ${\mathcal{J}:D\rightarrow \mathcal{D}}$. Then, we can choose a family of ${\epsilon}$-gapped ${X_K}$-morphisms ${\mathfrak{g}_K}$ such that ${\mathfrak{g}_K}$ is ${\epsilon}$-gapped ${X_K}$-homotopic to the restriction of ${\mathfrak{g}_{K+1}}$ as an $\epsilon$-gapped $X_K$-morphism. Our purpose is to construct an ${\epsilon}$-gapped ${X_{\infty}}$-morphism ${\mathfrak{g}}$ from this family. Corollary 36 implies that we can construct an ${\epsilon}$-gapped ${X_{\infty}}$-morphism ${\mathfrak{g}:C\rightarrow D}$ such that ${\mathfrak{g}}$ is ${\epsilon}$-gapped ${X_K}$-homotopic to each ${\mathfrak{g}_K}$. Next we prove the uniqueness of the underlying local ${X_{\infty}}$-morphism ${\overline{\mathfrak{g}}_{loc}:\overline{C}\rightarrow \overline{D}}$ up to ${X_{\infty}}$-homotopy. We see that two local $X_{\infty}$-morphisms are $X_{\infty}$-homotopic if and only if they are $X_K$-homotopic as restrictions of $X_K$-morphisms for any ${K<\infty}$. Let ${\mathfrak{g}':C\rightarrow D}$ be another ${\epsilon}$-gapped ${X_{\infty}}$-morphism which is ${\epsilon}$-gapped ${X_K}$-hotomopic to ${\mathfrak{g}}$ for any ${K<\infty}$. We fix a large positive integer ${S\in \mathbb{N}}$ so that 
\begin{gather}
|\textrm{deg}(x)|<\frac{1}{10}S   \tag{$*$}
\end{gather}
for any ${x\in \overline{C}\cup \overline{D}}$. Let ${\overline{\mathfrak{h}}_{loc}=\{\overline{h}_{j,loc}\}_{j=0}^S}$ be a family of maps ${\overline{h}_{j,loc}:\overline{C}\rightarrow \overline{D}}$ so that 
\begin{gather*}
\overline{g}_{loc}^{(S)}-\overline{g}'^{(S)}_{loc}\equiv
\overline{\mathfrak{h}}^{(S)}_{loc}\circ \overline{\mathfrak{d}}^{(S)}_{loc}+\overline{\mathfrak{l}}^{(S)}_{loc}\circ \overline{\mathfrak{h}}^{(S)}_{loc}
 \ \ \ \textrm{mod}(u^{S+1})
\end{gather*}
holds. This is possible because ${\overline{\mathfrak{g}}_{loc}}$ and ${\overline{\mathfrak{g}}'_{loc}}$ are ${X_S}$-homotopic. Then the condition ${(*)}$ implies that 
\begin{gather*}
\overline{g}^{(\infty)}_{loc}-\overline{g}'^{(\infty)}_{loc}=\overline{\mathfrak{h}}^{(\infty)}_{loc}\circ \overline{\mathfrak{d}}^{(\infty)}_{loc}+\overline{\mathfrak{l}}^{(\infty)}_{loc}\circ \overline{\mathfrak{h}}^{(\infty)}_{loc}
\end{gather*} 
holds. In particular, ${\overline{\mathfrak{g}}_{loc}}$ is ${X_{\infty}}$-homotopic to ${\overline{\mathfrak{g}}'_{loc}}$.
\begin{flushright}
    $\Box$
\end{flushright}

\begin{Rem}[Uniqueness of $\epsilon$-gapped $X_{\infty}$-morphism ${\mathfrak{g}:C\rightarrow D}$ ]
 An $\epsilon$-gapped $X_{\infty}$-morphism ${\mathfrak{g}:C\rightarrow D}$ constructed in Proposition 44(${\mathrm{ii}}$) has a uniqueness property as we explained in Remark 43. All we have to do is just add the word ``$\epsilon$-gapped". However, the underlying local $X_{\infty}$-morphism ${\overline{\mathfrak{g}}_{loc}:\overline{C}\rightarrow \overline{D}}$ has a stronger uniqueness property. Let $\mathcal{C}$ and ${\mathcal{D}}$ be directed families of $\epsilon$-gapped $X_K$-modules. Let ${C_1}$, ${C_2}$, ${D_1}$ and ${D_2}$ be $\epsilon$-gapped $X_{\infty}$-modules and let ${\mathcal{I}_1:C_1\rightarrow \mathcal{C}}$, ${\mathcal{I}_2:C_2\rightarrow \mathcal{C}}$, ${\mathcal{J}_1:D_1\rightarrow \mathcal{D}}$ and ${\mathcal{J}_2:D_2\rightarrow \mathcal{D}}$ be homotopy equivalences between directed families of $\epsilon$-gapped $X_K$-modules. As we explained in Remark 43, we get the following diagram.
 
 \[
\begin{tikzcd}
C_1 \arrow[ddddd, bend right=5, shift right=1,"\mathfrak{i}^{12}"'] \arrow[rrrrrrr, "\mathfrak{g}^{11}"] \arrow[rrrrrrrddddd, "\mathfrak{g}^{12}" near start]   &&&&&&& D_1 \arrow[ddddd, bend left=5, shift left=1,"\mathfrak{j}^{12}"]     \\ &&&&&&&   \\ &&&&&&& \\ &&&&&&& \\ &&&&&&&  \\ C_2 \arrow[uuuuu, bend right=5, shift right=1,"\mathfrak{i}^{21}"'] \arrow[rrrrrrr, "\mathfrak{g}^{22}"'] \arrow[rrrrrrruuuuu, "\mathfrak{g}^{21}" near start] &&&&&&& D_2 \arrow[uuuuu, "\mathfrak{j}^{21}", bend left=5,shift left=1]
\end{tikzcd}
\]

\end{Rem}
The diagram commutes up to $\epsilon$-gapped $X_{K}$-homotopy for any ${K<\infty}$. In particular, the diagram for underlying local $X_{\infty}$-morphisms also commutes up to $X_K$-homotopy for any ${K<\infty}$.

\[
\begin{tikzcd}
\overline{C}_1 \arrow[ddddd, bend right=5, shift right=1,"\overline{\mathfrak{i}}^{12}_{loc}"'] \arrow[rrrrrrr, "\overline{\mathfrak{g}}^{11}_{loc}"] \arrow[rrrrrrrddddd, "\overline{\mathfrak{g}}^{12}_{loc}" near start]   &&&&&&& \overline{D}_1 \arrow[ddddd, bend left=5, shift left=1,"\overline{\mathfrak{j}}^{12}_{loc}"]     \\ &&&&&&&   \\ &&&&&&& \\ &&&&&&& \\ &&&&&&&  \\ \overline{C}_2 \arrow[uuuuu, bend right=5, shift right=1,"\overline{\mathfrak{i}}^{21}_{loc}"'] \arrow[rrrrrrr, "\overline{\mathfrak{g}}^{22}_{loc}"'] \arrow[rrrrrrruuuuu, "\overline{\mathfrak{g}}^{21}_{loc}" near start] &&&&&&& \overline{D}_2 \arrow[uuuuu, "\overline{\mathfrak{j}}^{21}_{loc}", bend left=5,shift left=1]
\end{tikzcd}
\]

This implies that this diagram commutes up to $X_{\infty}$-homotopy because two local $X_{\infty}$-morphisms are $X_{\infty}$-homotopic if and only if they are $X_K$-homotopic for any ${K<\infty}$ (see the proof of Proposition 44(${\mathrm{ii}}$)).

\begin{Def}
Let ${(C,\mathfrak{d}=\{d_i\}_{i=0}^{\infty})}$ be an ($\epsilon$-gapped) ${X_{\infty}}$-module. Then, ${(C^{(\infty)},\mathfrak{d}^{(\infty)})}$ is a cochain complex (Here, $C^{(\infty)}$ is a $\Lambda_0$-module ${C^{(\infty)}=C\otimes\Lambda_0[[u]]}$ and $\mathfrak{d}^{(\infty)}$ is an infinite sum ${\mathfrak{d}^{(\infty)}=d_0+d_1u+d_2u^2+d_3u^3+\cdots}$.). We define the cohomology of an ($\epsilon$-gapped) ${X_{\infty}}$-module ${(C,\mathfrak{d})}$ to be the cohomology of  ${(C^{(\infty)},\mathfrak{d}^{(\infty)})}$ and denote it by ${H(C,\mathfrak{d})}$. We also define ${\widehat{H}(C,\mathfrak{d})}$ to be the cohomology of the cochain complex ${(C\otimes\Lambda_0[u^{-1},u]],\mathfrak{d}^{(\infty)})}$.
\end{Def}

If ${(C,\mathfrak{d})}$ is an ${\epsilon}$-gapped ${X_{\infty}}$-module, we can also consider the cohomologies of the underlying local $X_{\infty}$-module  ${H(\overline{C},\overline{\mathfrak{d}}_{loc})}$ and ${\widehat{H}(\overline{C},\overline{\mathfrak{d}}_{loc})}$. So, ${H(\overline{C},\overline{\mathfrak{d}}_{loc})}$ is the cohomology of ${(\overline{C}\otimes \mathbb{F}_p[[u]],\overline{\mathfrak{d}}^{(\infty)}_{loc})}$ and ${\widehat{H}(\overline{C},\overline{\mathfrak{d}}_{loc})}$ is the cohomology of ${(\overline{C}\otimes \mathbb{F}_p[u^{-1},u]],\overline{\mathfrak{d}}^{(\infty)}_{loc})}$.

Let ${\mathfrak{f}:(C,\mathfrak{d})\rightarrow(D,\mathfrak{l})}$ be an (${\epsilon}$-gapped) $X_{\infty}$-homotopy equivalence between ($\epsilon$-gapped) $X_{\infty}$-modules. Then, the cochain map
\begin{gather*}
    \mathfrak{f}^{(\infty)}:(C^{(\infty)},\mathfrak{d}^{(\infty)})\longrightarrow (D^{(\infty)},\mathfrak{l}^{(\infty)})
\end{gather*}
is a cochain homotopy equivalence. In particular, $\mathfrak{f}$ induces an isomorphism between cohomologies
\begin{gather*}
    \mathfrak{f}_*:H(C,\mathfrak{d})\longrightarrow H(D,\mathfrak{l}) 
\end{gather*}
\begin{gather*}
    \widehat{\mathfrak{f}}_*:\widehat{H}(C,\mathfrak{d})\longrightarrow \widehat{H}(D,\mathfrak{l}).
\end{gather*}

Let ${\mathcal{C}=\{C_{K,i},\iota_{K,i\rightarrow j},\tau_{K\rightarrow K+1,j}\}}$ be a directed family of ($\epsilon$-gapped) $X_{K}$-modules. Proposition 44(i) implies that we can construct an ($\epsilon$-gapped) $X_{\infty}$ module ${(C,\mathfrak{d})}$ from ${\mathcal{C}}$ and a homotopy equivalence between directed family of ($\epsilon$-gapped) $X_K$-modules
\begin{gather*}
    \mathcal{I}:C\longrightarrow \mathcal{C}.
\end{gather*}
Such ${(C,\mathfrak{d})}$ is unique up to ($\epsilon$-gapped) ${X_{\infty}}$-homotopy equivariance. So, ${H(C,\mathfrak{d})}$ and ${\widehat{H}(C,\mathfrak{d})}$ do not depend on the choice of ${(C,\mathfrak{d})}$. In particular, we can give the following definition.

\begin{Def}
    We define the cohomologies of a directed family of ($\epsilon$-gapped) $X_{K}$-modules ${\mathcal{C}=\{C_{K,i},\iota_{K,i\rightarrow j},\tau_{K\rightarrow K+1,j}\}}$ by
\begin{gather*}
    H(\mathcal{C})=H(C,\mathfrak{d}) 
\end{gather*}
\begin{gather*}
    \widehat{H}(\mathcal{C})=\widehat{H}(C,\mathfrak{d}).
\end{gather*}
If $\mathcal{C}$ is a directed family of ${\epsilon}$-gapped $X_K$-modules, we can define the cohomology of the underlying directed family of local $X_K$-modules ${\overline{C}}$ as follows:
\begin{gather*}
    H_{loc}(\mathcal{C})=H(\overline{C},\overline{\mathfrak{d}}_{loc})
\end{gather*}
\begin{gather*}
 \widehat{H}_{loc}(\mathcal{C})=\widehat{H}(\overline{C},\overline{\mathfrak{d}}_{loc}).
\end{gather*}
\end{Def}

Let ${\mathcal{C}}$ and ${\mathcal{D}}$ be directed families of ($\epsilon$-gapped) $X_K$-modules and let ${\mathcal{F}:\mathcal{C}\rightarrow \mathcal{D}}$ be a morphism between them. Let $C$ and $D$ be ($\epsilon$-gapped) $X_{\infty}$-modules and let ${\mathcal{I}:C\rightarrow \mathcal{C}}$ and ${\mathcal{J}:D\rightarrow \mathcal{D}}$ be homotopy equivalences constructed in Proposition 44(${\mathrm{i}}$). Then, we can construct an ($\epsilon$-gapped) $X_{\infty}$-morphism ${\mathfrak{g}:C\rightarrow D}$ such that the diagram

\[
\begin{tikzcd}
C \arrow[rrrr, "\mathcal{I}"] \arrow[ddd, "\mathfrak{g}"'] &  &  &  & \mathcal{C} \arrow[ddd, "\mathcal{F}"] \\
                                      &  &  &  &                     \\
                                      &  &  &  &                     \\ \ D \arrow[rrrr, "\mathcal{J}"]                   &  &  &  & \mathcal{D}               
\end{tikzcd}
\]
commutes up to homotopy (as morphisms between directed families). Such $\mathfrak{g}$ induces cochain maps ${(C^{(\infty)},\mathfrak{d}^{\infty})\rightarrow (D^{(\infty)},\mathfrak{l}^{(\infty)})}$ and ${(\overline{C}^{(\infty)},\overline{\mathfrak{d}}_{loc}^{\infty})\rightarrow (\overline{D}^{(\infty)},\overline{\mathfrak{l}}_{loc}^{(\infty)})}$. So we can define the following maps between homologies. We denote them by ${\mathcal{F}_*^\mathfrak{g}}$, ${\widehat{\mathcal{F}}_*^\mathfrak{g}}$, ${\widehat{\mathcal{F}}_{loc,*}^{\mathfrak{g}}}$ and ${\widehat{\mathcal{F}}_{loc,*}^\mathfrak{g}}$ in order to make it clear that they depend on the choice of $\mathfrak{g}$.

\begin{gather*}
    \mathcal{F}_*^{\mathfrak{g}}:H(\mathcal{C})\longrightarrow H(\mathcal{D})
\end{gather*}
\begin{gather*}
    \widehat{\mathcal{F}}_*^{\mathfrak{g}}:\widehat{H}(\mathcal{C})\longrightarrow \widehat{H}(\mathcal{D})
\end{gather*}
\begin{gather*}
    \mathcal{F}_{loc,*}^{\mathfrak{g}}:H_{loc}(\mathcal{C})\longrightarrow H_{loc}(\mathcal{D})
\end{gather*}
\begin{gather*}
    \widehat{\mathcal{F}}_{loc,*}^{\mathfrak{g}}:\widehat{H}_{loc}(\mathcal{C})\longrightarrow \widehat{H}_{loc}(\mathcal{D})
\end{gather*}
If ${\mathcal{F}:\mathcal{C}\rightarrow \mathcal{D}}$ is a homotopy equivalence, Proposition 33 implies that the induced ($\epsilon$-gapped) $X_{\infty}$-morphism ${\mathfrak{g}:C\rightarrow D}$ is an ($\epsilon$-gapped) ${X_{\infty}}$-homotopy equivalence because ${g_0:(C,d_0)\rightarrow (D,l_0)}$ is a cochain homogopy equivalence. So, ${H(\mathcal{C})}$, ${\widehat{H}(\mathcal{C})}$, ${H_{loc}(\mathcal{C})}$ and ${\widehat{H}_{loc}(\mathcal{C})}$ is determined by the homotopy equivalence class of $\mathcal{C}$.

Note that ${\mathcal{F}_{loc,*}^{\mathfrak{g}}}$ and ${\widehat{\mathcal{F}}_{loc,*}^{\mathfrak{g}}}$ do not depend on the choice of $\mathfrak{g}$ because ${\overline{\mathfrak{g}}_{loc}}$ is unique up to $X_{\infty}$-homotopy (Proposition 44(${\mathrm{ii}}$)). This uniqueness is sufficient for our purpose (see Remark 54).

\subsection{Geometric constructions}
In this subsection, we construct equivariant Floer cohomology on weakly monotone symplectic manifolds by applying $X_K$-modules and $X_K$-morphisms. Let ${(M,\omega)}$ be a closed weakly monotone symplectic manifold. We fix a homology class ${\gamma \in H_1(M:\mathbb{Z})/\mathrm{Tor}}$. Let ${H\in C^{\infty}(S^1\times M)}$ be a Hamiltonian function such that the number of $1$-periodic orbits of ${\phi_H}$ in $\gamma$ is finite. We choose a family of Hamiltonian functions ${\{H_k\in C^{\infty}(S^1\times M)\}_{k=1}^{\infty}}$ which satisfies the following conditions:

\begin{itemize}
\item ${H_k\longrightarrow H}$ in ${C^{\infty}}$-topology
\item ${(H_k,J)}$ is a Floer regular pair
\end{itemize}
Assume that ${P(H,\gamma)}$ is a finite set ${\{x_1,\cdots,x_l\}}$. For fixed ${1\le i\le l}$, $x_i$ splits into ${\{x_i^1,\cdots,x_i^{l_{i(k)}}\}\subset P(H_k,\gamma)}$. As in section 4, we slightly modify the Floer coboundary operator of ${(H_k,J)}$ as follows. 

Let ${v_i^j:[0,1]\times S^1\rightarrow M}$ be a small cylinder connecting ${v_i^j(0,t)=x_i(t)}$ and ${v_i^j(1,t)=x_i^j(t)}$. Let ${c(x_i,x_i^j)\in \mathbb{R}}$ be the action gap
\begin{equation*}
c(x_i,x_i^j)=\int_{[0,1]\times S^1}(v_i^j)^*\omega+\int_0^1H(t,x_i(t))-H_k(t,x_i^j(t))dt .
\end{equation*}
and let $\tau$ be a correction map defined as follows:
\begin{gather*}
\tau :CH(H_k,\gamma:\Lambda)\longrightarrow CH(H_k,\gamma:\Lambda)  \\
x_i^j\mapsto T^{c(x_i,x_i^j)}x_i^j.
\end{gather*}
The modified coboundary operator ${\widetilde{d_F}}$ was defined by
\begin{equation*}
\widetilde{d_F}=\tau^{-1}d_F\tau.
\end{equation*}
By using the modified Floer coboundary operator ${\widetilde{d_F}}$, we can define a modified coboundary operator 
\begin{equation*}
d_{\mathbb{Z}_p}:CF(H_k,\gamma:\Lambda_0)^{\otimes p}\otimes \Lambda_0[[u]]\langle \theta \rangle \longrightarrow CF(H_k,\gamma:\Lambda_0)^{\otimes p}\otimes \Lambda_0[[u]]\langle \theta \rangle
\end{equation*}
as follows:
\begin{gather*}
d_{\mathbb{Z}_p}(x\otimes 1)=\widetilde{d_F}(x)\otimes 1+(1-\tau)\otimes \theta \\
d_{\mathbb{Z}_p}(x\otimes \theta)=-\widetilde{d_F}(x)\otimes \theta +N(x)\otimes u\theta.
\end{gather*}
Then ${d_{\mathbb{Z}_p}}$ determines an ${X_{\infty}}$-module structure on 
\begin{equation*}
C_k=CF(H_k,\gamma:\Lambda_0)^{\otimes p}\otimes \Lambda_0\langle \theta \rangle .
\end{equation*}
We define a degree of $C_k$ as follows. For any periodic orbit ${x\in P(H_k,\gamma)}$, we can define the Conley-Zehnder index ${\mu_{CZ}(x)\in \mathbb{Z}_2}$ (we normalize ${\mu_{CZ}}$ so that the Conley-Zehnder index of a local maximum of a $C^2$-small Morse function is equal to $n$). We define the degree of $x$ by ${\textrm{deg}(x)=n-\mu_{CZ}(x)}$. The degree of ${\theta}$ is equal to $1$. Moreover, we fix a trivialization of ${x^*TM}$ so that ${\mu_{CZ}(x)\in \mathbb{Z}}$ and ${\textrm{deg}(x)\in \mathbb{Z}}$ are well-defined. Then, ${C_k}$ becomes an ${\epsilon}$-gapped $X_{\infty}$-module.

Note that ${(C_k,d_{\mathbb{Z}_p})}$ and ${(C_{k'},d_{\mathbb{Z}_p})}$ are chain homotopy equivalent (hence ${X_{\infty}}$-homotopy equivalent) for any ${k}$ and ${k'}$. In particular, ${C_{K,i}=C_i}$ becomes a directed family of ${\epsilon}$-gapped ${X_K}$-modules. We denote this directed family of  ${\epsilon}$-gapped ${X_K}$-modules by ${\mathcal{C}}$.

Next we construct a directed family of ${\epsilon}$-gapped ${X_K}$-modules for ${\mathbb{Z}_p}$-equivariant Floer cohomology. Assume that ${P(H^{(p)},p\gamma)}$ is a finite set. Let ${\epsilon_1}$ and ${\epsilon_2}$ be positive constants as follows:
\begin{gather*}
\epsilon_1=\min \Big\{\int_{\mathbb{CP}^1}u^*\omega\neq 0 \ \Big| \ \begin{matrix}u:\mathbb{CP}^1\rightarrow M \\ J\circ du=du\circ j_{\mathbb{CP}} \end{matrix}\Big\}  \\
\epsilon_2=\min \Big\{ E(u)\neq 0 \ \Big| \ \begin{matrix}u:\mathbb{R}\times S^1\rightarrow M \\ \partial_su+J(\partial_tu-X_{H^{(p)}}(u))=0\end{matrix} \Big\}.
\end{gather*}
Note that ${\epsilon_1}$ and ${\epsilon_2}$ are positive (See Lemma 5 and Lemma 6). We fix a positive constant ${\epsilon<\min \{\epsilon_1,\epsilon_2\}}$. Let ${\{G_{(K,i)}\}}$ be a family of Hamiltonian functions in ${C^{\infty}(S^1\times M)}$ which satisfies the following conditions:
\begin{itemize}
\item Every $G_{(K,i)}$ is sufficiently close to $H^{(p)}$ so that the situation in Remark 48 below applies.
\item ${G_{(K,i)}\rightarrow H^{(p)}}$ ${(i\rightarrow \infty)}$ in ${C^{\infty}}$-topology 
\item ${(G_{(K,i)},J)}$ is Floer regular
\end{itemize}
Let ${\mathcal{G}^{(K,i)}_{w,t}}$ be a family of Hamiltonian functions parametrized by ${(w,t)\in S^{2K+1}\otimes S^1\subset S^{\infty}\otimes S^1}$ which satisfies the following conditions:
\begin{itemize}
\item (locally constant at critical points) For all $w$ in a small neighborhood of ${Z_i^m\in S^{2k+1}}$,
\begin{equation*}
\mathcal{G}_{w,t}^{(K,i)}(x)=G_{(K,i)}(t-\frac{m}{p},x)
\end{equation*}
holds.
\item ($\mathbb{Z}_p$-equivariance) ${\mathcal{G}^{(K,i)}_{m\cdot w,t}=\mathcal{G}^{(K,i)}_{w,t-\frac{m}{p}}}$
\item (invariance under the shift $\tau$) ${\mathcal{G}^{(K,i)}_{\tau(w),t}=\mathcal{G}^{(K,i)}_{w,t}}$ holds.
\end{itemize}
Assume that ${\widetilde{x}, \widetilde{y}\in P(H^{(p)},p\gamma)}$ are periodic orbits which split into ${\{x_1,\cdots,x_{k_1}\}\subset P(G_{(K,i)},p\gamma)}$ and ${\{y_1,\cdots,y_{k_2}\}\subset P(G_{(K,i)},p\gamma)}$. As in the definition of the ${\mathbb{Z}_p}$-equivariant Floer coboundary operator for toroidally monotone symplectic manifolds, we consider the following equation for ${x_{l_1},y_{l_2}\in P(G_{(K,i)},p\gamma)}$, ${m\in \mathbb{Z}_p}$, ${\lambda\ge 0}$, ${\alpha \in \{0,1\}}$ and ${0\le l\le 2K+1}$:

\begin{gather*}
(u,v)\in C^{\infty}(\mathbb{R}\times S^1,M)\times C^{\infty}(\mathbb{R},S^{2K+1}) \\
\partial _su(s,t)+J(u(s,t))(\partial_tu(s,t)-X_{\mathcal{G}^{(K,i)}_{v(s),t}})=0  \\
\frac{d}{ds}v(s)-\textrm{grad}(\widetilde{F})=0  \\
\lim_{s\to -\infty}v(s)=Z_{\alpha}^0, \lim_{s\to +\infty}v(s)=Z_l^m, \lim_{s\to -\infty}u(s,t)=x_{l_1}(t), \lim_{s\to +\infty}u(s,t)=y_{l_2}(t-\frac{m}{p}) \\
\int_{\mathbb{R}\times S^1}\widetilde{u}^*\omega+\int_0^1H^{(p)}(t,\widetilde{x}(t))-H^{(p)}(t,\widetilde{y}(t))dt =\lambda.
\end{gather*}
Here, ${\widetilde{u}}$ is a cylinder obtained as follows. Let ${v_{x_{l_1}}}$ and ${v_{y_{l_2}}}$ be small cylinders as follows:
\begin{gather*}
v_{x_{l_1}}:[0,1]\times S^1\rightarrow M, \ \ \ v_{y_{l_2}}:[0,1]\times S^1\rightarrow M \\
v_{x_{l_1}}(0,t)=\widetilde{x}(t) \ \ \ v_{x_{l_1}}(1,t)=x_{l_1}(t) \\
v_{y_{l_2}}(0,t)=y_{l_2}(t) \ \ \ v_{y_{l_2}}(1,t)=\widetilde{y}(t).
\end{gather*}
Moreover, we assume that ${v_{x_{l_1}}}$ is contained in a small neighborhood of ${\widetilde{x}}$ and ${v_{y_{l_2}}}$ is contained in a small neighborhood of ${\widetilde{y}}$. Then, we define ${\widetilde{u}}$ by ${\widetilde{u}=v_{x_{l_1}}\sharp u\sharp v_{y_{l_2}}}$.
\begin{Rem}
If ${G_{(K,i)}}$ is sufficiently close to ${H^{(p)}}$, ${\lambda \ge 0}$ holds. The proof is the same as the proof of Lemma 11. Moreover we assume that ${G_{(K,i)}}$ is sufficiently close to ${H^{(p)}}$ so that ${\lambda=0}$ or ${\lambda \ge \epsilon}$ holds. See Lemma 5 and Lemma 6.
\end{Rem}
We denote the space of solutions modulo the natural ${\mathbb{R}}$-action by ${\mathcal{N}_{\alpha,l,m}^{\lambda,(K,i)}(x_{l_1},y_{l_2})}$. For all solutions ${(u,v)\in C^{\infty}(\mathbb{R}\times S^1,M)\times C^{\infty}(\mathbb{R},S^{2K+1})}$ of the above equation, we have to achieve surjectivity of the linearization of the equation:
\begin{gather*}
\mathcal{F}_{(u,v)}:W^{k_1,p}(\mathbb{R}\times S^1,u^*TM)\times W^{k_2,2}(\mathbb{R},v^*TS^{2K+1})\longrightarrow \\ W^{k_1-1,p}(\mathbb{R}\times S^1,u^*TM)\times W^{k_2-1,2}(\mathbb{R},v^*TS^{2K+1}).
\end{gather*}

This ${\mathcal{F}_{(u,v)}}$ is a Fredholm operator. We can achieve surjectivity of ${\mathcal{F}_{(u,v)}}$ by perturbing a family of Hamiltonian functions ${\mathcal{G}_{w,t}^{(K,i)}}$. If $u$ is a non-trivial solution, this follows from the standard arguments in Floer homology (see \cite{FHS}). The non-trivial part is the perturbation for trivial $u$. This part is proved in \cite{SZ}. The determinant line bundle of ${\mathcal{F}_{(u,v)}}$ has the so-called ``coherent orientation" (\cite{SZ,FH}). The coherent orientation determines an orientation of moduli spaces ${\mathcal{N}_{\alpha,l,m}^{\lambda,(K,i)}(x_{l_1},y_{l_2})}$. In particular, we can assign ${+1}$ or ${-1}$ to each dimension $0$ component of the moduli space. What remains is the problem of compactness of the moduli space. The compactness is broken when sphere bubbles happens in the limit of the moduli space. However, we can avoid sphere bubbles in the limit of moduli spaces of dimension $0$ or $1$ if ${(M,\omega)}$ is weakly monotone. We can prove this fact by applying dimension counting  arguments and it is exactly the same as in case of Floer cohomology for weakly monotone $(M,\omega)$ in \cite{HS}, only the words are changed.

We define ${d_{\alpha,l}^{(K,i)}}$ as follows:
\begin{gather*}
d_{\alpha,l}^{(K,i)}:CF(G_{(K,i)},p\gamma:\Lambda_0)\longrightarrow CF(G_{(K,i)},p\gamma:\Lambda_0) \\
x\mapsto \sum_{m\in \mathbb{Z}_p}\sum_{\lambda \ge 0, y\in P(G_{(K,i)},p\gamma)}\sharp \mathcal{N}_{\alpha,l,m}^{\lambda,(K,i)}(x,y)\cdot T^{\lambda}y.
\end{gather*}
Here, ${\sharp \mathcal{N}_{\alpha,l,m}^{\lambda,(K,i)}(x,y)}$ is the number of dimension $0$ components of the moduli space ${\mathcal{N}_{\alpha,l,m}^{\lambda,(K,i)}(x,y)}$. Each component is counted as ${+1}$ or ${-1}$. Now we define an ${\epsilon}$-gapped ${X_{K}}$-module structure on 
\begin{equation*}
D_{K,i}=CF(G_{(K,i)},p\gamma:\Lambda_0)\otimes \Lambda_0\langle \theta \rangle .
\end{equation*}
We determine ${\{\delta_l^{(K,i)}:D_{K,i}\rightarrow D_{K,i}\}_{l=0}^K}$ as follows:
\begin{gather*}
\delta_l^{(K,i)}(x\otimes 1)=d_{0,2l}^{(K,i)}(x)\otimes 1+d_{0,2l+1}^{(K,i)}(x)\otimes \theta \\
\delta_l^{(K,i)}(x\otimes \theta)=d_{1,2l}^{(K,i)}(x)\otimes 1+d_{1,2l+1}^{(K,i)}(x)\otimes \theta.
\end{gather*}

\begin{Lem}
${(D_{K,i},\{\delta_l^{(K,i)}\}_{l=0}^K)}$ is an ${\epsilon}$-gapped ${X_K}$-module.
\end{Lem}
\vspace{5mm}
\textbf{proof}:
It suffices to prove that 
\begin{gather*}
(\delta^{(K,i)})^{(K)}\circ (\delta^{(K,i)})^{(K)}=0 \ \ \textrm{mod}(u^{K+1})
\end{gather*}
holds. Let ${\mathcal{N}_{\alpha,l,m}^{\lambda,(K,i)}(x,y)^{(\mu)}}$ be the sum of the dimension ${\mu}$ components of the moduli spaces. As in the case of the Floer cohomology, these moduli spaces have the natural compactifications and we have the following decomposition of the boundary:
\begin{gather*}
\partial \overline{\mathcal{N}_{\alpha,l,m}^{\lambda,(K,i)}(x,y)^{(1)}}=\bigsqcup_{\alpha''=[l'],l=l'+l''+[l'],m=m'+m'',\lambda=\lambda'+\lambda''} \mathcal{N}_{\alpha,l',m'}^{\lambda',(K,i)}(x,z)^{(0)}\times \mathcal{N}_{\alpha'',l'',m''}^{\lambda'',(K,i)}(z,y)^{(0)}.
\end{gather*}
This implies that the above equality holds.
\begin{flushright}
    $\Box$
\end{flushright}

Next, we define ${\epsilon}$-gapped ${X_K}$-morphisms ${\iota_{(K,i\rightarrow j)}:D_{K,i}\rightarrow D_{K,j}}$ for ${i,j\in \mathbb{N}}$. We consider a family of Hamiltonian functions which connects ${\mathcal{G}_{w,t}^{(K,i)}}$ and ${\mathcal{G}_{w,t}^{(K,j)}}$. Let ${\mathcal{G}_{s,w,t}^{(K,i\rightarrow j)}}$ be a family of Hamiltonian functions parametrized by ${(s,w,t)\in \mathbb{R}\times S^{2K+1}\times S^1}$ which satisfies the following conditions:

\begin{itemize}
\item ${\mathcal{G}_{s,w,t}^{(K,i\rightarrow j)}(x)=\begin{cases} \mathcal{G}_{w,t}^{K,i}(x) & s\ll 0 \\ \mathcal{G}_{w,t}^{K,j}(x) & s\gg 0 \end{cases}}$
\item ($\mathbb{Z}_p$-equivariance) ${\mathcal{G}_{s,mw,t}^{(K,i\rightarrow j)}(x)=\mathcal{G}_{s,w,t-\frac{m}{p}}^{(K,i\rightarrow j)}(x) \ \ \ (\forall m\in \mathbb{Z}_p)}$
\item (invariance under the shift ${\tau}$) ${\mathcal{G}_{s,\tau(w),t}^{(K,i\rightarrow j)}(x)=\mathcal{G}_{s,w,t}^{(K,i\rightarrow j)}(x)}$
\end{itemize}

Assume that ${\widetilde{x}, \widetilde{y}\in P(H^{(p)},p\gamma)}$ are periodic orbits which split into ${\{x_1,\cdots,x_{k_1}\}\subset P(G_{(K,i)},p\gamma)}$ and ${\{y_1,\cdots,y_{k_2}\}\subset P(G_{(K,j)},p\gamma)}$. We consider the following equation for ${x_{l_1}\in P(G_{(K,i)},p\gamma)}$, ${y_{l_2}\in P(G_{(K,j)},p\gamma)}$, ${m\in \mathbb{Z}_p}$, ${\lambda \ge 0}$, ${\alpha\in \{0,1\}}$ and ${0\le l \le 2K+1}$:

\begin{gather*}
(u,v)\in C^{\infty}(\mathbb{R}\times S^1,M)\times C^{\infty}(\mathbb{R},S^{2K+1}) \\
\partial _su(s,t)+J(u(s,t))(\partial_tu(s,t)-X_{\mathcal{G}^{(K,i\rightarrow j)}_{s,v(s),t}})=0  \\
\frac{d}{ds}v(s)-\textrm{grad}(\widetilde{F})=0  \\
\lim_{s\to -\infty}v(s)=Z_{\alpha}^0, \lim_{s\to +\infty}v(s)=Z_l^m, \lim_{s\to -\infty}u(s,t)=x_{l_1}(t), \lim_{s\to +\infty}u(s,t)=y_{l_2}(t-\frac{m}{p}) \\
\int_{\mathbb{R}\times S^1}\widetilde{u}^*\omega+\int_0^1H^{(p)}(t,\widetilde{x}(t))-H^{(p)}(t,\widetilde{y}(t))dt =\lambda.
\end{gather*}

Here, ${\widetilde{u}}$ is a cylinder obtained as before (So $\widetilde{u}$ connects ${\widetilde{x}}$ and ${\widetilde{y}}$.). We denote the space of solutions by ${\mathcal{N}_{\alpha,l,m}^{\lambda,(K,i\rightarrow j)}(x_{l_1},y_{l_2})}$. We define ${\iota_{\alpha,l,(K,i\rightarrow j)}}$ as follows:

\begin{gather*}
\iota_{\alpha,l,(K,i\rightarrow j)}(x)=\sum_{m\in \mathbb{Z}_p,\lambda \ge 0,y\in P(G_{(K,j)},p\gamma)}\sharp \mathcal{N}_{\alpha,l,m}^{\lambda,(K,i\rightarrow j)}(x,y)\cdot T^{\lambda}y.
\end{gather*}
Here ${\sharp \mathcal{N}_{\alpha,l,m}^{\lambda,(K,i\rightarrow j)}(x,y)}$ is the number of dimension $0$ components of the moduli space ${\mathcal{N}_{\alpha,l,m}^{\lambda,(K,i\rightarrow j)}(x,y)}$. Now we define an ${\epsilon}$-gapped ${X_K}$-morphism ${\{\iota_{(K,i\rightarrow j),l}\}_{l=0}^K:D_{K,i}\rightarrow D_{K,j}}$ as follows:

\begin{gather*}
\iota_{(K,i\rightarrow j),l}(x\otimes 1)=\iota_{0,2l,(K,i\rightarrow j)}(x)\otimes 1+\iota_{0,2l+1,(K,i\rightarrow j)}(x)\otimes \theta  \\
\iota_{(K,i\rightarrow j),l}(x\otimes \theta)=\iota_{1,2l,(K,i\rightarrow j)}(x)\otimes 1+\iota_{1,2l+1,(K,i\rightarrow j)}(x)\otimes \theta.
\end{gather*}
We assume that ${\mathcal{G}_{s,w,t}^{(K,i\rightarrow j)}(x)}$ is sufficiently close to ${H^{(p)}(t,x)}$ so that ${\iota_{(K,i\rightarrow j)}}$ is an ${\epsilon}$-gapped ${X_K}$-morphism over ${\Lambda_0}$ (see Lemma 5 and Lemma 6).

\begin{Lem}
${\iota_{(K,i\rightarrow j)}}$ is an ${\epsilon}$-gapped $X_K$-morphism and ${\iota_{(K,j\rightarrow k)}\circ \iota_{(K,i\rightarrow j)}}$ is ${\epsilon}$-gapped ${X_K}$-homotopic to ${\iota_{(K,i\rightarrow k)}}$ for any ${i,j,k\in \mathbb{N}}$. In particular, ${\iota_{(K,j\rightarrow i)}\circ \iota_{(K,i\rightarrow j)}}$ is ${\epsilon}$-gapped ${X_K}$-homotopic to the identity and ${ \iota_{(K,i\rightarrow j)}}$ is an ${\epsilon}$-gapped ${X_K}$-homotopy equivalence between ${D_{K,i}}$ and ${D_{K,j}}$ for any ${i,j\in \mathbb{N}}$.
\end{Lem}
\vspace{5mm}
\textbf{proof}:
Let ${\mathcal{N}_{\alpha,l,m}^{\lambda,(K,i\rightarrow j)}(x,y)^{(\mu)}}$ be the sum of dimension ${\mu}$ components of ${\mathcal{N}_{\alpha,l,m}^{\lambda,(K,i\rightarrow j)}(x,y)}$. Then, the boundary of ${\mathcal{N}_{\alpha,l,m}^{\lambda,(K,i\rightarrow j)}(x,y)^{(1)}}$ is written as follows:

\begin{gather*}
\partial \overline{\mathcal{N}_{\alpha,l,m}^{\lambda,(K,i\rightarrow j)}(x,y)^{(1)}}
=\bigsqcup_{\alpha''=[l'],l=l'+l''+[l'],m=m'+m'',\lambda=\lambda'+\lambda''}\mathcal{N}_{\alpha,l',m'}^{\lambda',(K,i)}(x,z)^{(0)}\times \mathcal{N}_{\alpha'',l'',m''}^{\lambda'',(K,i\rightarrow j)}(z,y)^{(0)} \\
\sqcup \bigsqcup_{\alpha''=[l'],l=l'+l''+[l'],m=m'+m'',\lambda=\lambda'+\lambda''}\mathcal{N}_{\alpha,l',m'}^{\lambda',(K,i\rightarrow j)}(x,z)^{(0)}\times \mathcal{N}_{\alpha'',l'',m''}^{\lambda'',(K,j)}(z,y)^{(0)}.
\end{gather*}
This implies that 
\begin{gather*}
(\iota_{(K,i\rightarrow j)})^{(K)}\circ (\delta^{(K,i)})^{(K)}\equiv (\delta^{(K,j)})^{(K)}\circ (\iota_{(K,i\rightarrow j)})^{(K)} \ \ \ \textrm{mod}(u^{K+1})
\end{gather*}
holds and ${\iota_{(K,i\rightarrow j)}}$ is an ${\epsilon}$-gapped ${X_K}$-morphism. 

Next, we fix ${R>0}$ so that 
\begin{gather*}
\mathcal{G}_{s,w,t}^{(K,i\rightarrow j)}(x)=
\begin{cases}
\mathcal{G}_{w,t}^{(K,i)}(x) & s\le -R  \\
\mathcal{G}_{w,t}^{(K,i)}(x) & s\ge R
\end{cases}
\end{gather*}
and 
\begin{gather*}
\mathcal{G}_{s,w,t}^{(K,j\rightarrow k)}(x)=
\begin{cases}
\mathcal{G}_{w,t}^{(K,j)}(x) & s\le -R  \\
\mathcal{G}_{w,t}^{(K,k)}(x) & s\ge R
\end{cases}
\end{gather*}
hold. We define a family of Hamiltonian functions parametrized by 
\begin{gather*}
(s,\rho,w,t)\in \mathbb{R}\times \mathbb{R}_{\ge 0}\times S^{2K+1}\times S^1
\end{gather*}
as follows:
\begin{gather*}
\mathcal{G}_{(s,\rho,w,t)}(x)=
\begin{cases}  
\mathcal{G}_{w,t}^{(K,i)}(x) & s\ll 0  \\
\mathcal{G}_{w,t}^{(K,k)}(x) & s\gg 0
\end{cases}
\end{gather*}
\begin{gather*}
\mathcal{G}_{(s,0,w,t)}(x)=\mathcal{G}_{s,w,t}^{(K,i\rightarrow k)}(x)
\end{gather*}
\begin{gather*}
\mathcal{G}_{(s,\rho,w,t)}(x)=
\begin{cases}
\mathcal{G}_{s+\rho,w,t}^{(K,i\rightarrow j)}(x) & s\le 0, \ \rho \ge 2R  \\
\mathcal{G}_{s-\rho,w,t}^{(K,j\rightarrow k)}(x) & s\ge 0, \ \rho \ge 2R
\end{cases}
\end{gather*}
\begin{itemize}
\item ($\mathbb{Z}_p$-equivariance) ${\mathcal{G}_{(s,\rho,w,t)}(x)=\mathcal{G}_{(s,\rho,w,t-\frac{m}{p})}(x)}$ 
\item (invariance under the shift ${\tau}$) ${\mathcal{G}_{(s,\rho,\tau(w),t)}(x)=\mathcal{G}_{(s,\rho,w,t)}(x)}$
\end{itemize}
For ${x\in P(G_{(K,i)},p\gamma)}$, ${y\in P(G_{K,j},p\gamma)}$, ${m\in \mathbb{Z}_p}$, ${\lambda \ge 0}$, ${\alpha \in \{0,1\}}$ and ${0\le l\le 2K+1}$, we consider the following equation:

\begin{gather*}
(u,v,\rho)\in C^{\infty}(\mathbb{R}\times S^1,M)\times C^{\infty}(\mathbb{R},S^{2K+1})\times \mathbb{R}_{\ge 0} \\
\partial_su(s,t)+J(u(s,t))(\partial_tu(s,t)-X_{\mathcal{G}_{(s,\rho,v(s),t)}})=0  \\
\frac{d}{ds}f(s)-\textrm{grad}(\widetilde{F})=0  \\
\lim_{s\rightarrow -\infty}v(s)=Z_{\alpha}^0, \lim_{s\rightarrow +\infty}v(s)=Z_l^m, \ \lim_{s\rightarrow -\infty}u(s,t)=x(t), \ \lim_{s\rightarrow +\infty}u(s,t)=y(t-\frac{m}{p})  \\
\int_{\mathbb{R}\times S^1}\widetilde{u}^*\omega+\int_0^1H^{(p)}(t,\widetilde{x})-H^{(p)}(t,\widetilde{y})dt=\lambda.
\end{gather*}

Here, ${\widetilde{x},\widetilde{y}\in P(H^{(p)},p\gamma)}$ are periodic orbits which generate $x$ and ${y}$ respectively, and ${\widetilde{u}}$ is a cylinder obtained from $u$ as before. We denote the space of solutions by ${\mathcal{M}_{\alpha,l,m}^{\lambda}(x,y)}$. We define ${h_{\alpha,l}}$ as follows:

\begin{gather*}
h_{\alpha,l}(x)=\sum_{m\in \mathbb{Z}_p,\lambda \ge 0, y\in P(D_{K,k},p\gamma)}\sharp \mathcal{M}_{\alpha,l,m}^{\lambda}(x,y)\cdot T^{\lambda}y.
\end{gather*}
We define a family of maps ${\{h_l\}_{l=0}^K}$ as follows:
\begin{gather*}
h_l(x\otimes 1)=h_{0,2l}(x)\otimes 1+h_{0,2l+1}(x)\otimes \theta \\
h_l(x\otimes \theta)=h_{1,2l}(x)\otimes 1+h_{1,2l+1}(x)\otimes \theta.
\end{gather*}
Let ${\mathcal{N}_{\alpha,l,m}^{\lambda,(K,i)}(x,y)^{(\mu)}}$, ${\mathcal{N}_{\alpha,l,m}^{\lambda,(K,i\rightarrow j)}(x,y)^{(\mu)}}$ and ${\mathcal{M}_{\alpha,l,m}^{\lambda}(x,y)^{(\mu)}}$ be the sum of the dimension ${\mu}$ components of the moduli spaces. As in the case of the Floer cohomology, these moduli spaces have the natural compactifications. In particular, we have the following decomposition:
\begin{gather*}
\partial \overline{\mathcal{M}_{\alpha,l,m}^{\lambda}(x,y)^{(1)}}=\bigsqcup_{\alpha''=[l'],l=l'+l'+[l'],m=m'+m'',\lambda=\lambda'+\lambda''} \mathcal{N}_{\alpha,l',m'}^{\lambda',(K,i\rightarrow j)}(x,z)^{(0)}\times \mathcal{N}_{\alpha'',l'',m''}^{\lambda'',(K,j\rightarrow k)}(z,y)^{(0)} \\
\sqcup \bigsqcup_{\alpha''=[l'],l=l'+l''+[l'],m=m'+m'',\lambda=\lambda'+\lambda''} \mathcal{M}_{\alpha,l',m'}^{\lambda}(x,z)^{(0)}\times \mathcal{N}_{\alpha'',l'',m''}^{\lambda'',(K,k)}(z,y)^{(0)} \\
\sqcup \bigsqcup_{\alpha''=[l'],l=l'+l''+[l'],m=m'+m'',\lambda=\lambda'+\lambda''} \mathcal{N}_{\alpha,l',m'}^{\lambda',(K,i)}(x,z)^{(0)}\times \mathcal{M}_{\alpha'',l'',m''}^{\lambda''}(z,y)^{(0)} \\
\sqcup \ \mathcal{N}_{\alpha,l,m}^{\lambda,(K,i\rightarrow k)}(x,y)^{(0)}.
\end{gather*}
The identity is to be understood with orientations (see \cite{SZ}). This implies that 
\begin{gather*}
\iota_{(K,j\rightarrow k)}^{(K)}\circ \iota_{(K,i\rightarrow j)}^{(K)}+(\delta^{(K,k)})^{(K)}\circ h^{(K)}+h^{(K)}\circ (\delta^{(K,i)})^{(K)}-\iota_{(K,i\rightarrow k)}^{(K)}=0 \ \textrm{mod}(u^{K+1})
\end{gather*}
holds. In particular, ${\iota_{(K,j\rightarrow k)}\circ \iota_{(K,i\rightarrow j)}}$ is ${\epsilon}$-gapped ${X_K}$-homotopic to ${\iota_{(K,i\rightarrow k)}}$.
\begin{flushright}
    $\Box$
\end{flushright}

Next, we define ${\epsilon}$-gapped ${X_K}$-morphisms 
\begin{gather*}
\tau_{(K\rightarrow K+1,i)}:D_{K,i}\longrightarrow D_{K+1,i}
\end{gather*}
in the same way we defined ${\iota_{(K,i\rightarrow j)}}$. If $i$ is sufficiently large, parametrized Hamiltonian functions ${\mathcal{G}_{w,t}^{(K,i)}}$ and ${\mathcal{G}_{w,t}^{(K+1,i)}}$ are sufficiently close to ${H^{(p)}}$. Then, there is a constant ${N(K)\in \mathbb{N}}$ such that ${\epsilon}$-gapped ${X_K}$-morphisms which are constructed in the same way as ${\iota_{(K,i\rightarrow j)}}$
\begin{gather*}
\tau_{(K\rightarrow K+1,i)}:D_{K,i}\longrightarrow D_{K+1,i}
\end{gather*}
are well-defined over $\Lambda_0$ for ${i\ge N(K)}$ for the same reason I explained in Remark 48.

\begin{Lem}
${\tau_{(K\rightarrow K+1,j)}\circ \iota_{(K,i\rightarrow j)}}$ is ${\epsilon}$-gapped ${X_K}$-homotopic to ${\iota_{(K+1,i\rightarrow j)}\circ \tau_{(K\rightarrow K+1,i)}}$ for any ${i,j\ge N(K)}$. Here, ${\tau_{(K+1,i\rightarrow j)}}$ is the restriction of ${\tau_{(K+1,i\rightarrow j)}}$ to an $\epsilon$-gapped $X_K$-morphism.
\end{Lem}

The proof is essentially the same as the proof of Lemma 50. So, ${\{D_{K,i},\iota_{(K,i\rightarrow j)},\tau_{(K\rightarrow K+1,j)}\}}$ is a directed family of ${\epsilon}$-gapped $X_K$ modules. We denote it by ${\mathcal{D}}$.

Next, we define a morphism ${\{P_{(K,i)}:C_{K,i}\rightarrow D_{K,i}\}}$ between two directed families of $\epsilon$-gapped ${X_K}$-modules $\mathcal{C}$ and ${\mathcal{D}}$. By proposition 44, this morphism determines an ${\epsilon}$-gapped ${X_{\infty}}$-morphism and 
this is the ${\mathbb{Z}_p}$-equivariant pants product. First, we define a ${p}$-legged pants ${\Sigma_p}$ as follows:

\begin{gather*}
\Sigma_p=\Big(\bigsqcup_{0\le k\le p-1}\mathbb{R}\times [k,k+1]\Big)\Big/\sim.
\end{gather*}
The equivalence relation ${\sim}$ is defined as follows:
\begin{itemize}
\item ${(s,k)\in [0,\infty)\times [k-1,k]}$ is equivalent to ${(s,k)\in [0,\infty)\times [k,k+1]}$ for ${1\le k\le p-1}$.
\item ${(s,0)\in [0,\infty)\times [0,1]}$ is equivalent to ${(s,p)\in [0,\infty)\times [p-1,p]}$.
\item ${(s,k)\in (-\infty,0]\times [k,k+1]}$ is equivalent to ${(s,k+1)\in (-\infty,0]\times [k,k+1]}$ for ${0\le k\le p-1}$.
\end{itemize}
We determine a complex structure near ${[(0,0)]\in \Sigma_p}$. We give an explicit local coordinate of a neighborhood of  ${[(0,0)]\in \Sigma_p}$ as follows:
\begin{gather*}
w:\{z\in \mathbb{C} \ | \ |z|\le \frac{1}{2}\}\longrightarrow \Sigma_p
\end{gather*}
\begin{gather*}
w(z)=
\begin{cases}
[z^p] & 0\le \textrm{arg}(z)\le \frac{\pi}{p} \\
[z^p+k\sqrt{-1}] & \frac{(2k-1)\pi}{p}\le \textrm{arg}(z)\le \frac{(2k+1)\pi}{p} \\
[z^p+p\sqrt{-1}] & \frac{(2p-1)\pi}{p}\le \textrm{arg}(z)\le 2\pi.
\end{cases}
\end{gather*}
Here we identify ${\mathbb{R}^2}$ with ${\mathbb{C}}$. Without loss of generality, we assume that ${H(t,x)=0}$ near ${t=0}$. Let ${\mathcal{K}_{w,z}^{(K,i)}}$ be a family of Hamiltonian functions parametrized by ${(w,z)\in S^{2K+1}\times \Sigma_p}$ as follows:
\begin{itemize}
\item ${\mathcal{K}_{w,z}^{(K,i)}(x)=H_{K}([t],x)}$ \ \  ${(z=[(s,t)]\in \Sigma_p, s\ll 0)}$
\item ${\mathcal{K}_{w,z}^{(K,i)}(x)=\frac{1}{p}G_{(K,i)}(\frac{t}{p},x)}$ \ \ ${(z=[(s,t)]\in \Sigma_p, s\gg 0)}$
\item ${\mathcal{K}_{w,z}^{(K,i)}(x)\stackrel{i\rightarrow \infty}{\longrightarrow }H([t],x)}$  \ \ ${z=[(s,t)]}$
\item (${\mathbb{Z}_p}$-invariance) ${\mathcal{K}_{mw,z}^{(K,i)}(x)=\mathcal{K}_{w,z+m\sqrt{-1}}^{(K,i)}(x)}$
\item (invariance under the shift ${\tau}$) ${\mathcal{K}_{\tau(w),z}^{(K,i)}(x)=\mathcal{K}_{w,z}^{(K,i)}(x)}$
\end{itemize}
For ${\{x_k\}_{k=1}^p\subset P(H_K,\gamma)}$, ${y\in P(D_{K,i},p\gamma)}$, ${m\in \mathbb{Z}_p}$, ${\lambda \ge 0}$, ${\alpha\in \{0,1\}}$, and ${0\le l\le 2K+1}$, we consider the following equation:
\begin{gather*}
(u,v)\in C^{\infty}(\Sigma_p,M)\times C^{\infty}(\mathbb{R},S^{2K+1})  \\
\partial_su([(s,t)])+J(u([(s,t)]))(\partial_tu([(s,t)])-X_{\mathcal{K}_{v(s),[(s,t)]}^{(K,i)}})=0  \\
\frac{d}{ds}v(s)-\textrm{grad}(\widetilde{F})=0  \\
\lim_{s\rightarrow -\infty}v(s)=Z_{\alpha}^0, \lim_{s\rightarrow +\infty}v(s)=Z_l^m, \lim_{s\rightarrow -\infty}u([(s,t)])=x_i([t]) \ \ (t\in [i-1,i])  \\
\lim_{s\rightarrow +\infty}u([(s,t)])=y\big(\frac{t}{p}-m\big)  \\
\int_{\Sigma_p}\widetilde{u}^*\omega+\int_0^1H^{(p)}(t,\widetilde{y}(t))dt-\sum_{i=1}^p\int_0^1H(t,\widetilde{x}_i(t))dt=\lambda.
\end{gather*}
Here ${\{\widetilde{x}_i\}\subset P(H,\gamma)}$ and ${\widetilde{y}\in P(H^{(p)},p\gamma)}$ are periodic orbits as before and ${\widetilde{u}}$ is a map which connects ${\{\widetilde{x}_i\}}$ and ${\widetilde{y}}$. We denote the space of solutions by ${\mathcal{M}_{\alpha,l,m}^{\lambda,(K,i)}(x_1,\cdots,x_p:y)}$. We define ${f_{\alpha,l}^{(K,i)}}$ as follows:
\begin{gather*}
f_{\alpha,l}^{(K,i)}:C_{K,i}\longrightarrow D_{K,i}  \\
x_1\otimes \cdots \otimes x_p \mapsto \sum_{m\in \mathbb{Z}_p,\lambda \ge 0, y\in P(G_{(K,i)},p\gamma)} \sharp \mathcal{M}_{\alpha,l,m}^{\lambda,(K,i)}(x_1,\cdots,x_p:y)T^{\lambda}y.
\end{gather*}
Then, we define ${\mathfrak{f}^{(K,i)}=\{f_l^{(K,i)}\}_{l=0}^K}$ as follows:
\begin{gather*}
f_l^{(K,i)}((x_1\otimes \cdots \otimes x_p)\otimes 1)=f_{0,2l}^{(K,i)}(x_1\otimes \cdots \otimes x_p)\otimes 1+f_{0,2l+1}^{(K,i)}(x_1\otimes \cdots \otimes x_p)\otimes \theta  \\
f_l^{(K,i)}((x_1\otimes \cdots \otimes x_p)\otimes \theta)=f_{1,2l}^{(K,i)}(x_1\otimes \cdots \otimes x_p)\otimes 1+f_{1,2l+1}^{(K,i)}(x_1\otimes \cdots \otimes x_p)\otimes \theta.
\end{gather*}

\begin{Lem}
The system ${\mathcal{F}=\{\mathfrak{f}^{(K,i)}\}}$ determines a morphism between a directed family of ${\epsilon}$-gapped ${X_K}$-modules ${\{C_{K,i}\}}$ and a  directed family of ${\epsilon}$-gapped ${X_K}$-modules ${\{D_{K,i}\}}$.
\end{Lem}
\vspace{5mm}
\textbf{proof}:
If ${\mathcal{K}_{w,z}^{(K,i)}(x)}$ is sufficiently close to ${H(t,x)}$, ${\mathfrak{f}^{(K,i)}}$ becomes a morphism between ${\epsilon}$-gapped ${X_K}$-modules. This follows from the same argument as in the proof of Lemma 11. It suffices to prove that ${\mathfrak{f}^{(K,j)}\circ \iota_{(K,i\rightarrow j)}}$ and
${\iota_{(K,i\rightarrow j)}\circ \mathfrak{f}^{(K,i)}}$ are ${\epsilon}$-gapped ${X_K}$-homotopic and ${\mathfrak{f}^{(K+1,i)}\circ \tau_{(K\rightarrow K+1,i)}}$ and ${\tau_{(K\rightarrow K+1,i)}\circ \mathfrak{f}^{(K,i)}}$ are ${\epsilon}$-gapped ${X_K}$-morphism for sufficiently large $i\in \mathbb{N}$. However, these facts follow from the same arguments as in the proof of Lemma 50. 
\begin{flushright}
    $\Box$
\end{flushright}

We fix ${(K,i)}$ and denote ${H_K}$ by ${\overline{H}}$ and ${G_{(K,i)}}$ by ${\widetilde{H^{(p)}}}$. Proposition 44 implies that we have the following two ${\epsilon}$-gapped ${X_{\infty}}$-modules
\begin{gather*}
(C,\mathfrak{d})=((CF(\overline{H},\gamma)^{\otimes p}\otimes \Lambda_0\langle \theta \rangle, \{d_i\}) \\
(D,\mathfrak{l})=(CF(\widetilde{H^{(p)}},p\gamma)\otimes \Lambda_0\langle \theta \rangle, \{l_i\})
\end{gather*}
and homotopy equivalences between directed families of $\epsilon$-gapped $X_K$-modules 
\begin{gather*}
    \mathcal{I}:(C,\mathfrak{d})\longrightarrow \mathcal{C} \\
    \mathcal{J}:(D,\mathfrak{l})\longrightarrow \mathcal{D}
\end{gather*}
and an ${\epsilon}$-gapped ${X_{\infty}}$-morphism
\begin{gather*}
\mathfrak{g}=\{g_i\}:(C,\mathfrak{d})\longrightarrow (D,\mathfrak{l})
\end{gather*}
such that the following diagram (of directed families of $\epsilon$-gapped $X_K$-modules) is commutative up to homotopy.

\[
\begin{tikzcd}
C \arrow[rrrr, "\mathcal{I}"] \arrow[ddd, "\mathfrak{g}=\{g_i\}"'] &  &  &  & \mathcal{C} \arrow[ddd, "\mathcal{F}=\{\mathfrak{f}^{(K,i)}\}"] \\
                                      &  &  &  &                     \\
                                      &  &  &  &                     \\ \ D \arrow[rrrr, "\mathcal{J}"]                   &  &  &  & \mathcal{D}               
\end{tikzcd}
\]
Moreover, we can define morphisms between cohomologies of directed familes as follows:
\begin{gather*}
    \mathcal{F}_*^{\mathfrak{g}}:H(\mathcal{C})\longrightarrow H(\mathcal{D}) \\
    \widehat{\mathcal{F}}_*^{\mathfrak{g}}:\widehat{H}(\mathcal{C})\longrightarrow \widehat{H}(\mathcal{D})  \\
    \mathcal{F}_{loc,*}^{\mathfrak{g}}:H_{loc}(\mathcal{C})\longrightarrow H_{loc}(\mathcal{D}) \\
    \widehat{\mathcal{F}}_{loc,*}^{\mathfrak{g}}:\widehat{H}_{loc}(\mathcal{C})\longrightarrow \widehat{H}_{loc}(\mathcal{D}).
\end{gather*}

 The cochain complex ${(C\otimes \Lambda_0[[u]],\mathfrak{d}^{(\infty)})}$ is the ${\mathbb{Z}_p}$-equivariant cochain complex for ${CF(H,\gamma)^{\otimes p}}$ and the cochain complex ${(D\otimes \Lambda_0[[u]],\mathfrak{l}^{(\infty)})}$ is the ${\mathbb{Z}_p}$-equivariant Floer cochain complex for $H^{(p)}$. 

\begin{Rem}[Independence of the choices of directed families]
We defined a directed family of $\epsilon$-gapped $X_K$-modules $\mathcal{C}$ by taking a convergent sequence of Hamiltonian functions ${H_i\rightarrow H}$. If we take another convergent sequence ${K_i\rightarrow H}$, we get another directed family of $\epsilon$-gapped $X_K$-modules $\mathcal{C}'$. However, the usual Floer continuation map
\begin{gather*}
    CH(H_i,\gamma:\Lambda_0)\longrightarrow CH(K_i,\gamma:\Lambda_0)
\end{gather*}
induces a homotopy equivalence ${\mathcal{C}\rightarrow\mathcal{C}'}$. In particular, their cohomologies ${H}$, ${\widehat{H}}$, ${H_{loc}}$ and ${\widehat{H}_{loc}}$ are isomorphic. So, we can define equivariant cohomologies independent of the choice of directed families of $\epsilon$-gapped $X_K$-modules $\mathcal{C}$ (The same holds for $\mathcal{D}$.).
\end{Rem}

Now we define
\begin{gather*}
H(\mathbb{Z}_p,CF(H,\gamma:\Lambda_0)^{\otimes p})\stackrel{\mathrm{def}}{=}H(\mathcal{C}) \ \ \big(=H(C\otimes \Lambda_0[[u]],\mathfrak{d}^{(\infty)}) \ \big)  
\end{gather*}
\begin{gather*}
\widehat{H}(\mathbb{Z}_p,CF(H,\gamma:\Lambda_0)^{\otimes p})\stackrel{\mathrm{def}}{=}\widehat{H}(\mathcal{C}) \ \ \big(=H(C\otimes \Lambda_0[u^{-1},u]],\mathfrak{d}^{(\infty)}) \ \big)
\end{gather*}
\begin{gather*}
HF_{\mathbb{Z}_p}(H^{(p)},p\gamma)\stackrel{\mathrm{def}}{=}H(\mathcal{D}) \ \ \big(=H(D\otimes \Lambda_0[[u]],\mathfrak{l}^{(\infty)}) \ \big)
\end{gather*}
\begin{gather*}
\widehat{HF}_{\mathbb{Z}_p}(H^{(p)},p\gamma)\stackrel{\mathrm{def}}{=}\widehat{H}(\mathcal{D}) \ \ \big(=H(D\otimes \Lambda_0[u^{-1},u]],\mathfrak{l}^{(\infty)}) \ \big).
\end{gather*}
Remark 53 implies that these definitions are independent of the choices of directed families of ${\epsilon}$-gapped ${X_K}$-modules. The morphisms ${\mathcal{F}_*^{\mathfrak{g}}}$ and ${\widehat{\mathcal{F}}_*^{\mathfrak{g}}}$ determine the ${\mathbb{Z}_p}$-equivariant pair of pants products:

\begin{gather*}
\mathcal{P}=\mathcal{F}_*^{\mathfrak{g}}:H(\mathbb{Z}_p,CF(H,\gamma:\Lambda_0)^{\otimes p})\longrightarrow  HF_{\mathbb{Z}_p}(H^{(p)},p\gamma)  \\
\widehat{\mathcal{P}}=\widehat{\mathcal{F}}_*^{\mathfrak{g}}:\widehat{H}(\mathbb{Z}_p,CF(H,\gamma:\Lambda_0)^{\otimes p})\longrightarrow \widehat{HF}_{\mathbb{Z}_p}(H^{(p)},p\gamma).
\end{gather*}
Moreover, the underlying local morphisms
\begin{gather*}
\mathcal{P}_{loc}=\mathcal{F}_{loc,*}^{\mathfrak{g}}:\bigoplus_{x\in P(H,\gamma)}H(\mathbb{Z}_p,CF^{loc}(H,x)^{\otimes p})\longrightarrow \bigoplus_{x\in P(H,\gamma)}HF_{\mathbb{Z}_p}^{loc}(H^{(p)},x^{(p)}) \\
\widehat{\mathcal{P}}_{loc}=\widehat{\mathcal{F}}_{loc,*}^{\mathfrak{g}}:\bigoplus_{x\in P(H,\gamma)}\widehat{H}(\mathbb{Z}_p,CF^{loc}(H,x)^{\otimes p})\longrightarrow \bigoplus_{x\in P(H,\gamma)}\widehat{HF}_{\mathbb{Z}_p}^{loc}(H^{(p)},x^{(p)})
\end{gather*}
are equal to the local ${\mathbb{Z}_p}$-equivariant pair of pants product defined in \cite{SZ}. In particular, ${\widehat{\mathcal{P}}_{loc}}$ is an isomorphism. Using the same arguments as in the toroidally monotone case, we can apply this local isomorphism to prove that ${\widehat{\mathcal{P}}}$ is an isomorphism. The only non-trivial point is the proof of Lemma 16 in the weakly monotone case. In the proof of Lemma 16, we used a cochain homotopy ${\mathcal{L}}$ to prove that for any cocycle ${x\in F^{q+r}A}$ we can find ${y\in F^{q+1}A}$ such that ${d(y)=x}$ holds. In order to construct ${\mathcal{L}}$ for ${A=(D\otimes \Lambda_0[u^{-1},u]],\mathfrak{l}^{(\infty)})}$, we need the uniqueness of ${X_{\infty}}$-morphism up to ${X_{\infty}}$-hotomopy. However, we only proved the uniqueness up to ${X_K}$-homotopy for any ${K<\infty}$. So we have to give another proof in the weakly monotone case. 
\vspace{5mm}
\textbf{proof}(Lemma 16):
It suffices to prove that for any cocycle ${x\in (D\otimes \Lambda_0[u^{-1},u]],\mathfrak{l}^{(\infty)})}$ we can find ${y\in T^{-p||H||}(D\otimes \Lambda_0[u^{-1},u]])}$ so that ${\mathfrak{l}^{(\infty)}(y)=x}$ holds. Note that ${(D\otimes \Lambda,l_0)}$ is acyclic and we can construct a cochain homotopy ${\mathcal{H}}$
\begin{gather*}
\mathcal{H}:D\otimes \Lambda\longrightarrow D\otimes \Lambda  \\
\mathcal{H}(D)\subset T^{-p||H||}D \\
\textrm{Id}_{D\otimes \Lambda}=l_0\mathcal{H}+\mathcal{H}l_0.
\end{gather*}
Assume that $x$ is an infinite sum ${x=\sum_{k\ge 0}x_ku^k}$ (${x_k\in D}$). Our purpose is to construct ${y=\sum_{k\ge 0}y_ku^k}$ (${y_k\in T^{-p||H||}D}$) so that ${\mathfrak{l}^{(\infty)}(y)=x}$ holds. We define ${y_0=\mathcal{H}(x_0)}$. Then ${\mathfrak{l}^{(\infty)}(y_0)=x_0+\sum_{k\ge 1}x'_ku^k}$ holds. Moreover, ${x-\mathfrak{l}^{(\infty)}(y_0)=\sum_{k\ge 1}(x_k-x'_k)u^k}$ is a cocyle. Inductively, we can determine ${y_1,y_2,\cdots}$ so that ${\mathfrak{l}^{(\infty)}(y)=x}$ and ${y\in T^{-p||H||}D}$ hold. 
\begin{flushright}
    $\Box$
\end{flushright}

The rest of the proof of Theorem 2 for the weakly monotone case is the same as in the toroidally monotone case. 

\begin{Rem}
We can prove that ${\epsilon}$-gapped ${X_{\infty}}$-homotopy class of ${\mathfrak{g}}$ does not depend on the algebraic construction and approximations by applying ``homotopy of homotopy" theory. In particular, we can also prove that ${\mathcal{P}}$ and ${\widehat{\mathcal{P}}}$ does not depend on the algebraic construction and approximations. However, the uniqueness of ${\widehat{\mathcal{P}}_{loc}}$ is sufficient for our proof of Theorem 2. So, Proposition 44 is sufficient for our purpose. The proof of the uniqueness of ${\epsilon}$-gapped ${X_{\infty}}$-homotopy class of ${\mathfrak{g}}$ is beyond the scope of this paper.
\end{Rem}

\section*{Acknowledgements}
This work was carried out during my stay as a research fellow at National Center for Theoretical Sciences (NCTS) and Tokyo Metropolitan University. The author thanks NCTS and Tokyo Metropolitan University for a great research atmosphere and many supports. He also gratefully acknowledges his teacher Kaoru Ono for his continuous support. The author also deeply appreciates the significant contribution of the referees who carefully read the manuscript and provided appropriate advice.

\end{document}